\documentclass[11pt]{article}

\usepackage{amsfonts,amsmath,amssymb,comment,graphicx,hyperref,setspace,todonotes,xcolor}

\usepackage[normalem]{ulem}

\usepackage[margin=1in]{geometry}

\newtheorem{theorem}{Theorem}

\newtheorem{corollary}[theorem]{Corollary}

\newtheorem{example}[theorem]{Example}

\newtheorem{lemma}[theorem]{Lemma}

\newtheorem{proposition}[theorem]{Proposition}
\newtheorem{remark}[theorem]{Remark}

\newtheorem{assumption}[theorem]{Assumption}

\newcommand{\dd}{\,\mathrm{d}}

\newcommand{\stkout}[1]{\ifmmode\text{\sout{\ensuremath{#1}}}\else\sout{#1}\fi}

\DeclareMathOperator*{\esssup}{ess\,sup}

\newenvironment{proof}[1][Proof]{\noindent\textbf{#1.} }{\ \rule{0.5em}{0.5em}}

\newlength{\bibitemsep}
\let\oldthebibliography\thebibliography
\renewcommand\thebibliography[1]{%
  \oldthebibliography{#1}%
  \setlength{\parskip}{\bibitemsep}%
}

\title{\vskip -2.2cm 
Optimal Stopping with Randomly Arriving Stopping Opportunities\thanks{We are very grateful to Albert Menkveld, Peter Spreij, and Florian Wagener for helpful suggestions on an early draft of this paper. 
We are also grateful to conference and seminar participants at
the SIAM FM23 Conference in Philadelphia, 
the Stochastic Analysis and Stochastics of Financial Markets Seminar at Humboldt University Berlin,
the 26th International IME Congress in Edinburgh, 
the Amsterdam--Leuven--London joint PhD workshop 2024 in Leuven, 
the 15th AFMath Conference in Brussels, 
the ICCF24 Conference in Amsterdam, 
DynStoch 2024 in Kiel,
the 3rd Dutch Mathematical Finance Afternoon in Amsterdam, 
Stochastics and Computational Finance 2025 in Lisbon,
the 2026 
World Congress of the Bachelier Finance Society in Bologna,
the University of Amsterdam, 
the Tinbergen Institute and 
the Weierstrass Institute
for their comments and suggestions. 
This research was funded in part by
the Netherlands Organization for Scientific Research (NWO) under grant NWO-Vici 2020--2027 (Laeven)
and by the DFG Excellence Cluster Math+ Berlin, project AA4-2 (Schoenmakers).
}}
\author{Josha A.~Dekker\\
{\footnotesize Dept.~of Quantitative Economics}\\
{\footnotesize University of Amsterdam}\\
{\footnotesize and Tinbergen Institute}\\
{\footnotesize {\tt J.A.Dekker@uva.nl}}\\
\and Roger J.~A.~Laeven\thanks{Corresponding author. University of Amsterdam, 
PO Box 15867, 1001 NJ Amsterdam, The Netherlands.}\\
{\footnotesize Dept.~of Quantitative Economics}\\
{\footnotesize University of Amsterdam, EURANDOM}\\
{\footnotesize and CentER}\\
{\footnotesize {\tt R.J.A.Laeven@uva.nl}}\\[1mm]
\and John G.~M.~Schoenmakers\\
{\footnotesize Stochastic Algorithms \& Nonparametric Statistics}\\
{\footnotesize Weierstrass Institute Berlin}\\
{\footnotesize {\tt schoenma@wias-berlin.de}}\\
\and Michel H.~Vellekoop\\
{\footnotesize Dept.~of Quantitative Economics}\\
{\footnotesize University of Amsterdam}\\
{\footnotesize {\tt M.H.Vellekoop@uva.nl}}}

\vspace{-0.2cm}
\date{\today}

\begin{document}

\maketitle

\vspace{-0.7cm}
\begin{abstract}
We develop simulation-based methods to solve general optimal stopping problems with opportunities to stop that arrive randomly. 
Such problems occur naturally in a wide variety of applications with market frictions, such as limited liquidity, incomplete information and transaction costs. 
Our approach employs random time scales to map the original problem to an integer-time stopping problem with (possibly) infinite horizon. 
We introduce an infinite-horizon policy iteration algorithm to generate lower-biased estimates of the value function and establish a martingale dual representation to generate upper-biased estimates. 
Furthermore, we extend a family of backward dynamic programming methods that includes least-squares Monte Carlo. 
We illustrate the performance of our methods, their general applicability and the managerial implications of randomly arriving opportunities to stop in three examples: a stopped jump-diffusion that can be analyzed in closed form, a random-exercise version of a multi-dimensional max-call contract, and a real options problem in an illiquid market.
\end{abstract}

\textit{Secondary}: 93E24;
91G60.
\newpage

\section{Introduction}
Since the early work of Wald (\cite{W50}) and Snell (\cite{S52}), a vast literature has developed the theory of optimal stopping.
Its areas of application encompass a wide variety of fields, including statistics (sequential testing), operations research (secretary and search problems), and financial (Bermudan and American option pricing) and insurance mathematics (risk process optimization); see, e.g., \cite{PS06} for a textbook treatment.

A standard assumption in this literature is that opportunities to stop arise deterministically. 
For example, for financial options that allow for early exercise timing, i.e.,~exercise before or at a given maturity date, stopping is typically restricted to either pre-specified time instances (Bermudan options) or to any time within a pre-specified time interval (American options). 
More recently, motivated by market frictions, problems in which the stopping opportunities form a random set have been studied. 
Early examples assume the stopping opportunities to be generated by a standard Poisson process, independent of the reward process (\cite{AS01, DW02}). 

This paper develops simulation-based numerical methods to solve general optimal stopping problems with randomly arriving opportunities to stop.
Pivotal to our approach is to proceed on a random time scale, from one randomly arriving opportunity to the next.\footnote{That is, instead of deploying a --- possibly infinitesimally fine --- deterministic time grid, our approach consists in 
considering a random time scale.
On the random time scale, the order of the arrivals of stopping opportunities is encoded in the set of natural numbers 
rather than designated by the actual arrival times on the positive real line.}
This enables us to convert the original optimal stopping problem into a discrete-time optimal stopping problem, with respect to integer-valued stopping times with filtration generated by the consecutively arriving random times, on an infinite horizon.
To encompass the wide variety of problems in which randomly arriving stopping opportunities naturally occur, our setting allows for general dynamics of the stopping opportunities and general dynamics underlying the reward process, including rich forms of stochastic dependence between the stopping opportunities and the reward process. 
Moreover, we consider applications in which the stopping problem is multi-dimensional (easily up to dimension $100$) and in which state transition densities are not available. In such settings, existing numerical solution methods typically become intractable.

We provide both lower and upper bounds to the value function of the stopping problem.
To obtain a lower-biased estimate of the value function, we develop a policy improvement procedure using forward induction.
This non-trivially extends \cite{KS2006}, which was developed in a finite-horizon setting. 
More specifically, since the theory in \cite{KS2006} is entirely based on backward induction from a finite horizon, we extend their main statements to an infinite-horizon setting and give new, essentially different proofs. 
To obtain an upper bound to the value function, we first establish a suitable martingale dual representation, and next employ it in a corresponding Monte Carlo simulation procedure that we design.
We also extend a family of backward dynamic programming methods, including least-squares Monte Carlo (LSMC, \cite{LS01}). 
If the dimension of the underlying process is low, the resulting random times LSMC method is usually accurate and fast. 
It suffers, however, from the curse of dimensionality, in contrast to policy iteration. 
Therefore, for high-dimensional underlying processes, backward dynamic programming serves primarily as an input for our policy improvement approach. 

We illustrate the good performance of our methods in three applications: a geometric Brownian motion with jumps stopped at random times, 
which can be analyzed in explicit closed form; 
a random-exercise version of a multi-dimensional max-call option pricing problem; 
and a real options problem in a scarce, interconnected environment. 
In particular, our methods are shown to deliver (very) small duality gaps, and to be applicable under a wide range of problem specifications allowing for general multivariate dynamics generating the stopping opportunities and reward.
Furthermore, we show that, compared to optimal stopping problems with deterministic stopping opportunities, problems with randomly arriving stopping opportunities feature two additional sources of risk: 
\textit{timing risk}, referring to the change in value due to the randomness in the inter-arrival times between stopping opportunities; 
and \textit{opportunity risk}, referring to the change in value due to the possibility that no further opportunities to stop are seen. 
From a managerial point of view, ignoring either source of risk will lead to incorrect pricing and decision making.
Whereas \textit{opportunity risk} typically warrants a more prudent stopping decision and reduces the value function, we show that \textit{timing risk} may be beneficial or detrimental, depending on the specific application. 


Potential applications of our methods are numerous. 
When \emph{information} is \emph{scarce}, it is natural to consider optimal stopping problems where stopping is restricted to the random times at which an information signal is received; see e.g., \cite{DW02, L12, PY18, PPY21} for applications to American option valuation.
In such a setting, both cross-sectional and temporal dependence between the process generating the signal and the process generating the reward occurs naturally: 
a signal from the release of information may affect underlying prices and hence the reward, which in turn may generate further signals. 
Furthermore, in the presence of \emph{limited liquidity}, or \emph{transaction costs}, stopping opportunities may be linked to the arrival of (attractive) trading opportunities.
Indeed, \cite{RZ02} studies an optimal consumption and investment problem where trades are executed at a Poisson rate instead of continuously, representing, for example, the arrival of buyers. 
Due to the existence of asset-liquidity spirals (see, e.g., \cite{HM14}), dependence between underlying asset prices and liquidity risk, hence the arrival of trading opportunities, arises naturally in this setting. 

In \emph{real option} pricing, constraints of various forms can lead to randomly arriving stopping opportunities. 
Examples are the arrival of scarce goods or investment opportunities (\cite{LRS20, LT21}), or unexpected technological progress (\cite{AS01}). 
In these problems, it is again natural to allow for dependence between the constraints affecting the stopping opportunities and the reward process; we illustrate this in Section~\ref{sec:ex}, in the context of selling a valuable object, employing multivariate point processes allowing for feedback including self- and cross-excitation (\cite{H71,ACL15}) to generate the stopping opportunities and reward.
In the \emph{bank run} literature, \cite{HX12,LLW15} suppose that short-term credit matures randomly, at which point a creditor needs to make a stopping decision: either roll-over debt or run. 
This gives rise to a debt run model based on the classical staggered pricing model of \cite{C83} and the credit risk model by \cite{L98}. 
Other financial and environmental stopping problems also lend themselves naturally to randomly arriving opportunities. 
Examples in \emph{energy economics} include the exploration of natural resources (\cite{PSS88, D04}), weather derivatives (\cite{BLM15, BDL18}) or emissions reduction projects (\cite{FMV15}). 

Optimal stopping problems with constraints may be represented in terms of a Hamilton-Jacobi-Bellman (HJB) equation or a free boundary formulation; see e.g., \cite{MR16}, for a detailed account. 
For general Markov processes, an explicit solution of the continuous-time HJB equation is generally not available; the solution may then be obtained, theoretically and in some cases analytically, through a discrete-time HJB associated to the continuous-time problem.
Stopping problems with randomly arriving opportunities also naturally give rise to a sequence of partial differential equations (PDEs), which may then be solved iteratively.
Indeed, if we consider a sequence of problems in which stopping is only allowed at the first $k$ opportunities (in `problem $k$'), then this problem may be formulated as a PDE where the boundary condition is determined by the solution to `problem $k-1$'. 
This value function iteration approach is taken for infinite-horizon stopping problems in, e.g., \cite{LRS20, LT21}. 
Their `POST' algorithm exhibits geometric convergence of the iterated value functions to the true value function, and \cite{LRS20} indicate how it may be adapted for finite-horizon problems by including time explicitly as a state-variable. 
Value function iteration has also been used for control problems with random control opportunities, both in an infinite-horizon (\cite{PT08, PT09}) and a finite-horizon setting (\cite{FNR13, FNR14}).
The generality of our setting --- in terms of the processes generating the stopping opportunities, the reward processes, the dependences among them, or the dimensionality of our problem --- poses challenges to obtaining numerical solutions using HJB-based methods or value function iteration.

For specific dynamics of the reward and stopping opportunities, explicit solutions may be obtained for stationary (originally) infinite-horizon stopping problems. 
Indeed, for rewards driven by geometric Brownian motion and independent Poissonian arrivals of stopping opportunities (hence, inter-arrival times that are exponential), explicit solutions may be found through variational inequalities and smooth-fit principles (\cite{DW02,GL05}), resolvent methods (\cite{L12}) or Green representations (\cite{AS01, LLW15}). 
Furthermore, recent developments in the fluctuation theory of L\'evy processes may be used to solve stopping problems where rewards follow a L\'evy process, using scale functions and other properties of L\'evy processes; see, e.g., \cite{PY18, PPY21}. 
For a more extensive overview of recent developments in the Poissonian case, see \cite{S23}. 
When arrivals occur according to an inhomogenous Poisson process with state-dependent rate, stochastic coupling and a change of time scale may be used to derive theoretical solutions (\cite{H21}), while if the rate may be controlled, one may obtain a differential equation for the value function that may then be solved (\cite{HZ22}).

For finite-horizon problems, the explicit time-dependence makes it harder to obtain solutions, theoretically and numerically. 
For a problem with state-dependent Cox arrivals, \cite{LLW15} use Green representations to obtain a free-boundary problem for the value function; 
see also \cite{MR16,LRS20} discussed above, and 
\cite{EL09} who treat a related problem where stopping is only allowed while the state process is in a subset of the state space. 

We summarize our main contributions.
Methodologically, we show how general finite-horizon stopping problems with randomly arriving opportunities to stop can be transformed into integer-time, infinite-horizon problems, by mapping deterministic into random time scales, and establish equivalence between their solutions. 
For the induced infinite-horizon problem, we develop and formally substantiate a policy improvement procedure using forward induction; this significantly extends \cite{KS2006} that is entirely based on backward induction from a finite horizon. 
Furthermore, we establish and design a martingale dual representation and a corresponding Monte Carlo procedure.
Besides, we also 
provide a suitable truncation of the horizon that allows a class of backward dynamic programming methods to be extended to our problem, including the least-squares Monte Carlo method in \cite{LS01} and the threshold methods in \cite{A00}. 
Our simulation-based methods are demonstrated to yield good performance in a wide variety of problems and are generally applicable, also in settings in which existing methods become impractical.\footnote{Our simulation-based approach does not require solving HJBs or free boundary problems, sequences of (integro-)PDEs, or deriving explicit expressions for state transition densities; this may be hard (or even impossible) in cases where the reward process is inhomogeneous, the arrival process for opportunities to stop is non-trivial or not independent of the reward process or when the dimensionality of the problem becomes large.} 
From a managerial perspective, we show, across a variety of model specifications, how the randomness of the stopping opportunities impacts both the optimal value and the strategies of the stopping problem, by analyzing their \emph{timing} and \emph{opportunity} risk;
these can be (very) substantial.


This paper is organized as follows. 
In Section~\ref{sec:prel}, we formally introduce the problem and establish equivalence, Bellman and duality results. 
In Section~\ref{sec:PI}, we develop a policy iteration procedure that uses forward induction.
Section~\ref{sec:algs} describes the algorithms that we propose. 
In Section~\ref{sec:numerical}, we present extensive numerical illustrations to assess the accuracy and efficiency of our algorithms, and in Section~\ref{sec:ex}, we analyze a real options example to illustrate the relevance of randomly arriving opportunities for optimal stopping.
Conclusions are in Section~\ref{sec:con}.
Proofs that are not included in the main text, additional technical details and additional numerical results are contained in the Appendix.

\section{Optimal Stopping on Random Times}\label{sec:prel}
We consider a filtered probability space $\left(\Omega,\mathcal{F}, \mathbb{F}, \mathbb{P}\right)$ with filtration $\mathbb{F}:=\left(\mathcal{F}_{t}\right)_{0\leq t\leq\infty}$ satisfying ``the usual conditions'' and suppose that we are given a strictly increasing, random sequence of $\mathbb{F}$-stopping times satisfying
\begin{equation}
    \tau_{0}=0\,, \quad \tau_{i}<\tau_{i+1}<\infty \quad \text{for} \quad i\in \mathbb{N}_{0}\,, \quad \text{and}\quad \lim_{i\rightarrow\infty}\tau_{i}=\infty\,, \quad \text{a.s.} \label{eq:stop}
\end{equation}
Here, $\mathbb{N}$ denotes the (strictly positive) natural numbers and $\mathbb{N}_{0}$ the natural numbers including zero.
On each event $\omega$, we denote the first $K \in \mathbb{N} \cup \{\infty\}$ realized stopping times in \eqref{eq:stop} by $\mathcal{T}_{K}\left(\omega\right) :=\left\{\tau_{i}\left(\omega\right): i\in\mathbb{N}_{0},~i\leq K\right\}$. 

\subsection{Problem definition}
We consider an $\mathbb{F}$-adapted reward process $\left(Z_{t}\right)_{0\leq t<\infty}$ that is càdlàg on $[0,T]$ and satisfies $Z_t = 0$ a.s.~for $t > T$, where $T>0$ is the finite time horizon of our problem. 
Our objective is to optimally stop $Z$ on one of the elements of $\mathcal{T}_K\left(\omega\right)$, i.e., we consider the optimal stopping problem 
\begin{equation}
    Y_{0}^{(K)} := \sup_{\mathbb{F}\text{-}\mathsf{stopping}~\mathsf{time}~\mathfrak{t},~\mathfrak{t}\in\mathcal{T}_{K}}\mathbb{E}\left[Z_{\mathfrak{t}}\right]\,. \label{eq:stopproblem}
\end{equation}
For ease of exposition, we will assume throughout that $Z$ is non-negative. 
We impose the following (standing) integrability condition:
\begin{equation}
    \mathbb{E}\left[\sup_{i\in\mathbb{N}_{0}}Z_{\tau_{i}}^{\alpha}\right] \leq B^{\alpha}\,, \quad \mathrm{for\ some}\ B>0\ \mathrm{and}\ \alpha>1\,. \label{eq:boundZ}
\end{equation}
For $K < \infty$, it suffices to consider the weaker bound $\mathbb{E}\left[\sup_{i \leq K}Z_{\tau_{i}}^{\alpha}\right] \leq B^{\alpha}$. 

Depending on the specific application, the reward process $(Z_{t})$ may represent, e.g., the payoff of exercising a financial contract, 
the market value of a piece of real estate or an art object, 
the value of running when short-term credit matures,
or the payoff of exploiting a natural resources project. 
Furthermore, the first $K$ stopping times $\mathcal{T}_K$ may represent, for example, 
the arrival times of a financial exercise or trading opportunity, 
the release dates of valuable information, 
the arrival times of scarce goods or natural resources investment opportunities, 
or the expiry times of short-term credit. 

\subsection{An equivalent integer-time problem with infinite horizon}\label{sec: discrete-time}
To solve the optimal stopping problem \eqref{eq:stopproblem} on a random sequence of stopping times, we consider the induced discrete-time reward process, $(U_i)_{i\in \mathbb{N}_0}$, and study the corresponding infinite horizon integer-time stopping problem, where
\[ U_{i} := Z_{\tau_{i}}\,, \quad i\in\mathbb{N}_{0}\,. \]

More specifically, due to the order of the stopping times \eqref{eq:stop}, $\mathcal{G}_i := \mathcal{F}_{\tau_i} \subset \mathcal{F}$ defines a discrete filtration $\mathbb{G}:=(\mathcal{G}_{i})_{i\in\mathbb{N}_0}$, where the stop $\sigma$-algebras
\begin{equation}
    \mathcal{F}_{\tau_{i}}:=\left\{A\in\mathcal{F}_{\infty}:A\cap\left\{\tau_{i}\leq t\right\}\in\mathcal{F}_{t}\,,~t \geq 0\right\}\, \label{eq:stopsig}
\end{equation}
satisfy 
\[
\mathcal{F}_{\tau_{0}}\subset\mathcal{F}_{\tau_{1}}\subset\mathcal{F}_{\tau_{2}}\subset\ldots.
\]

As the discrete reward process $(U_i)_{i\in\mathbb{N}_{0}}$ is $\mathbb{G}$-adapted (see, e.g., \cite[Thm.~6]{P90} or \cite{F13}), we may formally relate our original stopping problem \eqref{eq:stopproblem} to a stopping problem over the set of discrete $\mathbb{G}$-stopping times, $\mathcal{N}$, given by
\[ \mathcal{N}:=\left\{\mathfrak{n}:\Omega\rightarrow\mathbb{N}_{0}\text{ such that }\left\{\mathfrak{n}=i\right\} \in \mathcal{G}_{i}\text{ for all }i\in\mathbb{N}_{0}\right\}\,. \]

\begin{proposition}\label{prop:proprandtime} 
    Recall \eqref{eq:stopproblem}. 
    We have that
    \begin{equation}
        Y_{0}^{(K)} = \sup_{\mathfrak{n}\in\mathcal{N},~\mathfrak{n}\leq K}\mathbb{E}\left[U_{\mathfrak{n}}\right]\,. \label{eq:equivstopproblem}
    \end{equation}
\end{proposition}
\begin{proof}
    See Appendix~\ref{app:properties-stop}.
\end{proof}

In view of \eqref{eq:equivstopproblem}, we may thus analyze \eqref{eq:stopproblem} by studying $\sup_{\mathfrak{n}\in\mathcal{N},~\mathfrak{n}\leq K}\mathbb{E}\left[U_{\mathfrak{n}}\right]$ instead. 
Furthermore, for the process $U$, the following bound follows from \eqref{eq:boundZ}:
\begin{equation}
    \mathbb{E}\left[\sup_{i\in\mathbb{N}_{0}}U_{i}^{\alpha}\right] \leq B^{\alpha}\,, \label{eq:boundZ1}
\end{equation}
for some $B>0$ and $\alpha>1$.

\begin{remark}
    The bound \eqref{eq:boundZ}, hence \eqref{eq:boundZ1}, may be determined by simulation. 
    An explicit, analytic bound may be obtained using a slightly stronger condition in conjunction with the Burkholder-Davis-Gundy inequality; see Remark~\ref{rem:analytic-BDG} in Appendix~\ref{app:properties-stop}. 
\end{remark}

Denote by $N_T$ the number of random times until we cross $T$, i.e.,
\begin{equation} 
N_{T} : =\min\left\{m\in\mathbb{N}:\tau_{m}>T\right\}\,. \label{NT}
\end{equation}
It is clear that, for $i \geq N_T\,,~U_i = Z_{\tau_i} = 0$ by our assumptions on $Z$, and in fact
\begin{equation}
    \mathbb{P}(\exists \kappa < \infty: \forall i > \kappa\,,~U_i = 0) = \mathbb{P}(N_T < \infty) = 1\,. \label{eq:stopU}
\end{equation}

In general, the sequence \eqref{eq:stop} may contain one or more deterministic exercise possibilities. 
For example, as a relevant case in applications, the decision-maker may be given the right to stop either at a random time before $T$ or at the terminal time $T$, which means that we need to ensure that \eqref{eq:stop} satisfies a.s.,
\begin{equation}
    \tau_{i\left(\omega\right)}(\omega)=T\text{ for a unique }i\left(\omega\right) \in \mathbb{N}\,. \label{detT}
\end{equation}
One may establish \eqref{detT} by intertwining the deterministic time $T$ with a particular sequence of increasing stopping times, see Lemma~\ref{lemma:insert-T} and Corollary~\ref{cor:insert-T} in Appendix~\ref{app:properties-stop}; these results also establish that \eqref{eq:boundZ} will continue to hold for the new sequence of stopping times if it held for the old sequence of stopping times and $\mathbb{E}[Z_T^\alpha]$ is finite. 

\subsection{Bellman principle and optimal stopping time}
To formulate a Bellman principle for Problem~\eqref{eq:equivstopproblem}, we begin by generalizing it to a stopping problem as seen from a generic point in discrete time $i\in\mathbb{N}_{0}$, which is given in terms of the Snell envelope $Y_i^{(K)}$ and continuation value process $C_i^{(K)}$ for the reward process $(U_i)_{i \in \mathbb{N}_0}$, with $K \in \mathbb{N} \cup \{\infty\}:$
\begin{align}
    Y_{i}^{(K)} &:=\esssup_{\mathfrak{n\in}\mathcal{N},~i\leq\mathfrak{n}\leq K}\mathbb{E}_{\mathcal{G}_{i}}\left[U_{\mathfrak{n}}\right]\,, \quad i\in\mathbb{N}_{0}\,, \label{eq:dcstopproblem} \\
    C_{i}^{(K)} &:=\esssup_{\mathfrak{n\in}\mathcal{N},~i+1 \leq \mathfrak{n}\leq K}\mathbb{E}_{\mathcal{G}_{i}}\left[U_{\mathfrak{n}}\right]\,, \quad i \in \mathbb{N}_{0}\,, \quad 0\leq i<K\,, \label{eq:cstop} \\
    C_{K}^{(K)} &:= 0\,, \quad K < \infty\,. \nonumber
\end{align}
It satisfies the following properties: 
\begin{proposition}\label{prop:bellman} 
Under assumption~\eqref{eq:boundZ}, we have for all $K\in\mathbb{N}\cup\{\infty\}$ and for all $i\in\mathbb{N}_{0}$ with $0\leq i\leq K$ that:
    \begin{itemize}
        \item[(i)] The continuation value process \eqref{eq:cstop} is connected with \eqref{eq:dcstopproblem} via the Bellman principle
              \[ Y_{i}^{(K)}=\max\left(U_{i},C_{i}^{(K)}\right)\,. \]
        \item[(ii)] The process $Y^{(K)}$ is a $\mathcal{G}$-super-martingale, i.e.,
              \begin{equation}
                  Y_{i}^{(K)}\geq\mathbb{E}_{\mathcal{G}_{i}}\left[Y_{i+1}^{(K)} \right]\,, \label{eq:YKsupermart}
              \end{equation}
              and it is the smallest super-martingale that dominates $U$.
        \item[(iii)] The random time $\mathfrak{n}_{i}^{K,\ast}$, defined by
              \begin{equation}
                  \mathfrak{n}_{i}^{K,\ast}:=\inf\left\{i\leq j\leq K:~Y_{j}^{(K)}=U_{j}\right\}\,, \label{iii1}
              \end{equation}
              is an optimal stopping time for $\mathcal{G}$ with $\mathfrak{n}_{i}^{K,\ast}<\infty$ a.s., hence,
              \begin{equation}
                  Y_{i}^{(K)}=\mathbb{E}_{\mathcal{G}_{i}}\left[U_{\mathfrak{n}_{i}^{K,\ast}}\right]\,. \label{iii2}
              \end{equation}
              In particular, for all $K<\infty$,
              \[ Y_{i}^{(K)}\leq\lim_{K^{\prime}\rightarrow\infty}Y_{i}^{(K^{\prime})} = Y_{i}^{\left(\infty\right)}=\esssup_{i\leq\mathfrak{n}<\infty}\mathbb{E}_{\mathcal{G}_{i}}[U_{\mathfrak{n}}]=\mathbb{E}_{\mathcal{G}_{i}}[U_{\mathfrak{n}_{i}^{\infty,\ast}}]\,, \]
              with $\mathfrak{n}_{i}^{\infty,\ast}$ given by \eqref{iii1} for $K=\infty$.
    \end{itemize}
\end{proposition}
\begin{proof}
    See Appendix \ref{app:properties-stop}.
\end{proof}

\begin{corollary}\label{cor:opstoptimes}
    Under assumption \eqref{eq:boundZ}, there exists an optimal stopping time for Problem~\eqref{eq:stopproblem}, which may be obtained by taking  $\mathfrak{t}^{\ast}:=\mathfrak{t}_{0}^{K,\ast}:=\tau_{\mathfrak{n}^{\ast}}:=\tau_{\mathfrak{n}_{0}^{K,\ast}}$.
\end{corollary}
\begin{proof}
    See Appendix~\ref{app:properties-stop}
\end{proof}

\subsection{Duality}
In this subsection, we present a dual representation of the discrete-time optimal stopping problem in our present setting, which essentially goes back to \cite{R02} and \cite{HK04} (see also the early \cite{DK94}). 
We first provide a result for the case $K = \infty$ in Proposition~\ref{prop:proprandtime}, which makes use of a uniform integrability condition that arises naturally for stopping problems on random times. 
The case $K<\infty$ follows analogously, and moreover does not require this extra uniform integrability argument; see Corollary~\ref{corsuff}.

The next proposition tells us that \eqref{eq:equivstopproblem} may be reformulated as a penalized optimization problem with full foresight, in the terminology of perfect information relaxation.

\begin{proposition}\label{dualprop} 
    Let $\mathcal{M}_{0}^{\mathsf{UI}}$ be the collection of uniformly integrable $\mathbb{G}$-martingales $\left(M_{i}\right)_{i\in\mathbb{N}_{0}}$ with $M_{0}=0$. Assume \eqref{eq:boundZ} and
    \begin{equation}
        \mathbb{E}\left[N_{T}\sup_{i\in\mathbb{N}_{0}}U_i\right] <\infty\,. \label{UIcon}
    \end{equation}
    We then have the following duality results. \\
    (i) \textit{Weak duality:} It holds that
    \begin{equation}
        Y_{0}\equiv Y_{0}^{(\infty)}=\inf_{M\in\mathcal{M}_{0}^{\mathsf{UI}}}\mathbb{E}\left[\sup_{i\geq0}\left(U_i-M_{i}\right)\right]\,.\label{dualobj}
    \end{equation}
    (ii) \textit{Strong duality:} Let us consider, for $Y_{i}\equiv Y_{i}^{(\infty)}$ given in \eqref{eq:dcstopproblem}, the $\mathcal{G}$-martingale, a.k.a.~the Doob martingale of the Snell envelope,
    \begin{equation}
        M_{i}^{\circ} := \sum_{j=1}^{i}Y_{j}-\mathbb{E}_{\mathcal{G}_{j-1}}\left[Y_{j}\right]\,, \quad i \geq 1\,, \quad M_{0}^{\circ}=0\,. \label{Mcirc}
    \end{equation}
    Then it holds that $M^{\circ}\in\mathcal{M}_{0}^{\mathsf{UI}}$ and \[Y_{0}=\sup_{i \geq 0}\left(U_{i}-M_{i}^{\circ}\right), \quad\mathrm{a.s.}\]
\end{proposition}

\begin{proof}
    See Appendix~\ref{app:duality}.
\end{proof}

\begin{remark}
Proposition~\ref{dualprop}-(ii) has a beneficial consequence in numerical applications: 
If the Doob martingale $M^\circ$ is approximated `close enough' by some martingale $M$, then the variance of the random variable $\sup_{i \geq 0}\left(U_{i}-M_{i}\right)$ in the dual representation \eqref{dualobj} will be `close' to zero. 
The quantification of this phenomenon is omitted, as it is similar to the finite horizon case (see, for example, \cite{BBS09,LSSS25}).
\end{remark}

\begin{corollary}\label{corsuff}
    If, for some $\alpha > 1$, \eqref{eq:boundZ} holds, a sufficient condition for \eqref{UIcon} is that the survival distribution of $N_T$ is sufficiently regular in the following sense: $\sum_{j=1}^\infty j^{\frac{1}{\alpha-1}}\mathbb{P}(N_T>j) < \infty$. 
    This holds in particular if one of the following holds:
    \begin{itemize}
        \item[(i)] $\mathbb{P}\left(N_{T}>j\right)\lesssim j^{-\gamma}$ for some $\gamma>\alpha/(\alpha-1)$ and $j\geq j_{0}$. 
        \item[(ii)] There is some $L\in\mathbb{N}$ such that $\mathbb{P}\left(\tau_{L}>T\right)=1$; condition \eqref{UIcon} is then trivially fulfilled. 
        However, we are then dealing with Problem~\eqref{eq:stopproblem} for $K=L$ in fact.
    \end{itemize}
    
    Furthermore, if the following condition:
    \[ 0\leq\sup_{i\in\mathbb{N}_{0}}U_i\leq B<\infty\text{ a.s.}\, \]
    is satisfied, condition \eqref{UIcon} holds if $\mathbb{E}\left[N_{T}\right] < \infty$. 
    (In fact, one then has \eqref{corcon} for $\alpha=\infty$.) 
    
    Similarly, if $K < \infty$, it suffices for Proposition~\ref{dualprop} that \eqref{eq:boundZ} holds with $\alpha = 1$, as a finite family of integrable martingales is uniformly integrable. 
\end{corollary}
\begin{proof}
    Follows from Corollary~\ref{cor:holder-UI} and Lemma~\ref{lem:check} in Appendix~\ref{app:duality}.
\end{proof}

\setcounter{equation}{0}
\section{Optimal Stopping via Policy Improvement}\label{sec:PI}

By passing to discrete time in Section~\ref{sec: discrete-time}, our problem becomes an integer-time, infinite-horizon problem. 
In this section, we introduce a method based on forward induction to solve \eqref{eq:equivstopproblem} that extends the iterative policy improvement procedure in \cite{KS2006}; see also \cite{BS06}.
As the results in \cite{KS2006} assume a finite horizon with a deterministic, fixed and finite number of opportunities to stop, and all proofs exploit the backward recursive nature with initialization at the last opportunity to stop, we extend their Theorem~3.1 and Propositions~4.1,~4.3 and~4.4 to the present infinite-horizon setting with essentially different proofs. 

Let us consider a generic stopping family $\sigma$, being a family of stopping times $\sigma := \left(\sigma_{i}\right)_{i\in\mathbb{N}_{0}}$, with respect to the discrete filtration $\mathbb{G}$ 
 that satisfies the consistency conditions
\begin{equation}
    i \leq \sigma_{i}<\infty \text{ and }\sigma_{i} >i\Longrightarrow\sigma_{i}=\sigma_{i+1}\,, \quad i\in\mathbb{N}_{0}\,.\label{con}
\end{equation}
For the stopping family $\sigma$, we define the process
\begin{equation}
    Y_{i}^{\sigma}:=\mathbb{E}_{\mathcal{G}_{i}}\left[U_{\sigma_{i}}\right]. \label{Y-sigma}
\end{equation}

Our aim is to construct, from the consistent stopping policy $\sigma$ in \eqref{con}, an improved policy.
Let us fix a window parameter $\kappa\in\mathbb{N}\cup\left\{\infty\right\}$ and introduce a new stopping family $\widehat{\sigma}$ $:=$ $\left(\widehat{\sigma}_{i}\right)_{i\in\mathbb{N}_{0}}$ 
by
\begin{equation}
    \widehat{\sigma}_{i}:=\inf\left\{j\geq i:~U_{j}\geq \max_{j\leq k\leq j+\kappa}\mathbb{E}_{\mathcal{G}_{j}}\left[U_{\sigma_{k}}\right]\right\}\,. \label{it}
\end{equation}
Clearly, the family $\left(\widehat{\sigma}_{i}\right)$ also satisfies \eqref{con}. 
Moreover, $\left(\widehat{\sigma}_{i}\right)$ improves over $\left(\sigma_{i}\right)$ in the following sense: 
\begin{proposition}\label{polit}
    Let $\left({\sigma}_{i}\right)$ satisfy \eqref{con} and $\left(\widehat{\sigma}_{i}\right)$ be given by \eqref{it}. 
    Then it holds that
    \[ Y_{i}^{\sigma}\leq \widetilde{Y}_{i}^{\sigma} \leq Y_{i}^{\widehat{\sigma}}\leq Y_{i}\,, \]
    with 
    \[ \widetilde{Y}_{j}^{\sigma}:=\max_{j\leq k\leq j+\kappa}\mathbb{E}_{\mathcal{G}_{j}}\left[  U_{\sigma_{k}}\right]\,. \]
\end{proposition}
\begin{proof}
See Appendix~\ref{app:policy-iteration} for details. 
The proof first applies a result from \cite{KS2006} to a suitably truncated reward, and then considers and analyzes a limiting procedure.
\end{proof}

\vskip 0.1cm 
The role of $\kappa$ is to provide the window over which we search for possible improvements, where Part (ii) of Corollary \ref{corimp} indicates that a larger window allows for more improvement:
\begin{corollary}\label{corimp}
    (i) It holds that $Y_{i}^{\widehat{\sigma}}\geq U_{i}$ for all $i\in\mathbb{N}_{0}$. \\
    (ii) Let $\widehat{\sigma}_{i}$ and $\widehat{\sigma}^{\prime}_{i}$ be improvements of $\sigma_{i}$ due to Proposition~\ref{polit} that correspond to $\kappa^{\prime}$ and $\kappa$, respectively. Then, if $\kappa^{\prime}>\kappa$, one has $Y_{i}^{\widehat{\sigma}^{\prime}}\geq Y_{i}^{\widehat{\sigma}}$. 
\end{corollary}
\begin{proof}
    See Appendix~\ref{app:policy-iteration}.
\end{proof}

To solve \eqref{eq:equivstopproblem}, we may construct a sequence of stopping families starting from some initial stopping family $\sigma^{(0)}$ that satisfies \eqref{con}. 
A myopic choice is $\sigma^{(0)}_i=i$, $i\in\mathbb{N}_{0}$, for example.  
The construction is done as follows: if $\sigma^{(m)}$, $m\in\mathbb{N}_{0}$, has already been constructed, then set
\begin{equation}
    \sigma_{i}^{(m+1)}:=\inf\left\{j \geq i:\text{ }U_{j}\geq\max_{j\leq k\leq j+\kappa}\mathbb{E}_{\mathcal{G}_{j}}\left[  U_{\sigma_{k}^{(m)}}\right]\right\}\,. \label{polit1}
\end{equation}
By Proposition~\ref{polit}, each stopping family improves over the previous one, and Theorem~\ref{impthm} below establishes that the sequence converges to an optimal stopping family.
\begin{theorem}\label{impthm} 
    Let $\left(\sigma_{i}^{\ast}\right)$ denote the family of (first) optimal stopping times. 
    That is,
    \[ \sigma_{i}^{\ast}:=\inf\left\{j\geq i:~U_{j}\geq\mathbb{E}_{\mathcal{G}_{j}}\left[Y_{j+1}\right]\right\}\,. \]
    For the sequence $\left(\sigma^{(m)}\right)_{m \in \mathbb{N}_0}$ defined by iterating \eqref{polit1}, we have
    \begin{equation}
        \sigma_{i}^{(m)} \leq \sigma_{i}^{(m+1)} \leq \sigma_{i}^{\ast}<\infty\,, \quad m\in\mathbb{N}_{0}\,, \label{siginc}
    \end{equation}
    and moreover it holds that
    \begin{equation}
        \sigma_{i}^{\ast}=\uparrow\lim_{m\rightarrow\infty}\sigma_{i}^{(m)}\text{ and }Y_{i}=\uparrow\lim_{m\rightarrow\infty} Y_{i}^{(m)},\, \label{pollim}
    \end{equation}
    with $Y_{i}^{(m)}\equiv$ $Y_{i}^{\sigma^{(m)}}$. 
    The up-arrows indicate that the respective sequences are non-decreasing.
\end{theorem}
\begin{proof}
    See Appendix~\ref{app:policy-iteration} for details.
    The proof covers the infinite-horizon case and therefore differs essentially from the corresponding one in \cite{KS2006}. 
    The latter proof was based on backward induction with trivial initialization at the finite horizon. 
\end{proof}

\setcounter{equation}{0}
\section{Simulation-Based Primal and Dual Methods}\label{sec:algs}
In this section, we introduce three simulation-based methods to solve, theoretically and pseudo-algorithmically, problems \eqref{eq:equivstopproblem} (primal) and \eqref{dualobj} (dual) in a Markovian environment specified by the following assumption. 
\begin{assumption}\label{SA}
    There exists a {càdlàg} process $X$ in $\mathbb{R}^{d}$ and a measurable map $(t, x) \mapsto Z(t, x)$ such that $Z_t = Z(t, X_t)$. 
    Furthermore, there exists an auxiliary càdlàg process $\Theta$ in $\mathbb{R}^{q}$ such that the process $\left(X,\Theta\right)$ in $\mathbb{R}^{d}\times\mathbb{R}^{q}$ generates the (augmented) filtration $\mathbb{F}$, and the $\mathbb{G}$-adapted discrete-time process
    \begin{equation*}
     (\tau_{i},X_{\tau_{i}},\Theta_{\tau_{i}})_{i\in\mathbb{N}_{0}}\,, \quad \left(\tau_{0},X_{\tau_{0}},\Theta_{\tau_{0}}\right) = \left(0,X_{0},\Theta_{0}\right)\,,
    \end{equation*}
    is a discrete-time Markov chain in the state space $\mathbb{R}_{\geq0}\times\mathbb{R}^{d}\times\mathbb{R}^{q}$.
\end{assumption}
Assumption~\ref{SA} may be seen as a kind of a strong Markov property,  
and its interpretation is that the joint process $(X,\Theta)$ determines the sequence of stopping times $(\tau_i)$ and that the reward process is given as $Z_{\tau_i}$ $=$  $Z(\tau_i,X_{\tau_i})$.  
In Appendix~\ref{app: markovian}, we establish that this assumption suffices to express conditional expectations 
as Borel functions of the Markov chain and provide some specific examples of possible structures that $(X,\Theta)$ may have. 

In Section~\ref{sec: fha and ls}, we propose two backward dynamic programming approaches to problem  \eqref{eq:equivstopproblem}.
In particular, we propose a least-squares Monte Carlo method in the spirit of \cite{LS01} and a threshold optimization method in the spirit of \cite{A00}, for the case $K < \infty$.
For the case $K = \infty$, we first establish a finite-horizon approximation to  \eqref{eq:equivstopproblem}, for an arbitrary, prespecified level of accuracy, and then proceed with the finite-horizon problem.
Next, in Section~\ref{sec:PIpseudo} we develop a simulation algorithm based on the policy improvement method presented in Section~\ref{sec:PI}
%
for problem~\eqref{eq:equivstopproblem}. 
Finally, for problem~\eqref{dualobj}, we extend the dual method of Andersen-Broadie \cite{AB04} to the infinite-horizon setting in Section~\ref{sec: dual alg}, by exploiting that each path has finite length due to \eqref{eq:stopU}. 

\subsection{Backward dynamic programming, finite-horizon approximation}\label{sec: fha and ls}
For numerically solving problem~\eqref{eq:equivstopproblem} under Assumption~\ref{SA} in the case $K<\infty$, we may take advantage of backward dynamic programming methods, such as the least-squares Monte Carlo algorithm in \cite{LS01}, or some stochastic optimization approach, such as the threshold-based approach in \cite{A00}. 
In this context, one carries out a backward recursion over the stopping opportunities,
trivially initialized at the last opportunity $K$. 
If $K = \infty$, there is no such natural terminal condition, however. 
We may then approximate the infinite-horizon problem by a finite-horizon problem, and next proceed with the latter one.
The following lemma shows that this can be done with arbitrary accuracy. 
\begin{lemma}\label{lem:FHA}
    Assume \eqref{eq:boundZ}.
    For any fixed tolerance level $\varepsilon > 0$, there exists a finite $K$ such that
    \[ Y_0^{(K)} \leq Y_0^{(\infty)} \leq Y_0^{(K)} + \varepsilon\,. \]
\end{lemma}
\begin{proof}
    See Appendix~\ref{app:FHA}. 
\end{proof}

Thus, for any pre-specified tolerance level $\varepsilon$, $Y_{0}^{(K)}$ yields a lower 
bound
 that is \textquotedblleft$\varepsilon$-close\textquotedblright\ to $Y^{(\infty)}_0$.
A suitable candidate for such a $K$, given a fixed tolerance level $\varepsilon > 0$, is the smallest $K'$ such that 
\begin{equation} \mathbb{E}\left[1_{\left\{\tau_{K'}\leq T\right\}} \sup_{i\in\mathbb{N}_{0},~i>K'}Z_{\tau_{i}}\right] < \varepsilon\,, \label{eq:FHA-sim} 
\end{equation}
which may be estimated through simulation.

\subsubsection{Random times least-squares Monte Carlo, pseudo-algorithm}\label{sec:lsmc}
Let us suppose we are in a Markovian environment $\left(  X_{t},\Theta
_{t}\right)  _{t\geq0},$ under Assumption~\ref{SA}, with a reward of the form
$Z_{t}:=Z(t,X_{t})$.
Then, $Y_{0}^{(K)}$, with $K<\infty$, may be numerically approximated by regression, based on a given set of measurable basis functions
\begin{align*}
\psi_{l}:\mathbb{R}_{\geq0}\times\mathbb{R}^{d}\mathbb{\times R}%
^{q}\rightarrow\mathbb{R}\text{, \ \ }l\in\mathbb{N}.
\end{align*}
We propose the following pseudo-algorithm for $n=1,\ldots,N$:

\begin{itemize}
\item[(i)] Simulate
$\left(  X_{\tau_{1}^{(n)}}^{(n)},\Theta_{\tau_{1}^{(n)}}^{(n)}\right)
,\ldots,\left(  X_{\tau_{K}^{(n)}}^{(n)},\Theta_{\tau_{K}^{(n)}}^{(n)}\right)
$.

\item[(ii)] Initialize $\widehat{Y}_{K}^{K,n}=U_{K}^{(n)}=Z(\tau_{K}%
^{(n)},X_{\tau_{K}^{(n)}}^{(n)}),$ $\widehat{\mathfrak{n}}_{K}^{(n)}=K$.

\item[(iii)] If, for $0<k\leq K$, the values $\widehat{Y}_{k}^{K,n}$ and
$\widehat{\mathfrak{n}}_{k}^{(n)}$
have been constructed, then solve the regression problem
\[
\left(  \widehat{\alpha}_{k-1,l}\right)  _{l=1,\ldots,L}:=\arg\min_{\alpha
\in\mathbb{R}^{L}}\sum_{n=1}^{N}\left(  \widehat{Y}_{k}^{K,n}-\sum_{l=1}%
^{L}\alpha_{l}\psi_{l}\left(  \tau_{k-1}^{(n)},X_{\tau_{k-1}^{(n)}}%
^{(n)},\Theta_{\tau_{k-1}^{(n)}}^{(n)}\right)  \right)  ^{2},
\]
and define
\[
\widehat{C}_{k-1}(t,x,\theta):=\sum_{l=1}^{L}\widehat{\alpha}_{k-1,l}\psi
_{l}\left(  t,x,\theta\right)  .
\]
Next, set
\[
\widehat{\mathfrak{n}}_{k-1}^{(n)}=\left\{
\begin{array}
[c]{l}%
k-1\text{ \ \ if \ \ }Z(\tau_{k-1}^{(n)},X_{\tau_{k-1}^{(n)}}^{(n)}%
)>\widehat{C}_{k-1}(\tau_{k-1}^{(n)},X_{\tau_{k-1}^{(n)}}^{(n)},\Theta
_{\tau_{k-1}^{(n)}}^{(n)})\\
\widehat{\mathfrak{n}}_{k}^{(n)}\text{ \ \ \ else.}%
\end{array}
\right.
\]

\end{itemize}

\medskip

Having constructed the functions $\widehat{C}_{0}(t,x,\theta),\ldots
,\widehat{C}_{K-1}(t,x,\theta)$ (with $\widehat{C}_{K}\equiv0$), one thus has
constructed a stopping policy
\begin{equation}
{\sigma}_{k}:=\inf\left\{  k\leq j\leq K:Z(\tau_{j},X_{\tau_{j}}%
)\geq\widehat{C}_{j}(\tau_{j},X_{\tau_{j}},\Theta_{\tau_{j}})\right\}  ,\text{
\ \ }0\leq k\leq K,\label{bdp:LSMC}%
\end{equation}
and we may generate, by a new simulation, i.i.d.\ copies
\[
\left(  \widetilde{\tau}_{1}^{(n)},\widetilde{X}_{\widetilde{\tau}_{1}^{(n)}%
}^{(n)},\widetilde{\Theta}_{\widetilde{\tau}_{1}^{(n)}}^{(n)}\right)
,\ldots,\left(  \widetilde{\tau}_{K}^{(n)},\widetilde{X}_{\widetilde{\tau}%
_{K}^{(n)}}^{(n)},\widetilde{\Theta}_{\widetilde{\tau}_{K}^{(n)}}%
^{(n)}\right)  ,\text{ \ \ for \ \ }n=1,\ldots,\widetilde{N},
\]
of the Markov chain
$
\left(  \tau_{1},X_{\tau_{1}},\Theta_{\tau_{1}}\right)  ,\ldots,\left(
\tau_{K},X_{\tau_{K}},\Theta_{\tau_{K}}\right)  ,
$
in order to obtain a lower-biased estimate $\widetilde{Y}_{0}^{(K)}$ for $Y_{0}^{(K)}$ by setting
\[
\widetilde
{\mathfrak{n}}^{(n)}:=\inf\left\{  0\leq j\leq K:Z(\widetilde{\tau
}_{j}^{(n)},\widetilde{X}_{\widetilde{\tau}_{j}^{(n)}}^{(n)})\geq
\widehat{C}_{j}(\widetilde{\tau}_{j}^{(n)},\widetilde{X}_{\widetilde{\tau}%
_{j}^{(n)}}^{(n)},\widetilde{\Theta}_{\widetilde{\tau}_{j}^{(n)}}%
^{(n)})\right\}  ,
\]
and defining
\[
\widetilde{Y}_{0}^{(K)}:=\frac{1}{\widetilde{N}}\sum_{n=1}^{\widetilde{N}%
}Z(\widetilde{\tau}_{\widetilde{\mathfrak{n}}^{(n)}}^{(n)},\widetilde{X}%
_{\widetilde{\tau}_{\widetilde{\mathfrak{n}}^{(n)}}^{(n)}}^{(n)}).
\]

Theorems~\ref{thm:conv-values} and~\ref{thm:conv-stopping} in Appendix~\ref{app: lsmc convergence} establish convergence to the true value of the optimal stopping problem and the optimal stopping rule, respectively.
The policy~\eqref{bdp:LSMC} can and will be used in Section~\ref{sec: dual alg} as input policy for the well-known Andersen-Broadie algorithm \cite{AB04} (see also \cite{Gl04}) in order to
obtain a dual estimate of an upper bound on $Y_{0}^{(K)}$ (in mean).
\begin{remark}
Over the past decades, various refinements have been proposed for the implementation of the standard, deterministic times LSMC algorithm. 
We have included in our random times LSMC implementation the (counterparts of the)  most obvious and usual refinements: (i) the cash-flow itself is included in the regression; and (ii) the regression is only carried out at in-the-money trajectories. 
\end{remark}

\subsubsection{Threshold optimization, pseudo-algorithm}\label{thresh-pol}

The threshold optimization method proposed by Andersen \cite{A00} can be seen as a special case of stochastic optimization of the stopping region for problem~\eqref{eq:equivstopproblem} with $K<\infty$. 
Formally and adapted to our setting, it aims at constructing an optimal ``threshold policy'' $\left( \sigma_{k}\right) $ due
to thresholds $h_{0},\ldots,h_{K}\geq0$ of the form
\begin{equation}
\sigma_{k}=\inf\left\{  k\leq j\leq K:Z_{\tau_{j}}\geq h_{j}\right\}  ,\text{
\ \ }k=0,\ldots,K,\label{bdp:thresh}%
\end{equation}
which is initialized by $h_{K}=0$ and $\sigma_{K}=K$. 
Then, given
$h_{k+1},\ldots,h_{K}$ and $\sigma_{k+1}$, one defines, recursively,
\begin{align*}
h_{k}  &  :=\arg\max_{h\geq0}\mathbb{E}\left[  Z_{\tau_{k}}1_{\{Z_{\tau_{k}%
}\geq h\}}+Z_{\tau_{\sigma_{k+1}}}1_{\{Z_{\tau_{k}}<h\}}\right]  ,\text{
\ \ and}\\
\sigma_{k}  &  :=k1_{\{Z_{\tau_{k}}\geq h_{k}\}}+\sigma_{k+1}1_{\{Z_{\tau_{k}%
}<h_{k}\}}.
\end{align*}
The Monte Carlo analogue of this method due to sample trajectories
$n=1,\ldots,N$ is obtained as follows. Initialize $h_{K,N}=0$ and
$\sigma_{K,n}=K,$ $n=1,\ldots,N$. Then, given $h_{k+1,N},\ldots,h_{K,N}$ and
$\sigma_{k+1,n}\in\left\{  k+1,\ldots,\bar{K}\right\}  $, $n=1,\ldots,N$, set
\begin{align*}
h_{k,N}  &  :=\arg\max_{h\geq0}\frac{1}{N}\sum_{n=1}^{N}\left(  Z_{\tau_{k}%
}^{(n)}1_{\{Z_{\tau_{k}}^{(n)}\geq h\}}+Z_{\tau_{\sigma_{k+1,n}}}%
^{(n)}1_{\{Z_{\tau_{k}}^{(n)}<h\}}\right)  ,\text{ \ \ and}\\
\sigma_{k,n}  &  :=k1_{\{Z_{\tau_{k}}^{(n)}\geq h_{k,N}\}}+\sigma
_{k+1,n}1_{\{Z_{\tau_{k}}^{(n)}<h_{k,N}\}},\text{ \ \ for \ \ }n=1,\ldots,N.
\end{align*}
It should be noted that in general the optimal threshold policy may be far away from a real optimal policy, but serves in many practical examples as a ``not too bad'' approximation of it.

\subsection{Policy improvement}\label{sec:PIpseudo}

Proposition~\ref{polit} and Theorem~\ref{impthm} suggest a natural algorithm: starting from some initial policy $\sigma$ that satisfies \eqref{con}, apply \eqref{it}
with some window $\kappa \geq 1$, where the conditional expectations are replaced by Monte Carlo estimates. 
The initial policy may be obtained, for example, by the least-squares approach of Section~\ref{sec:lsmc} or as a threshold policy constructed in Section~\ref{thresh-pol}. 
Depending on the underlying structure, other constructions are also possible. 
As the most trivial initial policy, one could take the myopic policy $\sigma_i=i$. 
In the next subsection, we present a corresponding pseudo-algorithm. 

\subsubsection{Pseudo-algorithm}\label{secPI-Alg}
Under Assumption~\ref{SA}, the conditional expectations in \eqref{it} take
the form
\[
\mathbb{E}_{\mathcal{G}_{j}}\left[  U_{\sigma_{k}}\right]  =\mathbb{E}%
_{\left(  \tau_{j},X_{\tau_{j}},\Theta_{\tau_{j}}\right)  }\left[
Z(\tau_{\sigma_{k}},X_{\tau_{\sigma_{k}}})\right]  ,\text{ \ \ }j\leq k\leq
j+\kappa,
\]
and can be estimated in the following way. 
Construct for $n=1,\ldots,N$ the \textquotedblleft outer\textquotedblright\ trajectories $(X_{t}^{(n)}%
,\Theta_{t}^{(n)})_{t\geq0}$, and for trajectory number ${n}$ the sequence
\[
\left(  \tau_{j},X_{\tau_{j}},\Theta_{\tau_{j}}\right)  ^{(n)}\equiv\left(
\tau_{j}^{(n)},X_{\tau_{j}^{(n)}}^{(n)},\Theta_{\tau_{j}^{(n)}}^{(n)}\right)
,\qquad j\geq0.
\]
Then simulate on the outer trajectory with number $n$, for each $m=1,\ldots,M$
and $j\geq0$, the sub-trajectory
\[
\left(  X_{s}^{n,m,j},\Theta_{s}^{n,m,j}\right)  _{s\geq\tau_{j}^{(n)}}\text{
\ \ with \ \ }\left(  X_{\tau_{j}^{(n)}}^{n,m,j},\Theta_{\tau_{j}^{(n)}%
}^{n,m,j}\right)  =\left(  X_{\tau_{j}^{(n)}}^{(n)},\Theta_{\tau_{j}^{(n)}%
}^{(n)}\right)  ,
\]
along with a strictly increasing sequence of $\mathbb{R}_{\geq0}$-valued
stopping times $\left(  \tau_{k}^{n,m,j}\right)  _{k\geq j}$ with $\tau
_{j}^{n,m,j}=\tau_{j}^{(n)}$. 
Next, construct on this sequence the integer-valued stopping times (or indices) $\sigma_{k}^{n,m,j},$ $j\leq k\leq j+\kappa$, according to the input stopping family \eqref{con} and set
\begin{equation}\label{Enjk}
E_{k}^{n,j}\equiv\frac{1}{M}\sum_{m=1}^{M}Z(\tau_{\sigma_{k}^{n,m,j}}%
,X_{\tau_{\sigma_{k}^{n,m,j}}}^{n,m,j})\approx\mathbb{E}_{\left(  \tau
_{j},X_{\tau_{j}},\Theta_{\tau_{j}}\right)  ^{(n)}}\left[  Z(\tau_{\sigma_{k}%
},X_{\tau_{\sigma_{k}}})\right]  ,\text{ \ \ }j\leq k\leq j+\kappa.
\end{equation}
The improved policy \eqref{it} is now path-wise estimated by
\[
\widehat{\sigma}_{0}^{(n)}:=\inf\left\{  j\geq0:Z(\tau_{j}^{(n)},X^{n,j})\geq \max_{{j\leq k\leq j+\kappa}}E_{k}^{n,j}\right\},
\]
and we finally obtain the estimate
\[
Y_{0}^{\sigma}\leq Y_{0}^{\widehat{\sigma}}=\mathbb{E}\left[
U_{\widehat{\sigma}_{0}}\right]  \approx\frac{1}{N}\sum_{n=1}^{N}%
Z(\widehat{\sigma}_{0}^{(n)},X_{\widehat{\sigma}_{0}^{(n)}}^{(n)}),
\]
which is naturally a lower-biased estimate in the sense that $Y_{0}^{\widehat{\sigma}}$ $\leq$ $Y_0$.

\begin{remark}\label{NumSub}
In the above pseudo-algorithm, the number of sub-simulations is taken to be equal for any $j\geq0$, for simplicity. 
Of course, the algorithm may be trivially refined in this respect, to further improve its efficiency.     
\end{remark}

\begin{remark}
In the implementation of the policy improvement algorithm, we have used the estimate due to the initial policy employing a form of variance reduction as proposed in Remark~6.2 in \cite{KS2006}.
\end{remark}

\subsection{Dual estimation}\label{sec: dual alg}
To obtain an upper-biased estimate of $Y_0$ from the dual representation~\eqref{dualobj}, we proceed by adapting the algorithm from \cite{AB04} to an infinite-horizon setting. 
For a given stopping family $\sigma=\left(  \sigma_{i}\right)  _{i\in
\mathbb{N}_{0}}$ satisfying \eqref{con}, we write the Doob martingale corresponding to \eqref{Y-sigma} as
\begin{equation}
M_{i}^{\sigma}:=\sum_{j=1}^{i}Y_{j}^{\sigma}-\mathbb{E}_{\mathcal{G}_{j-1}%
}\left[  Y_{j}^{\sigma}\right]  \,,\quad i\geq1\,,\quad M_{0}^{\sigma
}=0\,\label{doob-mart},%
\end{equation}
which satisfies the following recursion:
\begin{equation}
M_{j}^{\sigma}-M_{j-1}^{\sigma}=Y_{j}^{\sigma}-\mathbb{E}_{\mathcal{G}_{j-1}%
}\left[  Y_{j}^{\sigma}\right]  =Y_{j}^{\sigma}-Y_{j-1}^{\sigma}-1_{\left\{
\sigma_{j-1}=j-1\right\}  }\left(  \mathbb{E}_{\mathcal{G}_{j-1}}\left[
Y_{j}^{\sigma}\right]  -Y_{j-1}^{\sigma}\right)  \,,\label{doob-sigma-rec}%
\end{equation}
by the consistency of the family $\sigma$. 
We thus obtain $M_{i}^{\sigma}$ by summing \eqref{doob-sigma-rec} from $j=1$ to $j=i$, and then define
\begin{align}
\delta_{i}^{\sigma} &  :=U_{i}-M_{i}^{\sigma}-Y_{0}^{\sigma}=U_{i}%
-Y_{i}^{\sigma}+\sum_{j=1}^{i}1_{\left\{  \sigma_{j-1}=j-1\right\}  }\left(
\mathbb{E}_{\mathcal{G}_{j-1}}\left[  Y_{j}^{\sigma}\right]  -Y_{j-1}^{\sigma
}\right)  \nonumber\\
&  \ =1_{\left\{  \sigma_{i}>i\right\}  }\left(  U_{i}-\mathbb{E}%
_{\mathcal{G}_{i}}\left[  U_{\sigma_{i+1}}\right]  \right)  +\sum_{j=0}%
^{i-1}1_{\left\{  \sigma_{j}=j\right\}  }\left(  \mathbb{E}_{\mathcal{G}_{j}%
}\left[  U_{\sigma_{j+1}}\right]  -U_{j}\right).  \label{del-sig}
\end{align}
Now, a dual upper bound due to the policy $\sigma$ is obtained by Proposition~\ref{dualprop}(i), as
\begin{align}
Y_{0}^{\sigma,up} &  :=\mathbb{E}\left[  \sup_{i\geq0}(U_{i}-M_{i}^{\sigma
})\right]  =Y_{0}^{\sigma}+\mathbb{E}\left[  \sup_{i\geq0}\delta_{i}^{\sigma
}\right]  \nonumber\\
&  =:Y_{0}^{\sigma}+\Delta_{0}^{\sigma},\label{dual-sigma}%
\end{align}
provided that $M^{\sigma}$ is uniformly integrable. 
Note that
\begin{equation}
\delta_{i+1}^{\sigma}-\delta_{i}^{\sigma}=1_{\left\{  \sigma_{i+1}%
>i+1\right\}  }\left(  U_{i+1}-\mathbb{E}_{\mathcal{G}_{i+1}}\left[
U_{\sigma_{i+2}}\right]  \right)  +\mathbb{E}_{\mathcal{G}_{i}}\left[
U_{\sigma_{i+1}}\right]  -U_{i},\text{ \ \ }i\geq0, \label{deli-incr}
\end{equation}
which allows for a path-wise upward recursive computation of the $\delta
_{i}^{\sigma}$ in the pseudo-algorithm presented below.

\subsubsection{Pseudo-algorithm}\label{sec:dualpa}
Due to Assumption~\ref{SA}, we may estimate \eqref{dual-sigma} by using nested Monte Carlo in the spirit of Section~\ref{secPI-Alg}. 
In view of the monotonicity of the (conditional) expectation, it is easily verified that,
upon replacing the expectation and conditional expectations in \eqref{dual-sigma}, and hence in \eqref{del-sig}, by nested Monte Carlo simulation
due to outer trajectories and sub-trajectories, respectively, one obtains a ``biased high'' estimator. 
That is, one obtains an estimator with expectation larger than $Y_{0}^{\sigma,up}$. 
In pseudo-algorithmic terms, we proceed in the following way.
\begin{itemize}
\item 
     For $n=1$ to $N$:
    \begin{itemize}
        \item 
        Simulate trajectory number $n$ with rewards  $(U_{i}^{(n)})_{i=0,\ldots,N_{T}^{(n)}}$,  
      where $U_{i}^{(n)}=Z(\tau_i^{(n)},X^{(n)}_{\tau_i^{(n)}})$ and $N_{T}$ is path-wise defined in (\ref{NT});
       \item 
       Following the construction of (\ref{Enjk}) in the pseudo-algorithm of 
       Section~\ref{secPI-Alg}, estimate, for all $i,$ $0\leq
       i\leq N_{T}$, the conditional expectations in \eqref{del-sig} by
        \[
       E_{i+1}^{n,i}\approx\mathbb{E}_{\mathcal{G}_{i}}[U_{\sigma_{i+1}^{(n)}}^{(n)}],
       \]
       due to $M$ sub-simulations; 
      \item
       Construct initially $\widehat{\delta}_{0}^{\sigma,n}:=
       1_{\left\{  \sigma_{0}^{(n)}>0\right\}  }\left(  U_{0}%
      -E_{1}^{n,0}\right)  \approx\delta_{0}^{\sigma,n}.$ 
      Next, after having constructed $\widehat{\delta}_{i}^{\sigma,n}$ for $0\leq i<N_{T}^{(n)},$ 
      construct 
      \[
       \widehat{\delta}_{i+1}^{\sigma,n}:=\widehat{\delta}_{i}^{\sigma,n}+1_{\left\{
       \sigma_{i+1}^{(n)}>i+1\right\}  }\left(  U_{i+1}^{(n)}-E_{i+2}^{n,i+1}\right)
     +E_{i+1}^{n,i}-U_{i}^{(n)}\approx\delta_{i+1}^{\sigma,n}%
      \]
      in view of \eqref{deli-incr};   
     \end{itemize}
\item
    Estimate the duality gap in \eqref{dual-sigma} by
    \begin{equation}
    \widehat{\Delta}_{0}^{\sigma}:=\frac{1}{N}\sum_{n=1}^{N}\max_{0\leq i\leq
     N_{T}^{(n)}}\widehat{\delta}_{i}^{\sigma,n}\,.\label{gap-mc}%
\end{equation}
\end{itemize}
We now obtain as an upper-biased estimate for the dual upper bound,
\begin{equation}
\widehat{Y}_{0}^{\sigma,up}=\widehat{Y}_{0}^{\sigma}+\widehat{\Delta}%
_{0}^{\sigma}\,,\label{dual-mc}%
\end{equation}
where the lower-biased estimate $\widehat{Y}_{0}^{\sigma}$ of $Y_{0}^{\sigma}$ may be obtained by an unbiased non-nested Monte Carlo estimation due to a sufficiently large number of trajectories. 
We finally note that, concerning the number of sub-simulations involved, a remark similar to Remark~\ref{NumSub} also applies here.

\setcounter{equation}{0}

\section{Numerical Illustrations}\label{sec:numerical}
To illustrate the performance and general applicability of the numerical methods introduced in Section~\ref{sec:algs}, we begin by considering an example that admits an analytical solution in Section~\ref{ssec:explicit}.
Next, we extend in Section~\ref{ssec:maxcall} the max-call example from \cite{AB04} to a random exercise opportunities setting, allowing for various dependence structures between the process generating the stopping opportunities and the process generating the reward. 

\subsection{Stopped geometric Brownian motion with jumps}\label{ssec:explicit}
To analyze the performance of our algorithms for a problem admitting an explicit, analytical solution, which can serve as a benchmark, we consider a geometric Brownian motion with jumps: 
\begin{equation}
    X_t = X_0 e^{\left(\mu - \frac{1}{2}\sigma^2\right)t + \sigma W_t + N_t \ln(1+j)}\,, \label{eq:GBM-bench}
\end{equation}
where $X_0 > 0,~\mu\in\mathbb{R},~\sigma > 0$ and $j>-1$ are known constants. 
Here, $W$ is a standard Brownian motion and $N$ is a homogeneous Poisson process with intensity $\lambda>0$, assumed to be mutually independent. 
We consider the problem of stopping the reward
\begin{equation}
    Z_t = X_t^\eta 1_{\{t \leq T\}} ,\label{eq:Z-bench}
\end{equation}
for a given power preference parameter $\eta>0$ and time horizon $T>0$, on
$\{\tau_{i},\ldots,\tau_{i+K}\}$ with $i\in\mathbb{N}_{0}$, 
the set of $\tau_{i}$ and the first $K\geq 1$ jump times of $N$ that occur strictly after $i$. 
We wish to determine, for $i\in\mathbb{N}_{0}$, 
\begin{align*}
   Y_{i}^{(K)}   
 &:= \esssup_{i\leq{n}\leq i+K}\mathbb{E}_{\mathcal{G}_{i}}\left[Z_{\tau_n}\right]\,,
 \qquad
  C^{(K)}_i:=
 \esssup_{i+1\leq{n}\leq i+K}\mathbb{E}_{\mathcal{G}_{i}}\left[Z_{\tau_n}\right]\,,
\end{align*}
with $\tau_0:=0$. 
These definitions are similar to the ones in \eqref{eq:dcstopproblem}--\eqref{eq:cstop} but for ease of exposition we specify in this subsection $K$ remaining stopping opportunities instead of $K$ opportunities in total.
Letting  
\begin{align}
    \xi := \mu \eta + \frac{1}{2}\sigma^2\eta(\eta-1) - \lambda\,, \label{eq:bench-alpha} 
    \qquad\psi := \lambda(1+j)^\eta\,, 
\end{align}
a straightforward application of It\^o's formula shows that if $\xi + \psi \leq 0$, $(Z_{\tau_k})_{k\in\mathbb{N}}$ is a $\mathbb{G}$-super-martingale. 
Since it equals its own Snell envelope in this case, it is then always optimal to stop at the first opportunity.
Note that we do not need to impose that $Z$ is an $\mathbb{F}$-super-martingale at all times to conclude that we should stop as soon as possible, only that $(Z_{\tau_k})_{k\in\mathbb{N}}$ is a $\mathbb{G}$-super-martingale. 
If that is the case, the following is immediate:
\begin{lemma}
If $\xi + \psi \leq 0$, then  $\xi<0$ and, for any $K \in \mathbb{N} \cup \{\infty\}$, it is optimal to stop immediately, so $Y_i^{(K)}=Z_{\tau_i}$. 
Furthermore, we have that the continuation value $C_i^{(K)}$ equals the expected value of $Z$ when stopped at the first opportunity, so
\begin{equation*} 
    C_i^{(K)} = \mathbb{E}_{\mathcal{G}_{i}}\left[Z_{\tau_{i+1}}\right] = \int_0^{T-\tau_i}\psi e^{\xi s}X_t^\eta \,\mathrm{d}s = 
    X_t^\eta\frac{\psi}{\xi}\left(e^{\xi (T-\tau_i)} - 1 \right) \,. 
\end{equation*}
\end{lemma}
By induction on $K$ and a limiting argument, we may also find a solution for the case when $\xi > 0$:
\begin{theorem}\label{thm:bench-sub}
    Assume $\xi > 0$. 
    For finite $K \geq 1$ and $K=\infty$, it is optimal to stop at time $\tau_i\leq T$ if and only if $$\tau_i \ \geq\  t^*:=T-\frac{\ln(1+\frac\xi\psi)}{\xi}.$$ 
    The corresponding optimal values are $Y_i^{(K)}=X_{\tau_i}^\eta\vee  C_i^{(K)}$ with:
\begin{align*}
        C_{i}^{(1)} &= X_{\tau_i}^\eta \frac{\psi}{\xi}\left(e^{\xi (T-\tau_i)} - 1\right)\,, \\
        C_{i}^{(K)} 
        &=1_{\{\tau_i \geq  t^*\}} C_{i}^{(1)}\, +
        1_{\{\tau_i < t^*\}}\, X_{\tau_i}^\eta e^{\xi (t^* - \tau_i)} z_{K}(t^* - \tau_i)\,, \quad 1 < K < \infty ,\\
        C_{i}^{(\infty)} 
         &= 1_{\{\tau_i \geq  t^*\}} C_{i}^{(1)}\, +
        1_{\{\tau_i < t^*\}}\, X_{\tau_i}^\eta e^{(\xi+\psi) (t^* - \tau_i)}\,,
\end{align*}
    where 
    \[z_{K}(s) = \frac{\Gamma_{K}(\psi s)}{(K-1)!}e^{\psi s} + \frac{\gamma_{K}(-\xi s)}{(K-1)!}\left(-\frac{\psi}{\xi}\right)^{K}e^{-\xi s}\,,\]
    with $\Gamma_{K}$ and $\gamma_{K}$ the upper and lower Gamma functions of order $K$, respectively.
\end{theorem}
\begin{proof}
    See Appendix~\ref{app:benchmark}.
\end{proof}
\newline

\begin{figure}[t]
    \centering
    \includegraphics[width=0.8\textwidth]{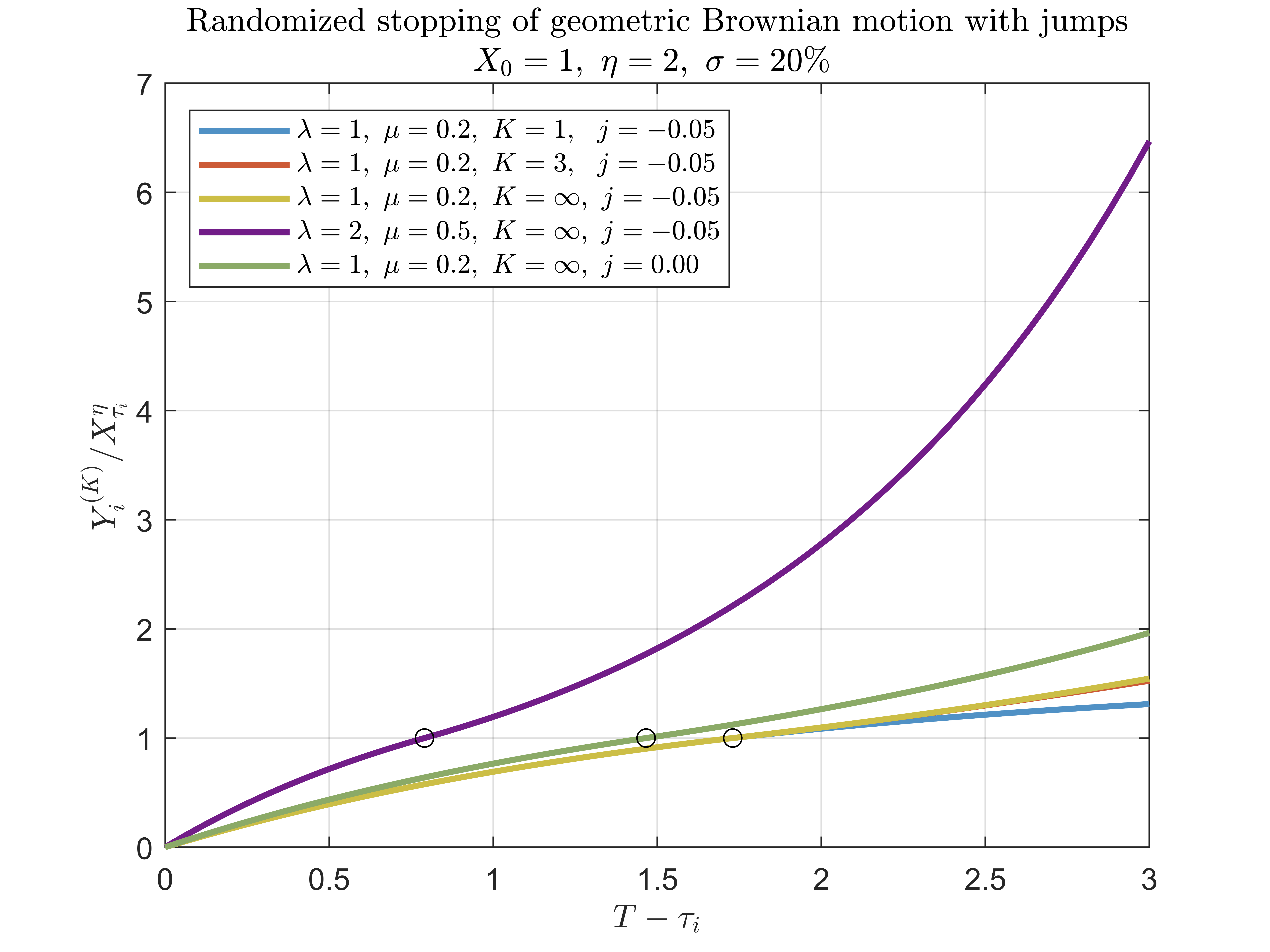}
    \caption{Randomized stopping of geometric Brownian motions with jumps. \\ 
    {\footnotesize\textit{Notes.} 
    The scaled value function $Y_i^{(K)}/X_{\tau_i}^\eta$ is displayed for different remaining times to maturity $T-\tau_i$. 
    The different lines correspond to the cases shown in Table~\ref{tab:benchmark}; see there for a more detailed description of the model and assumptions. 
    The black circles indicate the remaining time to maturities $T-t^*$ where it becomes optimal to stop, i.e., the times where $Y_i^{(K)} = X_{\tau_i}^\eta$. 
    We notice that the second and third cases (the red and yellow lines) are almost indistinguishable.}
    }
    \label{fig:gbmpoisson}
\end{figure}

In Figure~\ref{fig:gbmpoisson}, we visualize the scaled solutions of Theorem~\ref{thm:bench-sub} for different characteristics of the arrival process inducing the exercise opportunities   
and different parameters of the geometric Brownian motion with jumps.  
We notice that for the parameter values used here, the effect of allowing for at most $1$, $3$ or an unlimited number of randomly arriving remaining exercise opportunities is much smaller than the effect of changing the size of the downward jumps or the expected number of exercise opportunities before maturity. 
The optimal stopping times in terms of the remaining time to maturity have been indicated by circles in Figure~\ref{fig:gbmpoisson}. 
They correspond to the point in time at which the scaled value function equals one. 
As expected, stopping tends to take place later (hence, $T-t^*$ is smaller) when the expected number of remaining exercise opportunities is larger and/or when future values of the geometric Brownian motion with jumps have a higher expected value. 

The explicit form of the value function makes this optimal stopping problem a suitable benchmark to test our random times primal-dual approach.
In Table~\ref{tab:benchmark}, we display the corresponding estimates provided by the algorithms from Sections~\ref{sec:lsmc} and~\ref{sec:dualpa}. 
The estimates from the random times LSMC and dual algorithms are close to the exact value functions from Theorem~\ref{thm:bench-sub} across various specifications, and duality gaps are of the order 0.1\%. 
As the least-squares method requires a finite horizon, we have used the finite-horizon approximation due to Lemma~\ref{lem:FHA} for all problems with $K = \infty$ and for all methods, to ensure a fair comparison. 
Specifically, we have determined a suitable $K'$ by simulating \eqref{eq:FHA-sim} with $100\mathord{,}000$ simulated paths and $\varepsilon = 0.001$. 
All simulated paths have been generated in antithetic pairs.

\begin{table}[h!]
    \centering

{\small\begin{tabular}{rrrrr||rrrrr}
$j$ & $K$ & $\mu$ & $\lambda$ & $T$ & True & LSMC & (s.e.) & Dual & (s.e.) \\\hline\hline 
-0.05 & 1 & 0.2 & 1 & 3 & 1.31125 & 1.31103 & (0.00066) & 1.31103 & (0.00000) \\
-0.05 & 3 & 0.2 & 1 & 3 & 1.52503 & 1.52484 & (0.00096) & 1.52507 & (0.00006) \\
-0.05 & $\infty$ & 0.2 & 1 & 3 & 1.54483 & 1.54443 & (0.00099) & 1.54498 & (0.00012) \\
-0.05 & $\infty$ & 0.2 & 1 & 1 & 0.69104 & 0.69101 & (0.00042) & 0.69101 & (0.00000) \\
-0.05 & $\infty$ & 0.5 & 2 & 3 & 6.46847 & 6.47332 & (0.00482) & 6.47560 & (0.00041) \\
0 & $\infty$ & 0.2 & 1 & 3 & 1.96392 & 1.96316 & (0.00135) & 1.96399 & (0.00024) \\
\end{tabular}} 
    \caption{
    Randomized stopping of geometric Brownian motion with jumps.\\ 
    {\footnotesize\textit{Notes.} 
    This table displays estimated values of $Y_0^{(K)}$ under homogeneous Poisson opportunities with rate $\lambda$. 
    $X_t$ follows \eqref{eq:GBM-bench} with $X_0 = 1$, $\sigma = 20\%$ and various specifications for the remaining model parameters, whereas $Z_t$ follows \eqref{eq:Z-bench} with $\eta = 2$. 
    Column \rm{True} contains the true values (from Theorem~\ref{thm:bench-sub}). 
    Column \rm{LSMC} contains the lower-biased estimates based on the random times LSMC algorithm of Section~\ref{sec: fha and ls} with $200\mathord{,}000$ and $2\mathord{,}000\mathord{,}000$ independent paths for the backward and forward passes, respectively, along with its standard error. 
    The continuation values are estimated by using the basis functions $(t,x) \mapsto t^ix^j$, $i+j\leq 3$. 
    Column \rm{Dual} contains the upper-biased estimates based on the dual algorithm in Section~\ref{sec: dual alg} using $1\mathord{,}500$ paths and Doob martingales \eqref{doob-mart} based on the estimated LSMC policy, where conditional expectations have been estimated using $10\mathord{,}000$ sub-simulations, along with its standard error. 
    All simulated paths have been generated in antithetic pairs.}}
    \label{tab:benchmark}
\end{table}

\subsection{Max-call option}\label{ssec:maxcall}
In view of the good performance of our random times LSMC primal-dual algorithm in the benchmark example, we now analyze our algorithms in a more complex example for which analytic reference values are no longer available.
To illustrate the flexibility of our set-up, we consider several extensions of the max-call example in \cite{AB04}. 
The objective is to optimally exercise a max-call contract on $d$ assets with maturity date $T $, for which the payoff if exercised at time $t$ is given by
\begin{equation}
    Z_t = e^{-rt}\left(\max_{j=1,\ldots,d}X_t^j - \mathcal{K}\right)^+1_{\{t \leq T\}}\,. \label{eq:maxcall}
\end{equation}
Here, the $X_t^j$ are jump-diffusions, interpreted as underlying asset prices driving the value of potential investment projects, with dynamics governed by the following stochastic differential equations:
\begin{equation*}
    \dd{X}_t^j = (r - \delta^j - \mu_J^j\lambda_J^j)X_t^j\dd{t} + \sigma^jX_t^j\dd{W}_t^j + \mu_J^jX_{t-}^j\dd{N}_t^j\,, \quad j = 1,\ldots,d\,, 
\end{equation*}
where the $W_t^j$ are independent standard Brownian motions, the $N_t^j$ are independent homogeneous Poisson processes with arrival rates $\lambda_J^j \geq 0$, with the $W_t^j$ and the $N_t^j$ assumed to be mutually independent, and $r,~\delta^j,~\mu_J^j \in (-1,\infty)$ and $\sigma^j > 0$ are deterministic constants. 

\subsubsection{Reward-independent Poissonian arrivals of stopping opportunities}
We begin by considering the stopping problem \eqref{eq:stopproblem} for the reward \eqref{eq:maxcall} with stopping opportunities generated according to a reward-independent homogeneous Poisson process with rate $\lambda$:
\begin{equation}
    \tau_{i+1} = \tau_i + \eta_i\,, \quad \eta_i \sim \mathrm{Exp}(\lambda)\,~\mathrm{i.i.d}, \quad\tau_0 = 0\,. \label{eq:poisson}
\end{equation}
Hence, this example can be seen as (a multi-dimensional, jump-diffusion extension of) the finite-horizon counterpart of the perpetual call considered e.g., in \cite{DW02}.

Table~\ref{tab: LSMC poisson 5} displays value function estimates obtained by using the random times LSMC method of Section~\ref{sec: fha and ls}, the dual method of Section~\ref{sec: dual alg} and the policy iteration (PI) method of Section~\ref{sec:PIpseudo}, respectively, for the case of $d = 5$ assets with $T = 3,~r = 5\%,~\delta^j = 10\%$, $\sigma^j = 20\%$ and $\lambda_J^j = \mu_J^j = 0$. 
Note that the duality gaps, i.e., the differences between random times LSMC and dual estimates, are typically (very) small (less than 0.4\% of the option value), and are generally slightly increasing with the arrival rate $\lambda$ of stopping opportunities. 
The fact that random times LSMC already performs very well, leaving little room for further improvement, motivates our use of a trivial initial policy for the policy iteration method in this table: if we were to consider policy iteration over random times LSMC initial policies, standard errors would dominate the improvements of the policy iteration method. 
The trivial policy consists in stopping as soon as $Z_{\tau_i} > 0$ or $\tau_i \geq T$. 

We consider the case with jumps ($\lambda_J^j = 1,~\mu_J^j = 0.06$) in Table~\ref{tab: LSMC jumps poisson 5}. 
For comparison, we also provide Tables~\ref{tab: LSMC poisson 1}--\ref{tab: LSMC poisson 2} and~\ref{tab: LSMC jumps poisson 1}--\ref{tab: LSMC jumps poisson 2} for the cases with $d = 1$ and $d = 2$ assets in Appendix~\ref{app:tables}. 
The duality gaps and standard errors are (even) smaller there, which is to be expected when the dimension of the problem is smaller.
The results in Tables~\ref{tab: LSMC poisson 5} and \ref{tab: LSMC poisson 2} may be viewed as the random times counterpart of Table~2 in \cite{AB04}, which considers a Bermudan max-call option on the same underlying assets.
More specifically, let us consider the case $d=2,5$ in Tables~\ref{tab: LSMC poisson 5} and~\ref{tab: LSMC poisson 2} and compare our results to those in Table~2 of \cite{AB04}.
They consider 9 equidistant exercise opportunities on the interval $[0,3]$, whereas in our setting there are in expectation $\lambda T = \{3,9,15\}$ times at which the contract may be exercised. 
It is remarkable that, even for the case $\lambda = 5$, i.e., when there are $15$ exercise opportunities in expectation, our estimated values are substantially lower than those in \cite{AB04}: 
the difference is about $0.4\$$ ($d=2$) and $1.0\$$ ($d=5$), which amounts to $3\%$ to $4\%$ of the asset price. 
For $\lambda = 1,3$, this difference is even (much) larger. 
This entails that there is a substantial premium on ``market frictions''.

To analyze these effects in more detail, we compute the probability of not exercising, $\mathbb{P}(\tau^* > T)$, the expected time to maturity at exercise, $\mathbb{E}[T-\tau^* | \tau^* \leq T]$, and the \textit{opportunity risk} and \textit{timing risk} mentioned in the Introduction, and formally defined as
\begin{align*}
    R_O &:= \sup_{\mathbb{F}\text{-}\mathsf{stopping}~\mathsf{time}~\mathfrak{t},~\mathfrak{t}\in\mathcal{T}_{K}}\mathbb{E}\left[Z_{\mathfrak{t}}\right] - \sup_{\mathbb{F}\text{-}\mathsf{stopping}~\mathsf{time}~\mathfrak{t},~\mathfrak{t}\in\mathcal{T}_{K} \cup \{T\}}\mathbb{E}\left[Z_{\mathfrak{t}}\right]\,, \\ 
    R_T &:= \sup_{\mathbb{F}\text{-}\mathsf{stopping}~\mathsf{time}~\mathfrak{t},~\mathfrak{t}\in\mathcal{T}_{K} \cup \{T\}}\mathbb{E}\left[Z_{\mathfrak{t}}\right] - \sup_{\mathbb{F}\text{-}\mathsf{stopping}~\mathsf{time}~\mathfrak{t},~\mathfrak{t}\in t_1,\ldots,t_{\tilde{K}}}\mathbb{E}\left[Z_{\mathfrak{t}}\right]\,,\quad {\tilde{K}} = \left\lceil\mathbb{E}\left[\sum_{\tau_i}1_{\tau_i \leq T}\right]\right\rceil\,,
\end{align*}
where the problems $\sup_{\mathbb{F}\text{-}\mathsf{stopping}~\mathsf{time}~\mathfrak{t},~\mathfrak{t}\in\mathcal{T}_{K} \cup \{T\}}\mathbb{E}\left[Z_{\mathfrak{t}}\right]$ and $\sup_{\mathbb{F}\text{-}\mathsf{stopping}~\mathsf{time}~\mathfrak{t},~\mathfrak{t}\in t_1,\ldots,t_{\tilde{K}}}\mathbb{E}\left[Z_{\mathfrak{t}}\right]$ may be computed by the (deterministic-times) counterparts of the algorithms used. 
We pick $t_i = \frac{i}{\tilde{K}}T$, such that the last term in $R_T$ corresponds to a Bermudan stopping problem with equidistant stopping opportunities. 

Table~\ref{tab: LSMC poisson 5 risks} illustrates how the timing risk and the opportunity risk are sizeable for the max-call example, and vary with arrival rates and initial conditions. 
We note that for this problem, both opportunity and timing risk increase with $\lambda$, whereas the expected time to maturity at exercise and the probability of stopping after $T$ decrease with $\lambda$. 

\begin{table}[ht!]
    \centering

{\small \begin{tabular}{rr||rrrrrrrr}
$\lambda$ & $X_0$ & LSMC & (s.e.) & Dual & (s.e.) & Trivial & (s.e.) & Trivial-PI & (s.e.) \\\hline\hline 
1 & 90  & 11.67443 & (0.01032) & 11.68861 & (0.00307) & 9.16316 & (0.00838) & 11.65284 & (0.03108) \\
  & 100  & 19.53414 & (0.01284) & 19.55874 & (0.00417) & 14.74861 & (0.01012) & 19.50832 & (0.04418) \\
  & 110  & 28.64464 & (0.01503) & 28.69280 & (0.00669) & 24.45769 & (0.01185) & 28.59180 & (0.04894) \\
3 & 90  & 15.01530 & (0.01087) & 15.06641 & (0.01824) & 7.28141 & (0.00574) & 14.73748 & (0.04027) \\
  & 100  & 24.15326 & (0.01308) & 24.22701 & (0.02013) & 10.88004 & (0.00670) & 23.76937 & (0.05106) \\
  & 110  & 34.44991 & (0.01486) & 34.55357 & (0.02302) & 21.12589 & (0.00757) & 33.96516 & (0.05650) \\
5 & 90  & 15.76541 & (0.01098) & 15.81797 & (0.00450) & 5.87923 & (0.00451) & 15.23885 & (0.03972) \\
  & 100  & 25.09820 & (0.01317) & 25.19821 & (0.00827) & 8.74741 & (0.00528) & 24.30068 & (0.04866) \\
  & 110  & 35.56730 & (0.01488) & 35.70128 & (0.00905) & 19.02219 & (0.00595) & 34.61747 & (0.05376) \\
\end{tabular}
}
    \caption{Value of a max-call option on random times ($d=5$). \\ 
    {\footnotesize\textit{Notes.} 
    This table displays estimated values of $Y_0^{(\infty)}$ in the optimal stopping problem for payoff \eqref{eq:maxcall} with strike price $\mathcal{K} = 100$ and independent Poissonian arrivals of stopping opportunities defined in \eqref{eq:poisson} with $d = 5$ independent geometric Brownian motions (without jumps, $\mu_J^j = \lambda_J^j = 0$) and with varying arrival rates $\lambda$ and initial asset prices $X_0$. 
    Column \rm{LSMC} contains the lower-biased estimates based on the random times LSMC algorithm of Section~\ref{sec: fha and ls} and column \rm{Dual} contains the upper-biased estimates \eqref{dual-mc} based on the dual algorithm in Section~\ref{sec: dual alg} and the estimated LSMC policy; standard errors for the dual estimates are the standard errors of the duality gap \eqref{gap-mc}. 
    Both algorithms use the same settings as in Table~\ref{tab:benchmark}, except for the basis functions: as in \cite{AB04}, we use basis functions in the ordered asset prices $X_t^{(1)} \geq \ldots \geq X_t^{(d)}$. 
    The continuation values in the LSMC algorithm are estimated by using the basis functions $(t, x^{(1)}, x^{(2)}) \mapsto t^i\left(x^{(1)}\right)^j\left(x^{(2)}\right)^k$, $i + j + k \leq 5$. 
    Column \rm{Trivial} corresponds to stopping as soon as $Z_{\tau_i} > 0$ or $\tau_i \geq T$ and column \rm{Trivial-PI} contains the results of the policy iteration method \eqref{polit1} with $100\mathord{,}000$ outer paths, $500$ inner paths and a window length $\kappa = K'$. 
    All simulated paths have been generated in antithetic pairs.}}
    \label{tab: LSMC poisson 5}
\end{table}
\begin{table}[h!]
    \centering

{\small\begin{tabular}{rr||rrrrrrrr}
$\lambda$ & $X_0$ & LSMC & (s.e.) & Dual & (s.e.) & Trivial & (s.e.) & Trivial-PI & (s.e.) \\\hline\hline 
1 & 90  & 12.76951 & (0.01119) & 12.78354 & (0.00372) & 9.97336 & (0.00903) & 12.71436 & (0.03361) \\
  & 100  & 20.79173 & (0.01380) & 20.81643 & (0.00473) & 15.64744 & (0.01082) & 20.76846 & (0.04722) \\
  & 110  & 29.99001 & (0.01605) & 30.03046 & (0.00536) & 25.40365 & (0.01261) & 29.96612 & (0.05241) \\
3 & 90  & 16.37745 & (0.01179) & 16.41311 & (0.00401) & 7.85547 & (0.00616) & 16.11181 & (0.04371) \\
  & 100  & 25.70427 & (0.01408) & 25.76604 & (0.00716) & 11.50560 & (0.00716) & 25.26734 & (0.05468) \\
  & 110  & 36.12987 & (0.01597) & 36.21073 & (0.00697) & 21.76685 & (0.00808) & 35.60821 & (0.06066) \\
5 & 90  & 17.18631 & (0.01192) & 17.36115 & (0.09173) & 6.33664 & (0.00484) & 16.69437 & (0.04318) \\
  & 100  & 26.71960 & (0.01418) & 26.97080 & (0.10666) & 9.24555 & (0.00564) & 25.97468 & (0.05221) \\
  & 110  & 37.32898 & (0.01603) & 37.63602 & (0.12062) & 19.52885 & (0.00636) & 36.37071 & (0.05770) \\
\end{tabular}
}
    \caption{Value of a max-call option on random times with jumps ($d=5$). \\ 
    {\footnotesize\textit{Notes.} 
    This table displays estimated values of $Y_0^{(\infty)}$ in the optimal stopping problem for payoff \eqref{eq:maxcall} with strike price $\mathcal{K} = 100$ and independent Poissonian arrivals of stopping opportunities defined in \eqref{eq:poisson} with $d = 5$ independent geometric Brownian motions with jumps, $\mu_J^j = 0.06$ and $\lambda_J^j = 1$ and with varying arrival rates $\lambda$ and initial asset prices $X_0$. 
    The respective algorithms are carried out as described in the caption of Table~\ref{tab: LSMC poisson 5}. 
    }}
    \label{tab: LSMC jumps poisson 5}
\end{table}
\begin{table}[h!]
    \centering

{\small\begin{tabular}{rr||rrrrrr}
$\lambda$ & $X_0$ & $\mathbb{E}[T - \tau^{*} | \tau^{*} \leq T]$ & (s.e.) & $\mathbb{P}(\tau^{*} > T)$ & (s.e.) & $R_O$ & $R_T$ \\\hline\hline 
1 & 90  & 1.18156 & (0.00062) & 0.43962 & (0.00017) & -4.03581 & -0.28296 \\
  & 100  & 1.33435 & (0.00056) & 0.27743 & (0.00014) & -5.25015 & -0.46380 \\
  & 110  & 1.45362 & (0.00053) & 0.18089 & (0.00010) & -6.35926 & -0.63461 \\
3 & 90  & 0.99291 & (0.00053) & 0.34714 & (0.00016) & -1.34970 & -0.25629 \\
  & 100  & 1.11190 & (0.00049) & 0.18858 & (0.00011) & -1.59880 & -0.36390 \\
  & 110  & 1.20112 & (0.00046) & 0.09639 & (0.00006) & -1.81204 & -0.46686 \\
5 & 90  & 0.93253 & (0.00053) & 0.33563 & (0.00016) & -0.80653 & -0.19892 \\
  & 100  & 1.03873 & (0.00049) & 0.17846 & (0.00010) & -0.95009 & -0.27320 \\
  & 110  & 1.12509 & (0.00047) & 0.08612 & (0.00006) & -1.06166 & -0.34326 \\
\end{tabular}
}
    \caption{Various properties and risks of a max-call option on random times ($d=5$). \\ 
    {\footnotesize\textit{Notes.} 
    This table displays the expected time to maturity at exercise, $\mathbb{E}[T-\tau^* | \tau^* \leq T]$, the probability of not exercising, $\mathbb{P}(\tau^* > T)$, and the opportunity risk, $R_O$, and timing risk, $R_T$, for the stopping problem analyzed in Table~\ref{tab: LSMC poisson 5}.
    $\mathbb{E}[T-\tau^* | \tau^* \leq T]$ and $\mathbb{P}(\tau^* > T)$ are based on Monte Carlo estimates using the $2\mathord{,}000\mathord{,}000$ paths of the forward pass of the LSMC algorithm.
    $R_{O}$ and $R_{T}$ are computed as differences between LSMC estimates of both components in their respective definitions.
    All simulated paths have been generated in antithetic pairs.}}
    \label{tab: LSMC poisson 5 risks}
\end{table}
We note also that our methods are not limited to dimension $d=5$ and show in Table~\ref{tab: LSMC poisson 5 long} that good results may be obtained for dimensions easily up to $d=100$.
As expected, when the dimension increases, the duality gaps between random times LSMC and Dual become slightly larger, now leaving some room for improvement by invoking the policy iteration method. 
Specifically, we find that a single step of policy iteration over the random times LSMC based policy yields results that are (very) close to the Dual estimates.

\begin{table}[h!]
    \centering

{\small\begin{tabular}{rr||rrrrrr}
$\lambda$ & $d$ & LSMC & (s.e.) & LSMC-PI & (s.e.) & Dual & (s.e.) \\
\hline\hline
1 & 10  & 28.01203 & (0.01470) & 28.05150 & (0.02168) & 28.06987 & (0.00719) \\
  & 20  & 36.87827 & (0.01666) & 36.88446 & (0.02498) & 36.94684 & (0.00932) \\
  & 50  & 48.54796 & (0.01966) & 48.59752 & (0.02934) & 48.61604 & (0.00929) \\
  & 100  & 57.29451 & (0.02220) & 57.33810 & (0.03228) & 57.35563 & (0.00947) \\
3 & 10  & 34.81370 & (0.01437) & 34.88107 & (0.02419) & 34.88875 & (0.00653) \\
  & 20  & 46.15149 & (0.01539) & 46.24209 & (0.02698) & 46.24795 & (0.00892) \\
  & 50  & 61.23947 & (0.01692) & 61.28060 & (0.03034) & 61.34735 & (0.00969) \\
  & 100  & 72.62265 & (0.01847) & 72.61737 & (0.03395) & 72.75048 & (0.01125) \\
5 & 10  & 36.34189 & (0.01446) & 36.46260 & (0.02644) & 36.50452 & (0.01622) \\
  & 20  & 48.40340 & (0.01527) & 48.53107 & (0.02914) & 48.57309 & (0.01144) \\
  & 50  & 64.49434 & (0.01636) & 64.58681 & (0.03271) & 64.69898 & (0.01402) \\
  & 100  & 76.67618 & (0.01754) & 76.76847 & (0.03541) & 76.93326 & (0.02325) \\
\end{tabular}
}

    \caption{Value of a max-call option on random times (varying $d$). \\ 
    {\footnotesize\textit{Notes.} 
    This table displays estimated values of $Y_0^{(\infty)}$ in the optimal stopping problem for payoff \eqref{eq:maxcall} with strike price $\mathcal{K} = 100$ and independent Poissonian arrivals of stopping opportunities defined in \eqref{eq:poisson} with $d$ independent geometric Brownian motions (without jumps, $\mu_J^j = \lambda_J^j = 0$), with initial asset prices $X_0 = 100.0$, and with varying arrival rates $\lambda$. 
    Column \rm{LSMC} contains the lower-biased estimates based on the random times LSMC algorithm of Section~\ref{sec: fha and ls} and column \rm{Dual} contains the upper-biased estimates \eqref{dual-mc} based on the dual algorithm in Section~\ref{sec: dual alg} and the estimated LSMC policy; standard errors for the dual estimates are the standard errors of the duality gap \eqref{gap-mc}. 
    Both algorithms use the same settings as in Table~\ref{tab:benchmark}, except for the basis functions: as in \cite{AB04}, we use basis functions in the ordered asset prices $X_t^{(1)} \geq \ldots \geq X_t^{(d)}$. 
    The continuation values in the LSMC algorithm are estimated by using the basis functions $(t, x^{(1)}, x^{(2)}) \mapsto t^i\left(x^{(1)}\right)^j\left(x^{(2)}\right)^k$, $i + j + k \leq 5$. 
    Column \rm{LSMC-PI} contains the results of the policy iteration method \eqref{polit1}, starting from the LSMC estimates, with $100\mathord{,}000$ outer paths, $500$ inner paths and a window length $\kappa = K'$. 
    All simulated paths have been generated in antithetic pairs.}}
    \label{tab: LSMC poisson 5 long}
\end{table}

\subsubsection{Reward-dependent Poissonian arrivals of stopping opportunities}
We now introduce a form of dependence between the stopping opportunities $\tau_i$ and the process $X_t$ driving the reward.
Specifically, we consider the payoff \eqref{eq:maxcall} with \textit{potential} stopping opportunities generated according to a Poisson process with rate $\lambda$. 
The potential opportunities to stop constitute only \textit{genuine} opportunities to stop when the lowest asset price among $d$ asset prices, $X_t^{(d)}$, lies above a given threshold value $X_{min}>0$ at that time, i.e.:
\begin{align}
    \nu_{j+1} &= \nu_{j} + \eta_j\,, \quad \eta_j \sim\mathrm{Exp}(\lambda)\,~\mathrm{i.i.d.}, \quad \nu_0 = 0\,,\nonumber \\
    \tau_{i+1} &= \inf\{\nu_{j} : \nu_j > \tau_{i},~X_{\nu_j}^{(d)} \geq X_{min}\}\,, \quad \tau_0 = 0\,. \label{eq:poisson-dependent}
\end{align}

\begin{table}[h!]
    \centering

{\footnotesize\begin{tabular}{rrr||rrrrrrrr}
$d$ & $\lambda$ & $X_0$ & LSMC & (s.e.) & Dual & (s.e.) & Threshold & (s.e.) & Threshold-PI & (s.e.) \\\hline\hline 
5  & 5  & 90  & 10.70117 & (0.00872) & 11.16820 & (0.03041) & 10.23629 & (0.00819) & 11.05121 & (0.02260) \\
   &    & 100  & 21.02792 & (0.01158) & 21.99106 & (0.04818) & 19.90658 & (0.01098) & 21.45137 & (0.03356) \\
   &    & 110  & 32.59407 & (0.01380) & 33.76811 & (0.05761) & 30.77168 & (0.01311) & 32.93245 & (0.04233) \\
   & 10  & 90  & 11.39417 & (0.00874) & 12.10415 & (0.03554) & 10.84847 & (0.00819) & 11.75694 & (0.02403) \\
   &    & 100  & 21.80430 & (0.01154) & 23.16668 & (0.05569) & 20.66690 & (0.01095) & 22.19446 & (0.03490) \\
   &    & 110  & 33.43591 & (0.01373) & 35.07910 & (0.06711) & 31.74492 & (0.01311) & 33.69887 & (0.04304) \\
10  & 5  & 90  & 12.56657 & (0.00848) & 13.27405 & (0.04561) & 12.12012 & (0.00823) & 13.11287 & (0.02374) \\
   &    & 100  & 26.00339 & (0.01161) & 27.40261 & (0.06667) & 24.75869 & (0.01108) & 26.93033 & (0.03735) \\
   &    & 110  & 40.34155 & (0.01417) & 42.38168 & (0.08850) & 38.29985 & (0.01360) & 41.35958 & (0.04856) \\
   & 10  & 90  & 13.58689 & (0.00847) & 14.64311 & (0.05502) & 13.03969 & (0.00813) & 14.26980 & (0.02578) \\
   &    & 100  & 27.20783 & (0.01150) & 29.39729 & (0.08570) & 25.94181 & (0.01110) & 28.13476 & (0.03892) \\
   &    & 110  & 41.60750 & (0.01398) & 44.69199 & (0.10200) & 39.67943 & (0.01358) & 42.49327 & (0.04899) \\

\end{tabular}
}

    \caption{Value of a max-call option on reward-dependent random times (varying $d$). \\ 
    {\footnotesize\textit{Notes.} 
    This table displays estimated values of $Y_0^{(\infty)}$ in the optimal stopping problem for payoff \eqref{eq:maxcall} with strike price $\mathcal{K} = 100$ and reward-dependent Poissonian arrivals of stopping opportunities defined in \eqref{eq:poisson-dependent} with $d$ independent geometric Brownian motions (without jumps, $\mu_J^j = \lambda_J^j = 0$) and with varying arrival rates $\lambda$ and initial asset prices $X_0$, where $X_{min} = 60$. 
    All columns are determined as in Table~\ref{tab: LSMC poisson 5}, except for the policy iteration columns: \textrm{Threshold} and \textrm{PI-threshold} are based on the backward dynamic programming method with a threshold policy \eqref{bdp:thresh} rather than a trivial initial policy as in Table~\ref{tab: LSMC poisson 5}.  
    }}
    \label{tab: poisson dependent alt}
\end{table}

To illustrate the performance of the policy iteration method when the initial policy is already somewhat informed, we now use the backward dynamic programming procedure of Section~\ref{sec: fha and ls} with the threshold-based strategy in \eqref{bdp:thresh} to determine the initial policy. 
Table~\ref{tab: poisson dependent alt} displays the results. 

Table~\ref{tab: poisson dependent alt} shows that a single step of policy iteration over the threshold policy based on \eqref{bdp:thresh} may already lead to estimates that are close to those obtained by the random times LSMC method \eqref{bdp:LSMC} with many basis functions, and appear to be always larger than the random times LSMC lower bound. 
If the dimension cardinality of the underlying process is relatively low, the random times LSMC method will usually be computationally faster, since the simulation due to the iterated policy requires one degree of nesting. 
However, it is known that the LSMC method suffers from the curse of dimensionality, in contrast to the policy iteration method.\footnote{In this respect, we note that the relatively large duality gaps for LSMC may partially be attributed to the (somewhat coarse) choice of basis functions: although the arrivals of the opportunities depend on $X_t^{(d)}$, no basis functions in this term are used for $d \geq 3$. 
As the effect of this is less pronounced for small $\lambda$ and the PI method uses fewer paths due to its increased computational complexity, the associated larger standard deviations partially nullify the improvement for smaller $\lambda$, whereas the gains from policy iteration become more sizeable as $\lambda$, and correspondingly the duality gaps for LSMC, increases.} 

This phenomenon is in line with findings from the existing literature on deterministic stopping opportunities, and illustrates how policy iteration may help mitigate the curse of dimensionality. 
For instance, in \cite{BKS2008} it was demonstrated that for an exotic cancelable coupon swap with 20 underlying assets, the classical Longstaff-Schwartz LSMC algorithm with 300 basis functions did not give satisfactory results (duality gaps of around 5$\%$ in relative terms). 
This stands in stark contrast with a one step improvement of an input policy, obtained via the classical Longstaff-Schwartz algorithm with a low number of basis functions, which yields duality gaps of under 1.5$\%$ in relative terms, and moreover is up to 8 times faster. 

It is likely that for (very) high-dimensional underlying processes, a combination of random times LSMC with policy improvement may be (much) more efficient than running the random times LSMC algorithm alone with many basis functions. 
We consider this an interesting topic for a detailed further study, but beyond the scope of the present paper. 
\setcounter{equation}{0}
\section{Real options application}\label{sec:ex}
In this section, we analyze a real options pricing problem, to further illustrate the managerial implications of randomly arriving opportunities to stop. 
We suppose that an agent owns an illiquid asset, e.g., a piece of real estate or an art object, where illiquid means that a potential buyer is not always directly available and may arrive randomly. 
The agent attaches a monetary time-invariant value of $\mathcal{K}$ to the asset, but she may try to sell the asset at the prevailing time-varying market price $X_t$, if a buyer appears at that time instance, i.e., at $t=\tau_{i}$. 
Of course, she will only choose to do so if the market price is sufficiently large, and in any case larger than $\mathcal{K}$.\footnote{This application follows the spirit of \cite{DP94} but with randomly arriving stopping opportunities.}

The overall value of the asset for the agent thus consists of two parts: the monetary value $\mathcal{K}$ attached to it and the option value of selling the asset at a price $X_{\tau_{i}}$ beyond $\mathcal{K}$. 
We suppose that $X_t$ evolves according to a geometric Brownian motion with log-normal jump sizes $J_i$, with mean 1 and that occur at the arrival times of buyers, i.e., at $t=\tau_i$, hence:
\begin{equation}
    \dd{X}_t = (r-\delta)X_t\dd{t} + \sigma X_t \dd{W}_t + X_{t-} \sum_{t=\tau_i }(J_i-1)\,, \quad J_i \sim \mathrm{LogNormal}(1, \sigma_J)\,~\mathrm{i.i.d.}, \label{eq:house-jumps}
\end{equation}
with $X_{0}>$, $r$, $\delta$, $\sigma>0$, $\sigma_{J}>0$ known constants.
These dynamics underscore the advantage of simulation-based methods over Markov generator-based methods: the infinitesimal Markov generator for \eqref{eq:house-jumps} would be non-local, which means that its discretized matrix is non-sparse; in such cases, finite-difference methods typically become numerically costly in higher dimensions.\footnote{See \cite{CV05} for an in-depth analysis.} 

We will study the option value, which means that the reward process of interest in this section is given by
 $   Z_t = e^{-rt}(X_t - \mathcal{K})^+\,,  $ 
corresponding to the discounted payoff of a call option contract. 
We suppose that the agent is only considering to sell the asset over a three years time window, i.e., $T = 3$; 
if the asset is not sold by that time, the agent will retain the asset. 
In the previous section, we have seen that in a low-dimensional setting random times LSMC performs comparably to policy iteration. 
As this example is meant to illustrate the implications of randomly arriving stopping opportunities, and in particular the role of feedback, we focus only on lower- and upper-biased estimates based on the random times LSMC and dual algorithms, respectively. 

In this application, it is natural to consider various dynamics for the arrivals of prospective buyers and thus for the stopping opportunities in \eqref{eq:stop}. 
Buyers may arrive independently of everything else, or with a rate that depends on previous arrivals or the current market price of the asset. 
We will illustrate how the probability of not exercising, $\mathbb{P}(\tau^* > T)$, the expected time to maturity at exercise, $\mathbb{E}[T-\tau^*|\tau \leq T]$, and the opportunity and timing risks, $R_O$ and $R_T$, vary with the arrival dynamics of stopping opportunities and the dependence structure between arrivals of stopping opportunities and asset dynamics. 

\subsection{Reward-independent arrival of  buyers}
We begin by considering two cases in which buyers arrive either at fixed moments, say at the end of each year, i.e.,
\begin{equation}
    \tau_{i+1}^{\mathrm{Ber}} = \tau_i^{\mathrm{Ber}} + \frac{T}{K}\,, \quad \tau_0 = 0\,, \quad i=0,\ldots,K-1,\label{eq:ex-bermudan}
\end{equation}
with $K = T = 3$, or according to a Poisson process that is independent of the reward $Z_t$, i.e.,
\begin{equation}
    \tau_{i+1}^{\mathrm{Poi}} = \tau_i^{\mathrm{Poi}} + \eta_i \,, \quad \tau_0 = 0\,, \quad \eta_i \sim \mathrm{Exp}(\lambda)\,~\mathrm{i.i.d.}, \quad i=0,1,\ldots,\label{eq:ex-poisson}
\end{equation}
with $\lambda = 1$. 
Note that when opportunities follow \eqref{eq:ex-poisson}, there is both timing risk and opportunity risk involved in the stopping decision at time $\tau_i^{\mathrm{Poi}}$: the length of the time increment $\tau_{i+1}^{\mathrm{Poi}} - \tau_i^{\mathrm{Poi}}$ is \textit{ex ante} unknown and it may happen that $\tau_{i+1}^{\mathrm{Poi}}$ occurs after $T$. 

The presence of timing risk and opportunity risk does not entail that decision problems where opportunities follow \eqref{eq:ex-bermudan} are in general always to be preferred to decision problems where \eqref{eq:ex-poisson} applies: as a simple example, suppose $Z_t = 1 + \sin(\frac{2K\pi}{T}t + a)\,,~a \in \{0, \pi\}$, and consider $\tau_i^{\mathrm{Ber}} < T$. 
If $a = 0$, $Z_{\tau_{i}^{\mathrm{Ber}}} \equiv 0$, whereas if $a = \pi$, $Z_{\tau_{i}^{\mathrm{Ber}}} \equiv 2$, which are the minimal and maximal rewards possible, respectively. 
So, if $a = 0$, the agent clearly prefers Poissonian opportunities despite the presence of opportunity risk, as these attain a positive reward with positive probability. 
However, if $a = \pi$, the agent prefers Bermudan opportunities, which achieve higher rewards than Poissonian opportunities with probability one as both the timing risk and opportunity risk make Poissonian arrivals less favorable. 

Both forms of risk are also present, and sizeable, under reward-independent Poissonian arrivals in the context of the reward in the present section, as illustrated in the first two rows of Table~\ref{tab: example independent}, leading on average to stopping earlier. 

\begin{table}[ht!]
    \centering

{\footnotesize\begin{tabular}{r||rrrrrrrr}
Arrivals & LSMC & (s.e.) & Dual & (s.e.) & $\mathbb{E}[\mathrm{TtM}]$ 
& $\mathbb{P}(\tau^{*} > T)$ 
& $R_O$ & $R_T$ \\\hline \hline
Bermudan & 7.54039 & (0.00891) & 7.54071 & (0.00014) & 0.42816 
& 0.00000 
& 0.00000 & 0.00000 \\
Poisson & 6.04603 & (0.00792) & 6.04657 & (0.00022) & 1.58655 
& 0.66616 
& -1.19929 & -0.29507 \\
Hawkes(0.3) & 5.68566 & (0.00773) & 5.68705 & (0.00052) & 1.60190 
& 0.68344 
& -1.48159 & -0.37314 \\
Hawkes(0.9) & 4.08941 & (0.00656) & 4.09257 & (0.00144) & 1.74170 
& 0.76196 
& -2.79049 & -0.66049 \\
Feedback(-1) & 2.07946 & (0.00383) & 2.07947 & (0.00001) & 1.62801 
& 0.76158 
& -4.18977 & -1.27115 \\
Feedback(1) & 6.99131 & (0.00815) & 6.99558 & (0.00126) & 1.55637 
& 0.64199 
& -0.60036 & 0.05128 \\
\end{tabular}}
    \caption{Option value of the asset with buyers arriving at random times. \\ 
    {\footnotesize\textit{Notes.} 
        This table displays estimated values of $Y_0^{(K)}$ and its properties in the optimal stopping problem for reward $Z_t = e^{-rt}(X_t - 100)^+$ with price process \eqref{eq:house-jumps} (without jumps, $\sigma_J = 0\%$) and specifications \eqref{eq:ex-bermudan}--\eqref{eq:ex-feedback} for \eqref{eq:stop}. 
        Column \rm{LSMC} contains the lower-biased estimates based on the random times LSMC algorithm of Section~\ref{sec: fha and ls} and column \rm{Dual} contains the upper-biased estimates based on the dual algorithm in Section~\ref{sec: dual alg} and the estimated LSMC policy as input for \eqref{doob-mart}; standard errors for the dual estimates are the standard errors of the duality gap \eqref{gap-mc}. 
        Both algorithms use the same settings as in Table~\ref{tab:benchmark}, except for the choice of basis functions: denote by $d_{\mathrm{eff}}$ the effective state dimension (number of assets + time variable + any state variables needed for the arrival process), then we use $\{(y_1,\ldots,y_{d_{\mathrm{eff}}}) \mapsto \prod_{i=1}^{d_{\mathrm{eff}}}y^{j_i} : \sum_{i=1}^{d_{\mathrm{eff}}}j_i \leq 5\}$ as basis functions. 
        $\mathbb{E}[\mathrm{TtM}]:=\mathbb{E}[T-\tau^*|\tau \leq T]$ is the expected time to maturity at exercise, $\mathbb{P}(\tau^* > T)$ the probability of not exercising, $R_O$ measures the opportunity risk and $R_T$ the timing risk. 
        The measures $\mathbb{E}[\mathrm{TtM}]$, $\mathbb{P}(\tau^* > T)$, $R_O$ and $R_T$ are computed as in Table~\ref{tab: LSMC poisson 5 risks}.}}
    \label{tab: example independent}
\end{table}

\subsection{Arrival rates with feedback mechanisms}
Next, we analyze the case that the arrival of buyers is affected by previous arrivals or by the current market price of the asset.
This may be modeled by a c\`adl\`ag point process $N_t$ for the number of buyers with (left-continuous, predictable) stochastic jump intensity $\lambda_t$, where, as $\Delta\downarrow 0$,
\begin{align*}
    \mathbb{P}(N_{t+\Delta} = N_{t}|\mathcal{F}_t) &= 1 - \lambda_t\Delta + o(\Delta) \,,\quad \mathbb{P}(N_{t+\Delta} - N_{t} = 1|\mathcal{F}_t) = \lambda_t\Delta + o(\Delta)\,,\ \nonumber \\ \mathbb{P}(N_{t+\Delta} - N_{t} \geq 2|\mathcal{F}_t) &= o(\Delta)\,,\quad 
    \tau_{i}^{\lambda} = \inf\{t:N_t = i\}\,. 
\end{align*}

To let the arrival rate depend upon previous arrivals, we suppose that $(N_t,\lambda_t )$ follows a Hawkes process with exponential kernel (\cite{H71}), which introduces a self-exciting feature:
\begin{equation}
    \dd{\lambda}_t^{\mathrm{Haw}} = \alpha(\psi - \lambda_t^\mathrm{Haw})\dd{t} + \beta \dd{N}_t\,. \label{eq:ex-hawkes}
\end{equation}
That is, upon the occurrence of arrivals, the arrival intensity ramps up, by $\beta>0$, making future arrivals temporarily more likely, after which there is mean reversion at rate $\alpha>\beta$ to a level $\psi>0$.
For the sake of comparison, we consider two parameter specifications for $\lambda_t$ such that $\mathbb{E}[\lambda_t] = 1\,,~\forall t$, which means that the expected arrival rate at any point in time is the same as under the Poissonian arrivals \eqref{eq:ex-poisson}. 
In turn this means that $\mathbb{E}[N_T] = 3$, the (expected) number of buyers under the Bermudan arrivals \eqref{eq:ex-bermudan}.
In particular, we consider either $\alpha = 1,~\beta = 3/10,~ \psi = 7/10,~\lambda_0 = 1$ (denoted $\mathrm{Hawkes(0.3)}$) or $\alpha = 1,~\beta = 9/10,~\psi = 1/10,~\lambda_0 = 1$ (denoted $\mathrm{Hawkes(0.9)}$). 

\begin{table}[t]
    \centering

{\footnotesize\begin{tabular}{r||rrrrrrrr}
Arrivals & LSMC & (s.e.) & Dual & (s.e.) & $\mathbb{E}[\mathrm{TtM}]$ 
& $\mathbb{P}(\tau^{*} > T)$ 
& $R_O$ & $R_T$ \\\hline \hline
Bermudan & 8.79715 & (0.01045) & 8.79767 & (0.00023) & 0.40733 
& 0.00000 
& 0.00000 & 0.00000 \\
Poisson & 7.53311 & (0.00937) & 7.53705 & (0.00127) & 1.62419 
& 0.65759 
& -1.19040 & -0.07364 \\
Hawkes(0.3) & 7.12238 & (0.00914) & 7.13228 & (0.00537) & 1.64638 
& 0.67436 
& -1.46015 & -0.21462 \\
Hawkes(0.9) & 5.37431 & (0.00817) & 5.38934 & (0.00515) & 1.77864 
& 0.75885 
& -2.78202 & -0.64082 \\
\end{tabular}}
    \caption{Option value of the asset with buyers arriving at random times and where arrivals cause price jumps. \\ 
    {\footnotesize\textit{Notes.} 
    This table displays estimated values of $Y_0^{(K)}$ and its properties in the optimal stopping problem for reward $Z_t = e^{-rt}(X_t - 100)^+$ with price process \eqref{eq:house-jumps} (with jumps, $\sigma_J = 10\%$) and specifications \eqref{eq:ex-bermudan}--\eqref{eq:ex-hawkes} for \eqref{eq:stop}. 
    Details for the respective algorithms are provided in the caption of Table~\ref{tab: example independent}.
    $\mathbb{E}[\mathrm{TtM}]:=\mathbb{E}[T-\tau^*|\tau \leq T]$ is the expected time to maturity at exercise, $\mathbb{P}(\tau^* > T)$ the probability of not exercising, $R_O$ measures the opportunity risk and $R_T$ the timing risk.
    }}
    \label{tab: example jumps}
\end{table}

In Tables~\ref{tab: example independent} and~\ref{tab: example jumps}, we observe that the option value under Hawkes arrivals is lower than that under Poissonian arrivals, which is due to both a larger (absolute value of) opportunity and timing risk, both of which are moreover increasing with the degree of clustering. 
Indeed, due to the self-exciting nature of Hawkes arrivals and the condition that $\mathbb{E}[\lambda_t] = 1$, the variation of inter-arrival times increases with the degree of clustering, i.e., with $\beta$; hence, the timing risk is larger.
Furthermore, under the condition $\mathbb{E}[\lambda_t] = 1$, one may show that $\mathbb{P}(\tau_1^{\mathrm{Poi}} > T) < \mathbb{P}(\tau_1^{\mathrm{Haw}} > T)$, which explains why the opportunity risk is larger under Hawkes arrivals. 
From a managerial point of view, when only $\mathbb{E}[\lambda_t]$ or $\mathbb{E}[N_T]$ is available, it is thus important not to confuse self-exciting arrivals with Poissonian or Bermudan arrivals: 
Hawkes arrivals reduce the option value and require stopping earlier on average.
Upon comparing Tables~\ref{tab: example independent} and~\ref{tab: example jumps}, we see that the option value increases overall when arrivals of stopping opportunities lead to price jumps.
Indeed, the expected call option payoff benefits from upward potential due to the price jumps.

We may also let the arrival rate of stopping opportunities depend upon current market prices by setting
\begin{equation}
    \lambda_t^\phi := \lambda(t,X_t,\phi), \quad \text{where} \quad \lambda(t, x,\phi) := (1-\phi) + 2\phi F_{X_t}(x)\,,\quad \phi \in [-1,1]\,, \label{eq:ex-feedback} 
\end{equation}
assuming $\sigma_{J}\equiv0$ (i.e., no price jumps), with $F_{X_t}$ denoting the cumulative distribution of $X_t$.
This entails that $\lambda_t^\phi$ and $X_t$ move in the same direction if $\phi>0$ and in opposite directions if $\phi<0$. 
Note that $\mathbb{E}[\lambda_t^\phi] = 1$.
We denote the extremal cases $\phi = -1$ and $\phi = 1$ by $\mathrm{Feedback}(-1)$ and $\mathrm{Feedback}(1)$, respectively. 

The associated results are displayed in the last two rows of Table~\ref{tab: example independent}. 
Note that Poissonian arrivals correspond to $\mathrm{Feedback}(0)$.
Upon comparing the values we see that both the opportunity risk and the timing risk become more pronounced (i.e., more negative) as the feedback becomes more negative.
Under positive feedback, when the option becomes attractive (i.e., high market price of the asset), the corresponding increased arrival rate of stopping opportunities is beneficial. 
These insights underscore the importance of the ability to account for the relevant dependence structure and feedback mechanisms. 

\subsection{Cross-exciting arrivals}
We finally analyze a more complex form of dependence between arrivals and prices.
Specifically, we assume price jumps and arrivals of stopping opportunities are jointly driven by a two-dimensional Hawkes process with exponential kernel  (\cite{H71,ALP14,ACL15}), which allows for cross-excitation between the number of stopping opportunities, $N_t^1$, and the number of price jumps, $N_t^2$:
\begin{align}
    &\mathbb{P}(N_{t+\Delta}^i = N_{t}^i|\mathcal{F}_{t}) = 1 - \lambda_t^i\Delta + o(\Delta) \,,\quad \mathbb{P}(N_{t+\Delta}^i - N_{t}^i = 1|\mathcal{F}_{t}) = \lambda_t^i\Delta + o(\Delta)\,,\nonumber \\
    &\mathbb{P}(N_{t+\Delta}^i - N_{t}^i \geq 2|\mathcal{F}_{t}) = o(\Delta)\,, \quad \tau_i = \inf\{t:N_t^1 = i\}\,, \nonumber\\
    &\dd{\lambda}_t^i = \alpha(\psi - \lambda_t^i)\dd{t} + \beta \dd{N}_t^i + \rho \dd{N}_t^{3-i}\,, \quad i = 1, 2\,, \label{ex:hawkes-lambda}\\ 
    &\dd{X}_t = (r-\delta)X_t\dd{t} + \sigma X_t \dd{W}_t + X_{t-} J_{N_{t}^2}\dd{N}_t^2\,, \quad J_i \sim \mathrm{LogNormal}(1, \sigma_J)\,~\mathrm{i.i.d.} \label{ex:hawkes-X}
\end{align}
We set $r = 5\%,~\delta = 10\%,~\sigma = 20\%$ as before and fix $\alpha = 1$ and $\psi = 0.1$. 
We will vary the cross-excitation parameter $\rho>0$ while keeping $\mathbb{E}[\lambda_t^i] = 1$, $i=1,2$, which dictates that the self-excitation parameter $\beta = 0.9 - \rho$ (\cite{ACL15}). 
In Table~\ref{tab: example hawkes full}, we denote different parameters $(\beta, \rho)$ as $\mathrm{Hawkes}(\beta, \rho)$. 
Upon comparing rows 2--4, it is clearly visible that an increase in the degree of cross-excitation, i.e., larger $\rho$, leads to larger values for the stopping problem, in line with (but expanding the scope of) the positive feedback case in the previous setting. 

An interesting observation about the role of dependence can be made by comparing rows $\mathrm{Hawkes}(0.9)$ in Table~\ref{tab: example jumps} and $\mathrm{Hawkes}(0.9, 0)$ in Table~\ref{tab: example hawkes full}: in either case, the relevant jump times are the event times of a 1-dimensional exponential Hawkes process with the same parameters. 
In Table~\ref{tab: example jumps}, both the price jumps and stopping opportunities are generated by the same Hawkes process, whereas in Table~\ref{tab: example hawkes full}, both are generated by a different, independent (i.e., non-cross-exciting), Hawkes process; under this independence, the value of the stopping problem is substantially lower. 
This reinforces our earlier finding that the dependence between the reward and the arrival of stopping opportunities plays an important role in optimal stopping with randomly arriving stopping opportunities, and warrants the development of methods that can handle them.  
\begin{table}[t]
    \centering

{\footnotesize\begin{tabular}{r||rrrrrrrr}
Arrivals        
& LSMC   & (s.e.)  & Dual    & (s.e.)  & $\mathbb{E}[\mathrm{TtM}]$ 
& $\mathbb{P}(\tau^{*} > T)$ 
& $R_O$ & $R_T$ \\\hline \hline
Bermudan 
& 8.79715 & (0.01045) & 8.79767 & (0.00023) & 0.40733 
& 0.00000 
& 0.00000 & 0.00000 \\
Hawkes(0.9,0) 
& 4.67110 & (0.00773) & 4.68102 & (0.00423) & 1.71095 
& 0.76522 
& -3.40758 & -0.71847 \\
Hawkes(0.6,0.3) 
& 5.58933 & (0.00853) & 5.62453 & (0.01209) & 1.63582 
& 0.73352 
& -2.61380 & -0.59402 \\
Hawkes(0.3,0.6) 
& 5.85595 & (0.00864) & 5.88006 & (0.00655) & 1.61496 
& 0.72337 
& -2.39815 & -0.54305 \\
\end{tabular}}
    \caption{Option value of the asset with buyers arriving at random times where price and arrivals are driven by a Hawkes process. \\ 
    {\footnotesize\textit{Notes.} 
    This table displays estimated values of $Y_0^{(K)}$ and its properties in the optimal stopping problem for reward $Z_t = e^{-rt}(X_t - 100)^+$ with price process and arrivals given by \eqref{ex:hawkes-lambda}--\eqref{ex:hawkes-X}, where we set $\sigma_{J}=0.1$. 
    Details for the respective algorithms are provided in the caption of Table~\ref{tab: example independent}.
    $\mathbb{E}[\mathrm{TtM}]:=\mathbb{E}[T-\tau^*|\tau \leq T]$ is the expected time to maturity at exercise, $\mathbb{P}(\tau^* > T)$ the probability of not exercising, $R_O$ measures the opportunity risk and $R_T$ the timing risk.
    }}
    \label{tab: example hawkes full}
\end{table}
\section{Conclusion}\label{sec:con}

In this paper, we have studied general optimal stopping problems with randomly arriving opportunities to stop.
Such problems occur naturally when markets display frictions (e.g., limited liquidity, incomplete information, transaction costs) in a wide variety of settings (e.g., financial contracts, real options, bank runs, natural resources).
To encompass these, we allow for general dynamics of the stopping opportunities and the reward process, including general forms of stochastic dependence among them.
We have recast the original problem into an equivalent integer-time, infinite-horizon problem. 
In this context, we have developed, theoretically and pseudo-algorithmically, a policy improvement approach, a dual martingale method, and a regression procedure for optimal stopping problems with infinite horizon. 

We have illustrated the good performance of our algorithms in a variety of applications, efficiently yielding (very) small duality gaps, with cardinality of the state space that can easily go up to 100.
We have shown that, due in part to the presence of timing risk and opportunity risk, both the optimal value of the stopping problem and the optimal timing of the stopping decisions are substantially affected when stopping opportunities become risky. 
In the context of a real options problem, this is illustrated by analyzing several models of how prospective buyers may present themselves for a valuable asset: at deterministic times, following Poissonian arrivals, self-exciting (Hawkes) arrivals and price-dependent arrivals. 
Furthermore, we have shown that in our flexible framework, problems with feedback, in which the arrival of buyers may lead to sudden price changes, can also be treated. 

In order to deal with \textit{multiple} stopping with randomly arriving opportunities, one may follow a similar path. 
In particular, the corresponding approaches for multiple stopping with finite horizon in the literature, see for example \cite{BS06} for policy iteration, \cite{S12} for the dual martingale representation, and \cite{CT08} for regression methods, may be extended to infinite-horizon problems in the spirit of this paper.   



The methods in this paper were tailored to the specific properties of \textit{stopping} problems; in future work, we intend to extend our simulation-based approach to general \textit{control} problems with random control opportunities.

{\singlespacing
 }

\newpage 
\appendix
\onehalfspacing
\section{Proofs and Auxiliary Technical Results}
\subsection{Properties of Problem (\ref{eq:stopproblem})}\label{app:properties-stop}
\begin{proof}[Proof of Proposition~\ref{prop:proprandtime}]
    (i) For fixed $\mathfrak{n}\in\mathcal{N}$ with $\mathfrak{n}\leq K$ a.s., $\mathfrak{t} := \tau_{\mathfrak{n}}$ is an $\mathbb{F}$-stopping time with $\mathfrak{t}\in\mathcal{T}_{K}$, so $\sup_{\mathfrak{n}\in\mathcal{N},~\mathfrak{n}\leq K}\mathbb{E}\left[U_{\mathfrak{n}}\right] \leq Y_{0}^{(K)}$. Indeed, one trivially has that $\tau_{\mathfrak{n}\left(\omega\right)}(\omega)\in\mathcal{T}_{K}\left(\omega\right)$, and moreover,
    \[ \left\{\mathfrak{t}\leq t\right\}=\left\{\tau_{\mathfrak{n}} \leq t\right\} = \bigcup_{i=0}^{K}\left\{\tau_{\mathfrak{i}}\leq t\right\}\cap\left\{\mathfrak{n}=i\right\}\in\mathcal{F}_{t}\,, \]
    since, for any $i\in\mathbb{N}_{0}$, $\left\{\tau_{i}\leq t\right\}\in\mathcal{F}_{t}$ and $\left\{\mathfrak{n}=i\right\}\in\mathcal{F}_{\tau_{\mathfrak{i}}}$, and so $\left\{\tau_{\mathfrak{i}}\leq t\right\}\cap\left\{\mathfrak{n}=i\right\}\in\mathcal{F}_{t}$, due to \eqref{eq:stopsig}.

    (ii) Now let $\mathfrak{t}$ be any $\mathbb{F}$-stopping time with a.s.~$\mathfrak{t}\left(\omega\right)\in\mathcal{T}_{K}\left(\omega\right)$. Define $\mathfrak{n}:\Omega\rightarrow\mathbb{N}_{0}$ by
    \[ \mathfrak{n}\left(\omega\right)=\sum_{i=0}^{K}i1_{\left\{\mathfrak{t}\left(\omega\right)=\tau_{i}\left(\omega\right)\right\}}\,. \]
    Then, for any $j$, since both $\mathfrak{t}$ and $\tau_{j}$ are $\mathbb{F}$-stopping times, it holds that
    \[ \left\{\mathfrak{t}\left(\omega\right)=\tau_{j}\left(\omega\right)\right\}\in\mathcal{F}_{\mathfrak{t}}\cap\mathcal{F}_{\tau_{j}}\,; \]
    see e.g., \cite[Lemma 2.16]{KS98}. Hence,
    \[ \left\{\mathfrak{n}\left(\omega\right)=j\right\}=\left\{\mathfrak{t}\left(\omega\right)=\tau_{j}\left(\omega\right)\right\}\in\mathcal{F}_{\tau_{j}}=\mathcal{G}_{j}\,, \]
    i.e., $\mathfrak{n}$ is a discrete $\mathbb{G}$-stopping time with $\mathfrak{n}\leq K$ a.s., and moreover
    \[ \mathbb{E}\left[Z_{\mathfrak{t}}\right] =\mathbb{E}\left[\sum_{i=1}^{K}U_{\mathfrak{i}}1_{\left\{\mathfrak{t}=\tau_{\mathfrak{i}}\right\}}\right] =\mathbb{E}\left[\sum_{i=1}^{K}U_{i}1_{\left\{\mathfrak{n}=i\right\}}\right] = \mathbb{E}\left[U_{\mathfrak{n}}\right]\,. \]
    Thus, $Y_{0}^{(K)}\leq\sup_{\mathfrak{n}\in\mathcal{N},~\mathfrak{n}\leq K}\mathbb{E}\left[U_{\mathfrak{n}}\right]$.
\end{proof}
\begin{lemma}
\label{lemma:insert-T} Consider a sequence $\left(  \tau_{i}^{\prime}\right)
$ that satisfies \eqref{eq:stop}. Next, construct $\left(  \tau_{i}\right)  $
by setting $\tau_{0}=0$ and, for $i\geq1$,
\begin{equation}
\tau_{i}=\tau_{i}^{\prime}1_{\left\{  \tau_{i}^{\prime}<T\right\}
}+T1_{\left\{  \tau_{i-1}^{\prime}<T\,\wedge\,\tau_{i}^{\prime}\geq T\right\}
}+\tau_{i}^{\prime}1_{\left\{  \tau_{i-1}^{\prime}\geq T\right\}  }\,.
\label{detT1}%
\end{equation}
Then, $\left(  \tau_{i}\right)  $ defined by \eqref{detT1} satisfies both
\eqref{eq:stop} and \eqref{detT}.
\end{lemma}

\begin{proof}
We first show that \eqref{detT1} is a stopping time. 
For $s<T$, we have that
\begin{equation}
\left\{  \tau_{i}\leq s\right\}  =\left\{  \tau_{i}^{\prime}\leq s\right\}
\in\mathcal{F}_{s}\,, \label{ls}%
\end{equation}
since $\tau_{i}\in\left\{  T,\tau_{i}^{\prime}\right\}  $. 
Furthermore, we have that
\begin{align*}
\left\{  \tau_{i}=T\right\}   &  =\left\{  \tau_{i-1}^{\prime}<T\right\}
\cap\left\{  \tau_{i}^{\prime}<T\right\}  ^{c}\,\\
&  =\left(  \bigcup\limits_{n}\left\{  \tau_{i-1}^{\prime}\leq T-\frac{1}%
{n}\right\}  \right)  \cap\left(  \bigcup\limits_{n}\left\{  \tau_{i}^{\prime
}\leq T-\frac{1}{n}\right\}  \right)  ^{c}\,\\
&  \in\mathcal{F}_{T}\,,\quad\text{ hence }\\%
\left\{  \tau_{i}\leq T\right\}   &  =\left\{  \tau_{i}=T\right\}
\cup\left\{  \tau_{i}<T\right\}  \,\\
&  =\left\{  \tau_{i}=T\right\}  \cup\bigcup\limits_{n}\left\{  \tau_{i}\leq
T-\frac{1}{n}\right\}  \in\mathcal{F}_{T}\,,
\end{align*}
using \eqref{ls}. 
For $s>T$, it holds that
\begin{align*}
\left\{  \tau_{i}\leq s\right\}   &  =\left\{  \tau_{i}\leq T\right\}
\cup\left\{  T<\tau_{i}\leq s\right\}  \,\\
&  =\left\{  \tau_{i}\leq T\right\}  \cup\left(  \left\{  \tau_{i}^{\prime
}\leq T\right\}  ^{c}\cap\left\{  \tau_{i}^{\prime}\leq s\right\}  \right)
\in\mathcal{F}_{s}\,.
\end{align*}

We next show that $\tau_{i+1}>\tau_{i}$ for all $i\in\mathbb{N}_{0}$. 
Suppose that $\tau_{i+1}\leq\tau_{i}$. 
Since $\tau_{i}\in\left\{  T,\tau_{i}^{\prime
}\right\}  $, $\tau_{i+1}\in\left\{  T,\tau_{i+1}^{\prime}\right\}  ,$ and by
assumption $\tau_{i+1}^{\prime}>\tau_{i}^{\prime}$, we must have that (a)
$\left(  \tau_{i}=T\right)  \wedge\left(  \tau_{i+1}\leq T\right)  $ or that
(b) $\left(  \tau_{i+1}=T\right)  \wedge\left(  \tau_{i}\geq T\right)  .$ 
On the one hand, (a) implies $\tau_{i}^{\prime}\geq T$ by (\ref{detT1}), and this
in turn implies $\tau_{i+1}=$ $\tau_{i+1}^{\prime}$ by replacing $i$ by $i+1$ in
(\ref{detT1}). 
Hence, we get $\tau_{i+1}=$ $\tau_{i+1}^{\prime}>T$,
a contradiction. 
On the other hand, (b) implies $\tau_{i}^{\prime}<T$ by
(\ref{detT1}), and this in turn implies $\tau_{i}=\tau_{i}^{\prime}$ by (\ref{detT1}).
Hence, we get $\tau_{i}<T$, a contradiction. 
Finally, clearly, $\lim
_{i\rightarrow\infty}\tau_{i}=\infty$, and so \eqref{eq:stop} and \eqref{detT} follow.
\end{proof}

\begin{corollary}
\label{cor:insert-T} If in Lemma~\ref{lemma:insert-T} $(\tau_{i}^{\prime})$
satisfies \eqref{eq:boundZ} and $\mathbb{E}\left[  Z_{T}^{\alpha}\right]  $ is
finite, then $(\tau_{i})$ also satisfies \eqref{eq:boundZ}, albeit possibly
with a larger value of $B$.
\end{corollary}

\begin{proof}
Note that
\[
{\sup_{i\in\mathbb{N}_{0}}Z_{\tau_{i}}^{\alpha}\leq\sup_{i\in\mathbb{N}_{0}%
}Z_{\tau_{i}^{\prime}}^{\alpha}+Z_{T}^{\alpha}\,,}%
\]
and then take expectations on both sides. 
\end{proof}
\begin{remark}\label{rem:analytic-BDG}
    In practice, the bound \eqref{eq:boundZ}, hence \eqref{eq:boundZ1}, may be determined by simulation. 
    For an analytic alternative, let us observe that condition \eqref{eq:boundZ} is implied by the somewhat stronger condition
    \[ \mathbb{E}\left[\sup_{0\leq t\leq T}Z_{t}^{\alpha}\right] \leq B^{\alpha}\,, \quad \mathrm{for\ some}\ B>0\ \mathrm{and}\ \alpha>1\,. \]
    In the case in which the reward can be written in the form $Z_{t}=g(t,M_{t})$ for some $d$-dimensional martingale $M$, and the function $g$ satisfies
    \[ \left\vert g(t,M_{t})\right\vert \leq c_{1}\left\Vert M_{t}\right\Vert_{\infty}^{q}+c_{2}\,, \quad 0\leq t\leq T\,, \quad q\geq1\,, \]
    where $\left\Vert \cdot\right\Vert_{\infty}$ denotes the max-norm in $\mathbb{R}^{d}$, it follows that
    \begin{align*}
        \mathbb{E}\left[\sup_{0\leq t\leq T}Z_{t}^{\alpha}\right] & \leq\mathbb{E}\left[\sup_{0\leq t\leq T}\left(c_{1}\left\Vert M_{t}\right\Vert_{\infty}^{q}+c_{2}\right)^{\alpha}\right] \leq \mathbb{E}\left[\left(c_{1}\sup_{0\leq t\leq T}\left\Vert M_{t}\right\Vert_{\infty}^{q}+c_{2}\right)^{\alpha}\right]\, \\
        & \leq2^{\alpha-1}c_{1}^{\alpha}\mathbb{E}\left[\sup_{0\leq t \leq T}\left\Vert M_{t}\right\Vert_{\infty}^{\alpha q}\right] +2^{\alpha-1}c_{2}^{\alpha}\, \\
        \text{(with }M & =(M^{i})_{i=1,\ldots,d}\text{)}\leq2^{\alpha-1} c_{1}^{\alpha}d^{\alpha q-1}\sum_{i=1}^{d}\mathbb{E}\left[\sup_{0\leq t\leq T}\left\vert M_{t}^{i}\right\vert ^{\alpha q}\right] +2^{\alpha-1}c_{2}^{\alpha}\, \\
        & \leq 2^{\alpha-1}c_{1}^{\alpha}d^{\alpha q-1}C_{\alpha q}^{BDG}\sum_{i=1}^{d}\mathbb{E}\left[\langle M_{T}^{i},M_{T}^{i}\rangle^{\alpha q/2}\right] +2^{\alpha-1}c_{2}^{\alpha}\,,
    \end{align*}
    using the Burkholder-Davis-Gundy inequality. 
    The constant $C_{\alpha q}^{BDG}$ is universal and explicitly known for continuous martingales and jump martingales, respectively. 
    Furthermore, if $\alpha q=2$ for example, $\mathbb{E}\left[\langle M_{T}^{i}, M_{T}^{i}\rangle\right] $ may usually be estimated explicitly for a given specific model.
\end{remark}

\begin{lemma}
\label{lem:frac-m} Let
\[
\mathfrak{m}_{i}:=\min\left\{  j\geq i:\tau_{j}>T\right\}  \text{, \ for
\ }i\in\mathbb{N}_{0}.\,
\]
Then, 
\begin{itemize}\item[(i)] for all $i\in\mathbb{N}_{0}$, we have that $U_{j}=0$ a.s.\ for all
$j\geq\mathfrak{m}_{i}$, 
\item[(ii)] $\,i\leq\mathfrak{m}_{i}<\infty$ a.s., and
\item[(iii)] $\mathfrak{m}_{i}$ is a $\mathbb{G}$-stopping time.
\end{itemize}
\end{lemma}

\begin{proof}
(i) and (ii) follow from the assumption that $Z_{t}=0$ a.s.\ for
all $t>T$ and assumption \eqref{eq:boundZ}. \\
\indent(iii): Consider, for all $n\in\mathbb{N}_{0},$
\[
\left\{  \mathfrak{m}_{i}=n\right\}  =\left\{  \tau_{n}>T\right\}  \cap
\bigcap\limits_{i\leq m<n}\left\{  \tau_{m}\leq T\right\}  \,,
\]
and so, for all $s\geq0$,
\begin{align*}
&  \left\{  \mathfrak{m}_{i}=n\right\}  \cap\left\{  \tau_{n}\leq s\right\}
=\left\{  \tau_{n}>T\right\}  \cap\bigcap\limits_{i\leq m<n}\left\{  \tau
_{m}\leq T\right\}  \cap\left\{  \tau_{n}\leq s\right\}  \,\\
&  =%
\begin{cases}
\varnothing\in\mathcal{F}_{s}, & s<T\\
\left\{  \tau_{n}>T\right\}  \cap\bigcap\limits_{i\leq m<n}\left\{  \tau
_{m}\leq T\right\}  \cap\left\{  \tau_{n}\leq s\right\}  \in\mathcal{F}_{s}, &
s\geq T
\end{cases}
\,.
\end{align*}
That is, $\left\{  \mathfrak{m}_{i}=n\right\}  \in\mathcal{F}_{\tau_{n}%
}=\mathcal{G}_{n}$.
\end{proof}

\vskip 0.2cm
\begin{proof}
[Proof of Proposition~\ref{prop:bellman}]We note that \eqref{eq:boundZ1},
which follows from assumption \eqref{eq:boundZ}, implies that condition
\cite[Ch.~1, (1.1.3))]{PS06} is satisfied. Following \cite[Ch.~1]{PS06}, we
then have~(i) and~(ii), and that \eqref{iii1} is an optimal stopping time that
satisfies \eqref{iii2}.  We next
prove for the infinite horizon case $K=\infty$ that   $\mathfrak{n}_{i}^{\infty,\ast}<\infty$ a.s. 
Consider the event
\begin{align*}
A_{i} &  :=\left\{  \mathfrak{n}_{i}^{\infty,\ast}=\infty\right\}  =\left\{
\forall_{j\geq i}:U_{j}<\mathbb{E}_{\mathcal{G}_{j}}\left[  \esssup_{j+1\leq
\mathfrak{n<\infty}}\mathbb{E}_{\mathcal{G}_{j+1}}\left[  U_{\mathfrak{n}%
}\right]  \right]  \right\}  \\
&  \subset\left\{  \forall_{j\geq i}:U_{j}<\mathbb{E}_{\mathcal{G}_{j}}\left[
\sup_{j+1\leq r\mathfrak{<\infty}}U_{r}\right]  \right\}  =:\widetilde{A}_{i},
\end{align*}
for some $i\in\mathbb{N}_{0}$. 
Then, with $\mathfrak{m}_{i}$ defined in Lemma~\ref{lem:frac-m}, we have on the set  $\widetilde{A}_{i}$,
\begin{align*}
0=U_{\mathfrak{m}_{i}}=\sum_{l=i}^{\infty}1_{\left\{  \mathfrak{m}%
_{i}=l\right\}  }U_{l} &  <\sum_{l=i}^{\infty}1_{\left\{  \mathfrak{m}%
_{i}=l\right\}  }\mathbb{E}_{\mathcal{G}_{l}}\left[  \sup_{l+1\leq
r\mathfrak{<\infty}}U_{r}\right]  \\
&  =\sum_{l=i}^{\infty}\mathbb{E}_{\mathcal{G}_{l}}\left[  1_{\left\{
\mathfrak{m}_{i}=l\right\}  }\sup_{l+1\leq r\mathfrak{<\infty}}U_{r}\right]
=0\text{ \ \ a.s.}\,,
\end{align*}
using that $\mathfrak{m}_{i}$ is a stopping time due to Lemma~\ref{lem:frac-m}. 
This implies that $\mathbb{P}\left(  A_{i}\right)
\leq\mathbb{P}\left(  \widetilde{A}_{i}\right)  =0$. 
Hence, it follows that
$\mathfrak{n}_{i}^{\infty,\ast}<\infty$ a.s.
\end{proof}

\vskip 0.2cm
\begin{proof}[Proof of Corollary~\ref{cor:opstoptimes}]
    Indeed, from item~(i) in the proof of Proposition~\ref{prop:proprandtime}, we see that $\mathfrak{t}^{\ast}$ is an $\mathbb{F}$-stopping time that satisfies $\mathfrak{t}^{\ast}\in\mathcal{T}_{K}$ and that moreover satisfies
    \[ Y_{0}^{(K)}=\sup_{\mathfrak{n}\in\mathcal{N},~\mathfrak{n}\leq K}\mathbb{E}\left[U_{\mathfrak{n}}\right] =\mathbb{E}\left[U_{\mathfrak{n}^{\ast}}\right] = \mathbb{E}\left[Z_{\mathfrak{t}^{\ast}}\right]\,. \]
\end{proof}
\subsection{Duality}\label{app:duality}
\begin{proof}[Proof of Proposition \ref{dualprop}]
    By Proposition~\ref{prop:proprandtime} and Doob's sampling theorem we have, for any $M\in\mathcal{M}_{0}^{\mathsf{UI}}$,
    \begin{equation}
        Y_{0} =\sup_{\mathfrak{n}\in\mathcal{N}}\mathbb{E}\left[U_{\mathfrak{n}}\right] = \sup_{\mathfrak{n}\in\mathcal{N}}\mathbb{E}\left[U_{\mathfrak{n}}-M_{\mathfrak{n}}\right] \leq \mathbb{E}\left[\sup_{i\geq0}\left(U_{i}-M_{i}\right)\right]\,. \label{eq:weak}
    \end{equation}
    Note that for any $i \geq 0$ (empty sums being zero),
    \begin{align}
        U_{i}-M_{i}^{\circ} & =U_{i}-\sum_{j=1}^{i}Y_{j}+\sum_{j=1}^{i}\mathbb{E}_{\mathcal{G}_{j-1}}\left[{Y}_{j}\right]\, \label{eq:y*0} \\
        & =Y_{0}+U_{i}-Y_{i}-\sum_{j=0}^{i-1}Y_{j}+\sum_{j=0}^{i-1}\mathbb{E}_{\mathcal{G}_{j}}\left[Y_{j+1}\right]\, \nonumber \\
        & \leq Y_{0}\,,\nonumber
    \end{align}
    due to \eqref{eq:YKsupermart}. Hence, we have
    \begin{equation}
        \sup_{i\geq0}\left(U_{i}-M_{i}^{\circ}\right)\leq Y_{0},\label{eq:y*1}
    \end{equation}
    by \eqref{eq:y*0}. 
    We now only need to show that $M^{\circ}\in\mathcal{M}_{0}^{\mathsf{UI}}$. Indeed, if $M^{\circ}\in\mathcal{M}_{0}^{\mathsf{UI}}$, we have
    \[ Y_{0}\leq\mathbb{E}\left[\sup_{i\geq0}\left(U_{i}-M_{i}^{\circ}\right)\right], \]
    by inequality \eqref{eq:weak}, and then by \eqref{eq:y*1} we obtain statement~(ii) due to the sandwich property. Finally, (ii) combined with \eqref{eq:weak} yields statement~(i). Let us show that $M^{\circ}\in\mathcal{M}_{0}^{\mathsf{UI}}$. Note that, by \eqref{Mcirc},
    \begin{equation}
        \left\vert M_{i}^{\circ}\right\vert \leq\sum_{j=1}^{\infty}\left\vert Y_{j}-\mathbb{E}_{\mathcal{G}_{j-1}}\left[Y_{j}\right] \right\vert\quad\text{a.s.\ for all }i \geq 0\,.\label{Mcircmaj}
    \end{equation}
    Since $Z_{t} \geq 0$ for $t \geq 0$, and $Z_{t} = 0$ for $t>T$, we have that $Y_{j}=Y_{j}1_{\left\{\tau_{j-1}\leq T\right\}}\geq0$ a.s. Hence, by the monotone convergence theorem,
    \begin{align*}
        \mathbb{E}\left[\sum_{j=1}^{\infty}\left\vert Y_{j}-\mathbb{E}_{\mathcal{G}_{j-1}}\left[Y_{j}\right] \right\vert \right] &= \mathbb{E}\left[\sum_{j=1}^{\infty}\left\vert Y_{j}-\mathbb{E}_{\mathcal{G}_{j-1}}\left[Y_{j}\right] \right\vert 1_{\left\{\tau_{j-1}\leq T\right\}}\right]\, \\
        & \leq\mathbb{E}\left[\sum_{j=1}^{\infty}Y_{j}1_{\left\{\tau_{j-1}\leq T\right\}}+\mathbb{E}_{\mathcal{G}_{j-1}}\left[Y_{j}1_{\left\{\tau_{j-1}\leq T\right\}}\right]\right]\, \\
        &= 2\mathbb{E}\left[\sum_{j=1}^{\infty}Y_{j}1_{\left\{\tau_{j-1}\leq T\right\}}\right] \leq 2\mathbb{E}\left[\sup_{i\in\mathbb{N}_{0}}U_i\sum_{j=1}^{\infty}1_{\left\{\tau_{j-1}\leq T\right\}}\right]\, \\
        &= 2\mathbb{E}\left[N_{T} \sup_{i\in\mathbb{N}_{0}}U_i\right]<\infty,
    \end{align*}
    due to \eqref{UIcon} and
    using that $\left\{  \tau_{j-1}\leq T\right\}  =\left\{  \mathfrak{m}_{j}\leq
j-1\right\}  ^{c}$ is $\mathcal{G}_{j-1}$-measurable since $\mathfrak{m}_{j}$
is a stopping time by  Lemma\thinspace
\ref{lem:frac-m}.
    So the family $\left(M_{i}^{\circ}\right)_{i\geq0}$ is uniformly integrable since it has an integrable dominator \eqref{Mcircmaj}.
\end{proof}

\begin{corollary}\label{cor:holder-UI}
    Condition~\eqref{eq:boundZ1} implies that $\sup_{i\in\mathbb{N}_0}U_i \in L^\alpha$. By Hölder's inequality,
    \[ \mathbb{E}\left[N_{T}\sup_{i\in\mathbb{N}_{0}}U_i\right] \leq B\mathbb{E}\left[N_{T}^{\alpha'}\right]^{1/\alpha'}\,, \]
    where $\alpha':=\alpha/(\alpha-1)>1$ is the Hölder conjugate of $\alpha$. So a sufficient condition for \eqref{UIcon} is that $N_T \in L^{\alpha'}$ i.e., 
    \begin{equation}
        \mathbb{E}\left[N_{T}^{\frac{\alpha}{\alpha-1}}\right] < \infty\,. \label{corcon}
    \end{equation}
\end{corollary}

The following lemma provides a convenient way to check (\ref{corcon}) by checking the tail probabilities of $N_T$. 
\begin{lemma}\label{lem:check}
    Let $N$ be an $\mathbb{N}$-valued random variable and let 
    $\alpha > 1$. 
    Then, 
    \[ N \in L^{\frac{\alpha}{\alpha-1}} < \infty \text{ if and only if } \sum_{j=1}^\infty j^{\frac{1}{\alpha-1}}\mathbb{P}(N > j) < \infty\,. \]
\end{lemma}
\begin{proof}
     We may write (since $N\geq1$ a.s.),
    \begin{align*}
        \mathbb{E}\left[N^{\frac{\alpha}{\alpha-1}}\right] &= \sum_{j=1}^\infty j^{\frac{\alpha}{\alpha - 1}}\mathbb{P}(N = j)\, \\
        &= \sum_{j=1}^\infty j^{\frac{\alpha}{\alpha - 1}}\left(\mathbb{P}(N > j - 1) - \mathbb{P}(N > j)\right)\, \\
        &=1+ \sum_{j=1}^\infty \left((j+1)^{\frac{\alpha}{\alpha-1}}-j^{\frac{\alpha}{\alpha-1}}\right)\mathbb{P}(N>j)\,.
    \end{align*} 
Let us now write, with $\alpha^{\prime}=\alpha/(\alpha-1),$
\[
\left(  j+1\right)  ^{\frac{\alpha}{\alpha-1}}-j^{\frac{\alpha}{\alpha-1}%
}=j^{\alpha^{\prime}}\left(  \left(  1+\frac{1}{j}\right)  ^{\alpha^{\prime}%
}-1\right)  .
\]
Then, due to elementary inequalities, we have that
\[
\frac{\alpha^{\prime}}{j}\leq\left(  1+\frac{1}{j}\right)  ^{\alpha^{\prime}%
}-1\leq e^{\frac{\alpha^{\prime}}{j}}-1\leq\frac{2\alpha^{\prime}}{j}\text{
\ \ for }j>\alpha^{\prime}.%
\]
Hence, $\left(  j+1\right)  ^{\frac{\alpha}{\alpha-1}}-j^{\frac{\alpha}%
{\alpha-1}}\asymp j^{\alpha^{\prime}-1}=j^{\frac{1}{\alpha-1}},$
    which shows the equivalence. 
\end{proof}
\subsection{Policy iteration}\label{app:policy-iteration}
\begin{proof}[Proof of Proposition~\ref{polit}]
    Consider for $N\in\mathbb{N}$ the truncated reward $U_{i}^{N}:=U_{i}1_{\{i< N\}}$ and the stopped policy $\sigma_{i}^{N}:=\sigma_{i}\wedge N$.
    It is easy to see that, for $0\leq i\leq N$,
    \[ i\leq\sigma_{i}^{N}\leq N\,\quad \text{and} \quad \sigma_{i}^{N}>i \quad \Longrightarrow \quad \sigma_{i}^{N}=\sigma_{i+1}^{N}\,, \]
    i.e., $\sigma_{i}^{N}$ is consistent in the sense of \eqref{con}. We next consider the iteration step
    \[ \widehat{\sigma}_{i}^{N}:=\inf\left\{  j\geq i:U_{j}^{N}\geq\max_{j\leq k\leq j+\kappa}\mathbb{E}_{\mathcal{G}_{j}}\left[  U_{\sigma_{k}^{N}}^{N}\right]\right\} \]
    and define
    \begin{align*}
        Y_{i}^{\sigma^{N}} &:= \mathbb{E}_{\mathcal{G}_{i}}\left[U_{\sigma_{i}^{N}}^{N}\right]\,, \\
        \widetilde{Y}_{i}^{\sigma^{N}} &:= \max_{i\leq k\leq i+\kappa}\mathbb{E}_{\mathcal{G}_{i}}\left[  U_{\sigma_{k}^{N}}^{N}\right]\,, \\
        Y_{i}^{\widehat{\sigma}^{N}} &:= \mathbb{E}_{\mathcal{G}_{i}}\left[U_{\widehat{\sigma}_{i}^{N}}^{N}\right]\,.
    \end{align*}
    Note that $i\leq\widehat{\sigma}_{i}^{N}\leq N$. Due to \cite[Theorem~3.1]{KS2006}, it holds that
    \begin{equation}
        Y_{i}^{\sigma^{N}}\leq\widetilde{Y}_{i}^{\sigma^{N}} \leq Y_{i}^{\widehat{\sigma}^{N}}\leq Y_{i}^{(N)}\,, \quad \text{for} \quad\leq i\leq N\,, \label{KSit}
    \end{equation}
    where $\left(Y_{i}^{(N)}\right)$ is the Snell envelope for the truncated reward $\left(U_{i}^{N}\right)$. We now fix $i$ and show that for $N\rightarrow\infty$, a.s.,

    (i) $Y_{i}^{\sigma^{N}}\rightarrow Y_{i}^{\sigma}$,

    (ii) $\widetilde{Y}_{i}^{\sigma^{N}}\rightarrow\widetilde{Y}_{i}^{\sigma}$,

    (iii) $Y_{i}^{\widehat{\sigma}^{N}}\rightarrow Y_{i}^{\widehat{\sigma}}$,\\
    \noindent from which the statement follows by \eqref{KSit} and the obvious fact that $Y_{i}^{(N)}\leq Y_{i}$.

    Ad (i)+(ii): For any fixed $k$ one has $\sigma_{k}^{N}\rightarrow\sigma_{k}<\infty$ a.s. for $N\rightarrow\infty$, hence
    \[ U_{\sigma_{k}^{N}}^{N}=U_{\sigma_{k}^{N}}1_{\{\sigma_{k}^{N}< N\}}=U_{\sigma_{k}}1_{\{\sigma_{k}< N\}}\rightarrow U_{\sigma_{k}}\,,~\text{\ a.s.} \]
    It then follows by conditional dominated convergence that
    \begin{align*}
        Y_{i}^{\sigma^{N}} &= \mathbb{E}_{\mathcal{G}_{i}}\left[U_{\sigma_{i}^{N}}^{N}\right]  \rightarrow\mathbb{E}_{\mathcal{G}_{i}}\left[U_{\sigma_{k}}\right] = Y_{i}^{\sigma}\,,\quad \text{and} \\
        \widetilde{Y}_{i}^{\sigma^{N}} &= \max_{i\leq k\leq i+\kappa}\mathbb{E}_{\mathcal{G}_{i}}\left[U_{\sigma_{k}^{N}}^{N}\right]\rightarrow\max_{i\leq k\leq i+\kappa}\mathbb{E}_{\mathcal{G}_{i}}\left[U_{\sigma_{k}}\right] = \widetilde{Y}_{i}^{\sigma}\,.
    \end{align*}

    Ad (iii): Note that by the definition of $U_{i}^{N}$,
    \[ \widehat{\sigma}_{i}^{N} := \inf\left\{j \geq i:U_{j}1_{\{j\leq N\}}\geq\max_{j \leq k \leq j+\kappa}\mathbb{E}_{\mathcal{G}_{j}}\left[U_{\sigma_{k}}1_{\{\sigma_{k}< N\}}\right]\right\}\,. \]
    Let us suppose that we are on the set $\left\{\omega:\widehat{\sigma}_{i}\left(\omega\right)=l\right\}$.
    For $N > l$, one then has
    \[ U_{l}\left(\omega\right)1_{\{l< N\}}=U_{l}\left(\omega\right) \geq \max_{l\leq k\leq l+\kappa}\mathbb{E}_{\mathcal{G}_{l}}\left[  U_{\sigma_{k}}\right]\left(\omega\right)\geq\max_{l\leq k\leq l+\kappa}\mathbb{E}_{\mathcal{G}_{l}}\left[U_{\sigma_{k}}1_{\{\sigma_{k}< N\}}\right]\left(\omega\right) \]
    (recall that by assumption $U_\cdot\geq0$). Hence, for $N\geq l$, one has $\widehat{\sigma}_{i}^{N}\left(\omega\right)\leq l$. Since $l$ was arbitrary, we conclude that
    \[ \widehat{\sigma}_{i}^{N}\leq\widehat{\sigma}_{i}~\text{a.s.} \quad \text{for}~N \geq \mathcal{N}\left(\omega\right)\quad \text{say}\,. \]
    Now assume that $\widehat{\sigma}_{i}^{N}\nrightarrow\widehat{\sigma}_{i}$ with positive probability. Then with positive probability, there must exist numbers $r$ and $l$ with $0\leq r<l$ and a sequence $\left(N_{m}\right)_{m\in\mathbb{N}}$ with $N_{m}\uparrow\infty$ for $m\uparrow\infty$, such that $\widehat{\sigma}_{i}=l$ and $\widehat{\sigma}_{i}^{N_{m}}=r<l$ for $m\in\mathbb{N}$\,. This implies
    \[ U_{r}1_{\{r\leq N_{m}\}}\geq\max_{r\leq k\leq r+\kappa}\mathbb{E}_{\mathcal{G}_{r}}\left[U_{\sigma_{k}}1_{\{\sigma_{k}< N_{m}\}}\right]\,,\quad \forall m\in\mathbb{N}\,, \]
    on a set of positive probability. By conditional dominated convergence it then follows that
    \[ U_{r}\geq\max_{r\leq k\leq r+\kappa}\mathbb{E}_{\mathcal{G}_{r}}\left[U_{\sigma_{k}}\right] \]
    on this set, which contradicts the definition of $\widehat{\sigma}_{i}$. That is, $\widehat{\sigma}_{i}^{N}\rightarrow\widehat{\sigma}_{i}$ a.s., and then (iii) follows by conditional dominated convergence.
\end{proof}

\begin{proof}[Proof of Corollary \ref{corimp}]
    Ad (i): By Proposition~\ref{polit} we may write
    \begin{align*}
        1_{\left\{Y_{i}^{\widehat{\sigma}}<U_{i}\right\}}Y_{i}^{\widehat{\sigma}} &= 1_{\left\{  Y_{i}^{\widehat{\sigma}}<U_{i}\right\}}1_{\left\{
        \widetilde{Y}_{i}^{\sigma}\leq U_{i}\right\}}\mathbb{E}_{\mathcal{G}_{i}}\left[U_{\widehat{\sigma}_{i}}\right]\, \\
        \text{(definition of }\widehat{\sigma}_{i}\text{)} &= 1_{\left\{Y_{i}^{\widehat{\sigma}}<U_{i}\right\}}1_{\left\{\widetilde{Y}_{i}^{\sigma}\leq U_{i}\right\}  }U_{i}=1_{\left\{  Y_{i}^{\widehat{\sigma}}<U_{i}\right\}}U_{i}\,, \quad \text{a.s.}
    \end{align*}
    That is, on the set $\left\{Y_{i}^{\widehat{\sigma}}<U_{i}\right\}$ one has $Y_{i}^{\widehat{\sigma}}= U_{i}$, a.s., which implies $\mathbb{P}(Y_{i}^{\widehat{\sigma}}<U_{i})=0$.

    Ad (ii): First note that $\widehat{\sigma_i}^{\prime}\geq\widehat{\sigma_i}$ by construction. 
    Hence, we may write
    \begin{align*}
        Y_{i}^{\widehat{\sigma}^{\prime}} &= \mathbb{E}_{\mathcal{G}_{i}}\left[U_{\widehat{\sigma}_{i}^{\prime}}\right] = \sum_{m=i}^{\infty}\mathbb{E}_{\mathcal{G}_{i}}\left[  1_{\left\{  \widehat{\sigma}_{i}=m\right\}}U_{\widehat{\sigma}_{i}^{\prime}}\right]\, \\
        \text{(using \eqref{con})} &= \sum_{m=i}^{\infty}\mathbb{E}_{\mathcal{G}_{i}}\left[1_{\left\{\widehat{\sigma}_{i}=m\right\}}U_{\widehat{\sigma}_{m}^{\prime}}\right]  = \sum_{m=i}^{\infty}\mathbb{E}_{\mathcal{G}_{i}}\left[1_{\left\{  \widehat{\sigma}_{i}=m\right\}}Y_{m}^{\widehat{\sigma}^{\prime}}\right]\, \\
        & \geq\sum_{m=i}^{\infty}\mathbb{E}_{\mathcal{G}_{i}}\left[1_{\left\{  \widehat{\sigma}_{i}=m\right\}}U_{m}\right] = \mathbb{E}_{\mathcal{G}_{i}}\left[U_{\widehat{\sigma}_{i}}\right]
        = Y_{i}^{\widehat{\sigma}}\,,
    \end{align*}
    where the already proved statement~(i) is used.
\end{proof}

\vskip 0.2cm\begin{proof}[Proof of Theorem~\ref{impthm}]
    The fact that $\sigma_{i}^{\ast}<\infty$ follows from Proposition~\ref{prop:bellman}. Suppose that $\sigma_{i}^{(m)}>\sigma_{i}^{\ast}$ for some $m\in\mathbb{N}_{0}$. 
    Then we have
    \begin{align*}
        U_{\sigma_{i}^{\ast}} &< \max_{\sigma_{i}^{\ast}\leq k\leq\sigma_{i}^{\ast}+\kappa}\mathbb{E}_{\mathcal{G}_{\sigma_{i}^{\ast}}}\left[U_{\sigma_{k}^{(m)}}\right] =\sum_{j=i}^{\infty}1_{\left\{  \sigma_{i}^{\ast}=j\right\}}\max_{j\leq k\leq j+\kappa}\mathbb{E}_{\mathcal{G}_{j}}\left[U_{\sigma_{k}^{(m)}}\right]\, \\
        \text{(by Proposition~\ref{polit})} &\leq \sum_{j=i}^{\infty}1_{\left\{
        \sigma_{i}^{\ast}=j\right\}}Y_{j}=Y_{\sigma_{i}^{\ast}}\,,
    \end{align*}
    which contradicts the optimality of $\sigma_{i}^{\ast}.$ Suppose next that
    $\sigma_{i}^{(m+1)}<\sigma_{i}^{(m)}$ for some $m\in\mathbb{N}_{0}$. Then, since $\sigma_{i}^{(0)}=i$, we must have that $m \geq 1$, and due to the definition of $\sigma_{i}^{(m)}$ that
    \begin{align*}
        U_{\sigma_{i}^{(m+1)}} & <\max_{\sigma_{i}^{(m+1)}\leq k\leq\sigma_{i}^{(m+1)}+\kappa}\mathbb{E}_{\mathcal{G}_{\sigma_{i}^{(m+1)}}}\left[U_{\sigma_{k}^{(m-1)}}\right]\, \\
        & =\sum_{j=i}^{\infty}1_{\left\{\sigma_{i}^{(m+1)}=j\right\}}\max_{ j\leq k\leq j+\kappa}\mathbb{E}_{\mathcal{G}_{j}}\left[  U_{\sigma_{k}^{(m-1)}}\right]\,.
    \end{align*}
    On the other hand, due to the definition of $\sigma_{i}^{(m+1)}$ we must have that
    \begin{align*}
        U_{\sigma_{i}^{(m+1)}} &\geq \max_{\sigma_{i}^{(m+1)}\leq k\leq\sigma_{i}^{(m+1)}+\kappa}\mathbb{E}_{\mathcal{G}_{\sigma_{i}^{(m+1)}}}\left[U_{\sigma_{k}^{(m)}}\right]\, \\
        & =\sum_{j=i}^{\infty}1_{\left\{\sigma_{i}^{(m+1)}=j\right\}}\max_{j \leq k \leq j+\kappa}\mathbb{E}_{\mathcal{G}_{j}}\left[  U_{\sigma_{k}^{(m)}}\right]\,,
    \end{align*}
    which in turn yields
    \[ \sum_{j=i}^{\infty}1_{\left\{  \sigma_{i}^{(m+1)}=j\right\}  }\max_{j\leq k\leq j+\kappa}\mathbb{E}_{\mathcal{G}_{j}}\left[  U_{\sigma_{k}^{(m)}}\right] <\sum_{j=i}^{\infty}1_{\left\{  \sigma_{i}^{(m+1)}=j\right\}}\max_{j\leq k\leq j+\kappa}\mathbb{E}_{\mathcal{G}_{j}}\left[U_{\sigma_{k}^{(m-1)}}\right]\,. \]
    However, since for all $j\geq i$,
    \[ \max_{j\leq k\leq j+\kappa}\mathbb{E}_{\mathcal{G}_{j}}\left[  U_{\sigma_{k}^{(m)}}\right] \geq \mathbb{E}_{\mathcal{G}_{j}}\left[  U_{\sigma_{j}^{(m)}}\right] = Y_{j}^{(m)}\geq\max_{j\leq k\leq j+\kappa}\mathbb{E}_{\mathcal{G}_{j}}\left[U_{\sigma_{k}^{(m-1)}}\right] \quad \text{a.s.}\,, \]
    by Proposition \ref{polit}, we get a contradiction. 
    Thus, \eqref{siginc} is proved.

    Due to \eqref{siginc} and Proposition~\ref{polit} we may define
    \[ \sigma_{i}^{\infty}=\uparrow\lim_{m\rightarrow\infty}\sigma_{i}^{(m)} \quad \text{and} \quad Y_{i}^{\infty} = \uparrow\lim_{m\rightarrow\infty}Y_{i}^{(m)} \quad \text{a.s.} \]
    Note that $\sigma_{i}^{\infty}\leq\sigma_{i}^{\ast}<\infty$, and that, since $\sigma_{i}^{(m)}$ is integer valued, one must have that $\sigma_{i}^{(m)}\left(\omega\right) = \sigma_{i}^{(\infty)}\left(  \omega\right)$ for $m>N\left(i,\omega\right)$ say. 
    Hence, $U_{\sigma_{i}^{(m)}}\rightarrow U_{\sigma_{i}^{\infty}}$ a.s., and then by dominated convergence we get
    \[ Y_{i}^{\infty}=\uparrow\lim_{m\rightarrow\infty}Y_{i}^{(m)}=\uparrow \lim_{m\rightarrow\infty}\mathbb{E}_{\mathcal{G}_{i}}\left[  U_{\sigma_{i}^{(m)}}\right] = \mathbb{E}_{\mathcal{G}_{i}}\left[  U_{\sigma_{i}^{\infty}}\right]\,. \]
    Obviously, we have that $Y_i\geq Y_{i}^{\infty}\geq U_{i}$ a.s. Since $\left(Y_{i}\right)$ is the smallest supermartingale that dominates the reward $\left(U_{i}\right)$ by Proposition~\ref{prop:bellman}, we can prove that $\sigma_{i}^{\infty}=\sigma_{i}^{\ast}$ and $Y_{i}^{\infty}=Y_{i}$, and hence \eqref{pollim}, by showing that $\left(Y_{i}^{\infty}\right)$ is a supermartingale. It is easy to see that $\left(\sigma_{i}^{\infty}\right) $ is a consistent $\mathcal{G}$-stopping family, i.e., $\left(  \sigma_{i}^{\infty}\right)$ satisfies \eqref{con}. So for any $i \geq 0$ one has
    \begin{align}
        \mathbb{E}_{\mathcal{G}_{i}}\left[  Y_{i+1}^{\infty}\right] &= \mathbb{E}_{\mathcal{G}_{i}}\left[  U_{\sigma_{i+1}^{\infty}}\right] = Y_{i}^{\infty}+\mathbb{E}_{\mathcal{G}_{i}}\left[\left(  U_{\sigma_{i+1}^{\infty}}-U_{\sigma_{i}^{\infty}}\right)\right]\, \nonumber \\
        &= Y_{i}^{\infty}+\mathbb{E}_{\mathcal{G}_{i}}\left[1_{\left\{\sigma_{i}^{\infty}=i\right\}}\left(U_{\sigma_{i+1}^{\infty}}-U_{\sigma_{i}^{\infty}}\right)\right]\, \nonumber \\
        \text{(by \eqref{con})} &= Y_{i}^{\infty}+1_{\left\{  \sigma_{i}^{\infty}=i\right\}}\mathbb{E}_{\mathcal{G}_{i}}\left[  U_{\sigma_{i+1}^{\infty}}\right] - 1_{\left\{\sigma_{i}^{\infty}=i\right\}}U_{i}\,. \label{supm}
    \end{align}
    Note that by \eqref{siginc} and the definition of $\sigma_{i}^{\infty}$ one has $\left\{\sigma_{i}^{\infty} = i\right\} = \bigcap\limits_{m=0}^{\infty}\left\{\sigma_{i}^{(m)}=i\right\}$. 
    Hence, by the definition of $\sigma_{i}^{(m^{\prime})}$, for all $m^{\prime} \geq 1$,
    \begin{align*}
        1_{\left\{\sigma_{i}^{\infty}=i\right\}}U_{i} &\geq 1_{\left\{\sigma_{i}^{\infty}=i\right\}}\max_{i\leq k\leq i+\kappa}\mathbb{E}_{\mathcal{G}_{i}}\left[U_{\sigma_{k}^{(m^{\prime}-1)}}\right]\, \\
        &\geq 1_{\left\{\sigma_{i}^{\infty}=i\right\}}\mathbb{E}_{\mathcal{G}_{i}}\left[  U_{\sigma_{i}^{(m^{\prime}-1)}}\right]\,.
    \end{align*}
    Since $\uparrow\lim_{m^{\prime}\uparrow\infty}\mathbb{E}_{\mathcal{G}_{i}}\left[  U_{\sigma_{i}^{(m^{\prime}-1)}}\right]  =\mathbb{E}_{\mathcal{G}_{i}}\left[  U_{\sigma_{i}^{\infty}}\right]$ a.s., it then follows that
    \[ 1_{\left\{\sigma_{i}^{\infty}=i\right\}}U_{i}\geq1_{\left\{\sigma_{i}^{\infty}=i\right\}}\mathbb{E}_{\mathcal{G}_{i}}\left[  U_{\sigma_{i}^{\infty}}\right]\,, \]
    hence by \eqref{supm} we get $\mathbb{E}_{\mathcal{G}_{i}}\left[Y_{i+1}^{\infty}\right] \leq Y_{i}^{\infty}$. That is, $\left(Y_{i}^{\infty}\right)$ is a super-martingale.
\end{proof}
\subsection{Finite-horizon approximation}\label{app:FHA}
\begin{proof}[Proof of Lemma \ref{lem:FHA}]
    The first inequality is trivial. 
    To prove the second inequality, we write:
    \begin{align*}
        Y_0^{(\infty)} &= \sup_{\mathfrak{n} \in \mathcal{N}}\mathbb{E}\left[U_{\mathfrak{n}}1_{\{\mathfrak{n}\leq K\}} + U_{\mathfrak{n}}1_{\{\mathfrak{n} > K\}}\right]\, \\
        &\leq \sup_{\mathfrak{n} \in \mathcal{N}}\mathbb{E}\left[U_{\mathfrak{n}\wedge K}1_{\{\mathfrak{n}\leq K\}}\right] + \sup_{\mathfrak{n} \in \mathcal{N}}\mathbb{E}\left[U_{\mathfrak{n}\vee K}1_{\{\mathfrak{n} > K\}}\right]\, \\
        &\leq \sup_{\mathfrak{n} \in \mathcal{N}, 0\leq \mathfrak{n} \leq K}\mathbb{E}\left[U_{\mathfrak{n}}\right] + \sup_{\mathfrak{n} \in \mathcal{N}, \mathfrak{n} > K}\mathbb{E}\left[U_{\mathfrak{n}}\right] \leq Y_0^{(K)} + \sup_{\mathfrak{n} \in \mathcal{N}, \mathfrak{n} > K}\mathbb{E}\left[U_{\mathfrak{n}}\right]\,,
    \end{align*}
    by Proposition~\ref{prop:proprandtime}. 
    Now, for any $\mathfrak{n}$ with $\mathfrak{n} > K$, we have
        \begin{align*}
        \mathbb{E}\left[U_{\mathfrak{n}}\right] & = \mathbb{E}\left[Z_{\tau_{\mathfrak{n}}}1_{\{\tau_{\mathfrak{n}} \leq T\}}\right] + \mathbb{E}\left[Z_{\tau_{\mathfrak{n}}}1_{\{\tau_{\mathfrak{n}} > T\}}\right]\, \\
        & \leq \mathbb{E}\left[Z_{\tau_{\mathfrak{n}}}1_{\{\tau_{K} \leq T\}}\right] + 0\, \\
        & \leq \mathbb{E}\left[1_{\left\{\tau_{K}\leq T\right\}} \sup_{i\in\mathbb{N}_{0},~i>K}Z_{\tau_{i}}\right]\,\\
        &\leq B\mathbb{P}\left(\tau_K \leq T\right)^{1-1/\alpha}\,,
        \end{align*}  
        by a H\"older estimate using \eqref{eq:boundZ}.
         Because, by \eqref{eq:stop}, $\tau_k \to \infty$ a.s.\ as $k\to \infty$, dominated convergence implies that $\mathbb{P}(\tau_K\leq T) \to 0$ as $K\to\infty$, $K\in\mathbb{N}$. 
         Hence, to bound $\sup_{\mathfrak{n} \in \mathcal{N}, \mathfrak{n} > K}\mathbb{E}\left[U_{\mathfrak{n}}\right]$ by $\varepsilon$, it suffices to pick a value of $K$ such that $\mathbb{P}(\tau_K \leq T) < (\varepsilon/B)^{\alpha/(\alpha-1)}$, with $B$ and $\alpha$ specified in \eqref{eq:boundZ}. 
\end{proof}
\subsection{Markovian structure}\label{app: markovian}
The following results are readily obtained.
\begin{lemma}\label{lemSMP}
    If Assumption~\ref{SA} is fulfilled, then for any non-negative Borel measurable $g:$ $(t,x,\theta) \in \mathbb{R}\times\mathbb{R}^{d}\times\mathbb{R}^{q} \rightarrow g(t,x,\theta)$, we have that
    \begin{align*}
        \mathbb{E}_{\mathcal{F}_{\tau_{k}}}\left[g(\tau_{k+1},X_{\tau_{k+1}},\Theta_{\tau_{k+1}})\right] &= \mathbb{E}_{(\tau_{k},X_{\tau_{k}},\Theta_{\tau_{k}})}\left[g(\tau_{k+1},X_{\tau_{k+1}},\Theta_{\tau_{k+1}})\right]\, \\
        &= c(\tau_{k},X_{\tau_{k}},\Theta_{\tau_{k}})\,,
    \end{align*}
    for some non-negative Borel function $c(\cdot,\cdot,\cdot)$.
\end{lemma}
\begin{lemma}
    Suppose that Assumption~\ref{SA} applies, $X$ and $\Theta$ are independent,  and $\tau_{k+1}-\tau_{k}$ is independent of $\Theta_{\tau_{k}}$. Then we have that
    \[ \mathbb{E}_{\mathcal{F}_{\tau_{k}}}\left[g(\tau_{k+1},X_{\tau_{k+1}})\right] =c(\tau_{k},X_{\tau_{k}})\,, \]
    for some non-negative Borel function $c(\cdot,\cdot)$.
\end{lemma}

\begin{proof}
    By Assumption~\ref{SA},
    \begin{align*}
        \mathbb{E}_{\mathcal{F}_{\tau_{k}}}\left[g(\tau_{k+1},X_{\tau_{k+1}})\right] & =\mathbb{E}_{\left(\tau_{k},X_{\tau_{k}},\Theta_{\tau_{k}}\right)}\left[g(\tau_{k}+\tau_{k+1}-\tau_{k},X_{\tau_{k}+\tau_{k+1}-\tau_{k}})\right]\, \\
        & =\mathbb{E}_{\left(\tau_{k},X_{\tau_{k}}\right)}\left[g(\tau_{k+1},X_{\tau_{k+1}})\right] =: c(\tau_{k},X_{\tau_{k}})\,.
    \end{align*}
\end{proof}
\begin{example}
    Take $X$ independent of $\Theta$, where $\Theta$ is a continuous-time Markov chain with state space $\mathbb{N}_{0}$, generator
    \begin{equation*}
        \left[
            \begin{array}
                [c]{cccc}
                -\lambda_{1} & \lambda_{1}  & 0            & \ldots \\
                0            & -\lambda_{2} & \lambda_{2}  & \ldots \\
                0            & 0            & -\lambda_{3} & \ldots \\
                \ldots       & \ldots       & \ldots       & \ldots
            \end{array}
            \right]\,, \quad \lambda_{i}>0\,, \quad i=1,2,\ldots\,,
    \end{equation*}
    and $\Theta_{0}=0$. One thus has that $\tau_{k}=\inf\left\{s\geq0:\Theta_{s}=k\right\}$.
\end{example}
\begin{example}
    As a more general example, consider $\left(X,\Theta\right)$ to be given by a jump-diffusion, more specifically, a strong Markov process of the form
    \begin{align*}
        \mathrm{d}X_{t} &= \sigma(t,X_{t})X_{t}\mathrm{d}W\,,\quad X_{0}=x_{0}\,, \\
        \mathrm{d}\Theta_{t} &= \int_{\mathbb{R}^{q}}z\mathcal{N}\left(X_{t},\mathrm{d}t,\mathrm{d}z\right)\,, \quad \Theta_{0}=0\,,
    \end{align*}
    where $\mathcal{N}\left(x,\mathrm{d}t,\mathrm{d}z,\omega\right)$ is some Poisson random measure on $\mathbb{N}_{0}$, independent of $W$, with
    \[ \mathbb{P}\left(\mathcal{N}\left(x,(s,t],B\right)=k\right):=\exp(-v(x,B)\left( t-s\right))\frac{v(x,B)^{k}\left(t-s\right)^{k}}{k!},\quad k\in\mathbb{N}_{0}\,, \]
    where for simplicity $v(x,B)$ is of finite activity in the sense that
    \[ v(x,\{0\})=0 \quad \text{and} \quad \int v(x,\mathrm{d}z)<\infty\,, \quad x\in\mathbb{R}^{d}\,. \]
    The stopping times are determined by $\tau_{0}=0$ and $\tau_{k}=\inf\left\{s>0:\Theta_{s}=k\right\}$, $k > 0$.

    (i): Simple case with $d=q=1$ and $v(x,B)= \lambda\delta_{1}(B) := \lambda1_{B}(1)$:
    \[ \mathbb{P}\left(\mathcal{N}\left(x,(0,t],B\right)=k\right)=\exp\left(-\lambda t\delta_{1}(B)\right)\frac{\lambda^{k}t^{k}\delta_{1}(B)^{k}}{k!}\,, \]
    with $0^{0} := 1$, and hence we get
    \begin{gather*}
        \Theta_{t}=\int_{(0,t]}\int_{\mathbb{R}}z\mathcal{N}\left(\mathrm{d}s,\mathrm{d} z\right)=\mathcal{N}\left((0,t],\left\{1\right\}\right)\,, \quad \text{where}\\
        \mathbb{P}\left(\mathcal{N}\left((0,t],\left\{1\right\}  \right)=k\right)
        =\exp\left(-\lambda t\right)\frac{\lambda^{k}t^{k}}{k!}\,,
    \end{gather*}
    i.e., $\Theta_{t}=\sharp\,\text{jumps}$ in the interval $(0,t]$ observed by the agent. With $\tau_{0}=0$ and $\tau_{k}=\inf\left\{s>0:\Theta_{s}=k\right\}$, $k > 0$, one has that $\tau_{k}-\tau_{k-1}$, $k \geq 1$, is $\mathrm{Exp}(\lambda)$ distributed and independent of $\tau_{k-1}$. In particular, the conclusion of Lemma~\ref{lemSMP} applies.

    (ii): Suppose $d=q=1$ and that the Poisson random measure is only active during times that the process $X$ is visiting a certain interval $I^{\pm} := \left[c_{-}, c_{+}\right]$ with $0 < c_{-} < c_{+}$. That is, we consider
    \begin{align*}
        v(x,B) &:= \lambda1_{\left[c_{-},c_{+}\right] }(x)\delta_{1}(B)\,, \quad \text{hence} \\
        \mathbb{P}\left(\mathcal{N}\left(x,(0,t],B\right)=k\right) &= \exp\left( -\lambda t1_{\left[c_{-},c_{+}\right]}(x)\delta_{1}(B)\right)\frac{\lambda^{k}t^{k}1_{\left[c_{-},c_{+}\right]}(x)^{k}\delta_{1}(B)^{k}}{k!}
    \end{align*}
    with $0^{0} := 1$, and
    \begin{align*}
        \Theta_{t} & =\int_{(0,t]}1_{\left[c_{-},c_{+}\right]}(X_{s-})\int_{\mathbb{R}}z\mathcal{N}\left(X_{s-},\mathrm{d}s,\mathrm{d}z\right)\, \\
                   & =\int_{(0,t]}1_{\left[c_{-},c_{+}\right]}(X_{s})\mathcal{N}_{0}\left(\mathrm{d}s,\{1\}\right)\,,
    \end{align*}
    with $\mathcal{N}_{0}$ as in (i). Note that $X$ is continuous. 
    So $\Theta_{t}=\sharp\,\text{jumps}$ in the interval $(0,t]$ while $X$ is in the interval $\left[c_{-},c_{+}\right]$, which is observable by the agent. 
    Let $\tau_{0} = 0$, and $\tau_{k} = \inf\left\{s>0:\Theta_{s}=k\right\}$, $k>0$. 
    Then, obviously, $\Theta_{\tau_{k}}=k$, $k \geq 0$, and
    \begin{align*}
        \tau_{k+1}-\tau_{k} & \in\sigma\left\{\left(X_{s},\mathcal{N}_{0}\left(\left(\tau_{k},s\right],\{1\}\right)\right)1_{\left\{s\geq\tau_{k}\right\}}\right\}\, \\
                            & \subset\sigma\left\{\left(X_{s},\Theta_{s}\right)1_{\left\{s\geq\tau_{k}\right\}}\right\}\,.
    \end{align*}
    Hence, in particular Assumption~\ref{SA} applies, i.e.,  by the Markov property,
    \begin{align*}
         &\mathbb{E}_{\mathcal{F}_{\tau_{k}}}\left[g(\tau_{k+1},X_{\tau_{k+1}},\Theta_{\tau_{k+1}})\right]\, \\
         &= \mathbb{E}_{\mathcal{F}_{\tau_{k}}}\left[g(\tau_{k}+\tau_{k+1}-\tau_{k},X_{\tau_{k}+\tau_{k+1}-\tau_{k}},\Theta_{\tau_{k}+\tau_{k+1}-\tau_{k}})\right]\, \\
         &= \mathbb{E}_{\mathcal{F}_{\tau_{k}}}\left[g(\tau_{k}+\tau_{k+1}-\tau_{k},X_{\tau_{k}+\tau_{k+1}-\tau_{k}},\Theta_{\tau_{k}+\tau_{k+1}-\tau_{k}})\right]\, \\
         &= \mathbb{E}_{\left(\tau_{k},X_{\tau_{k}},\Theta_{\tau_{k}}\right)}\left[g(\tau_{k}+\tau_{k+1}-\tau_{k},X_{\tau_{k}+\tau_{k+1}-\tau_{k}},\Theta_{\tau_{k}+\tau_{k+1}-\tau_{k}})\right]\, \\
         &= \mathbb{E}_{\left(\tau_{k},X_{\tau_{k}},\Theta_{\tau_{k}}\right)}\left[g(\tau_{k+1},X_{\tau_{k+1}},\Theta_{\tau_{k+1}})\right]\, \\
         &=: c(\tau_{k},X_{\tau_{k}},\Theta_{\tau_{k}})\,.
    \end{align*}
    So, this is a simple example in which the continuation functions are dependent on the evaluation of both $X$ and $\Theta$.
\end{example}
\subsection{Convergence of random times least-squares Monte Carlo}\label{app: lsmc convergence}
In the setting of Assumption~\ref{SA}, consider measurable functions 
\begin{align}
    \psi_{l}:\mathbb{R}_{\geq0}\times\mathbb{R}^{d}\mathbb{\times R}^{q}\rightarrow\mathbb{R}\,,\quad l\in\mathbb{N}\,.
    \label{eq:basis}
\end{align}
To establish convergence results, we require the following assumptions on the set of basis functions that are used in our regression:

\begin{assumption}\label{CLPass}
    For the collection of measurable functions in \eqref{eq:basis}, let $P_{k}^{L}$ denote the $L^{2}$-projection of $\sigma\left\{\left(\tau_{k},X_{\tau_{k}},\Theta_{\tau_{k}}\right)\right\}$-measurable random variables onto the subspace
    \[ {\cal S}_L = \overline{\mathsf{span}\left\{\psi_{l}\left(\tau_{k},X_{\tau_{k}},\Theta_{\tau_{k}}\right),~l=1,\ldots,L\right\}}\,. \]
    We assume the following. For each $k=1,\ldots,K-1$:
    \begin{itemize}
        \item[1.] The sequence $\left(\psi_{l}\left(\tau_{k},X_{\tau_{k}},\Theta_{\tau_{k}}\right)\right)_{l\in\mathbb{N}}$ is total in $L^{2}\left(\Omega,\sigma\left\{\left(\tau_{k},X_{\tau_{k}},\Theta_{\tau_{k}}\right)\right\},\mathbb{P}\right)$, i.e.,\\ ${\cal S}_\infty =L^{2}\left(\Omega,\sigma\left\{\left(\tau_{k},X_{\tau_{k}},\Theta_{\tau_{k}}\right)\right\},\mathbb{P}\right)$.
        \item[2.] For $L$ $\in\mathbb{N}$, if $\sum_{l=1}^{L}\lambda_{l}\psi_{l}\left(\tau_{k},X_{\tau_{k}},\Theta_{\tau_{k}}\right) = 0$ a.s., then $\lambda_{l}=0$ for $l=1,\ldots L$.
        \item[3.] We have
              \[ \mathbb{P}\left(\sum_{l=1}^{L}\alpha_{k,l}\psi_{l}\left(\tau_{k},X_{\tau_{k}},\Theta_{\tau_{k}}\right)=Z\left(\tau_{k},X_{\tau_{k}}\right)\right) = 0\,, \]
              where
              \[ \sum_{l=1}^{L}\alpha_{k,l}\psi_{l}\left(\tau_{k},X_{\tau_{k}}, \Theta_{\tau_{k}}\right)=P_{k}^{L}\left(Z(\tau_{\mathfrak{n}_{k}^{L}}, X_{\tau_{\mathfrak{n}_{k}^{L}}})\right)\,, \]
              and $\left(\mathfrak{n}_{k}^{L}\right)_{k=0,\ldots,K}$ is a sequence of discrete stopping times with values in $\left\{1,\ldots,K\right\}$ solving the recursion $\mathfrak{n}_{K}^{L}=K$,
              \[ \mathfrak{n}_{k}^{L} = \begin{cases}
                      k, & Z\left(\tau_{k},X_{\tau_{k}}\right)\geq P_{k}^{L}\left(Z(\tau_{\mathfrak{n}_{k}^{L}},X_{\tau_{\mathfrak{n}_{k}^{L}}})\right) \\
                      \mathfrak{n}_{k+1}^{L}, & Z_{k}<P_{k}^{L}\left(Z(\tau_{\mathfrak{n}_{k}^{L}},X_{\tau_{\mathfrak{n}_{k}^{L}}})\right)
                  \end{cases} \,, \quad 1\leq k<K\,. \]
    \end{itemize}
\end{assumption}

Under Assumptions~\ref{SA} and~\ref{CLPass}, we have the following convergence result:

\begin{theorem}\label{thm:conv-values}
    (\cite[Theorem 3.2]{CLP02}) Let $\mathfrak{n}_{k}^{\ast}$ be an optimal stopping index after the (physical) times $\tau_{1},\ldots,\tau_{k-1}$ have passed. 
    It then holds a.s.\ that
    \[ \lim_{L\rightarrow\infty}\lim_{N\rightarrow\infty}\frac{1}{N}\sum_{n=1}^{N}Z(\widehat{\tau}_{\widehat{\mathfrak{n}}_{k}^{L,n}}^{(n)},X_{\widehat{\tau}_{\widehat{\mathfrak{n}}_{k}^{L,n}}^{(n)}}^{(n)}) = \mathbb{E}\left[Z(\tau_{\mathfrak{n}_{k}^{\ast}},X_{\tau_{\mathfrak{n}_{k}^{\ast}}})\right]\,. \]
\end{theorem}

\noindent \cite{CLP02} further implies that one may obtain a sequence $\widehat{C}_{k}^{(r)}(t,x,\theta),$ $r\in\mathbb{N}$, of approximate continuation functions from a sequence of LS-training procedures due to $L_{r}$ basis functions and $N_{r}$ training trajectories, such that
\begin{equation}
    \widehat{C}_{k}^{(r)}(\tau_{k},X_{{\tau}_{k}},\Theta_{{\tau}_{k}})\overset{P}{\longrightarrow}C_{k}^{\ast}(\tau_{k},X_{\tau_{k}},\Theta_{\tau_{k}})\,, \quad k=1,\ldots,K\,,\label{cinp}
\end{equation}
for $r\rightarrow\infty\,,$, where $C^*_k$ denotes the optimal continuation function at time $k$, and the process $(\tau_{k}, X_{\tau_{k}}, \Theta_{\tau_{k}})_{k=1,\ldots,K}$ is independent of the training procedures. 
Under the same assumptions, one next has the following convergence result on stopping times.

\begin{theorem}\label{thm:conv-stopping}
    For each $r\in\mathbb{N}$ there exists an optimal stopping time $\left(\mathfrak{n}^{\ast,r}\right)$ such that
    \[ \mathbb{P}(\mathfrak{n}^{(r)}\neq\mathfrak{n}^{\ast,r})\rightarrow0\text{ for }r\rightarrow\infty\,, \]
    where
    \[ \mathfrak{n}^{(r)}:=\min\left\{k \geq 0 : Z(\tau_{k}, X_{\tau_{k}}) \geq \widehat{C}_{k}^{(r)}(\tau_{k}, X_{{\tau}_{k}}, \Theta_{\tau_{k}})\right\}\,, \quad r\in\mathbb{N}\,. \]
    If in particular $\mathfrak{n}^{\ast}$ is unique, one has that $\mathfrak{n}^{(r)}\overset{P}{\longrightarrow}\mathfrak{n}^{\ast}$.
\end{theorem}

\begin{proof}
    The result follows, in principle, from the convergence \eqref{cinp} of the continuation values. 
    For a detailed proof in the context of robust optimal stopping, we refer to \cite[Section 6]{KLLSS18}.
\end{proof}
\subsection{Solution for stopped geometric Brownian motion with jumps}\label{app:benchmark}
\begin{proof}[Proof of Theorem~\ref{thm:bench-sub}]
    We begin by solving the problem for $K < \infty$ and will denote the value function at a time $t\in [0,T)$ by $v_K(t,x,d)$, where $x\in\mathbb{R}$ is the value of the state variable at time $t$, $d$ is a binary variable that indicates whether we may stop $(d = 1)$ or not $(d = 0)$ at this given moment and $K$ denotes the number of stopping opportunities remaining. 
    This implies that if ${\cal T}_{K,t}$ denotes the set of the first $K$ stopping times that occur strictly after  time $t$, we can write
    \begin{equation*}
    v_K(t,x,0) = \sup_{\tau\in {\cal T}_{K,t} }{\mathbb{E}}[Z_\tau|X_t=x],
    \qquad
    v_K(t,x,1) = \sup_{\tau\in ({\cal T}_{K,t}\cup\{t\}) }{\mathbb{E}}[Z_\tau|X_t=x].
    \end{equation*}
    
    The values that are characterized in the theorem, $Y_{i}^{(K)}$ and $C_i^{(K)}$, are given by $v_K(\tau_i, X_{\tau_i}, 1)$ and $v_K(\tau_i, {X_{\tau_i}}, 0)$, respectively.
    By definition, $v_0 \equiv 0$. 
    If we are allowed to stop at time $t$, i.e., $d = 1$, we either receive $x^\eta$ or continue with the problem with $K - 1$, so for $K \geq 1$,
    \begin{align}
        v_K(t,x,1) &= x^\eta1_{\{t\leq T\}}\vee v_{K-1}(t,x,0)\,. \label{eq:bench1}
        \intertext{If we are not allowed to stop at time $t$, the value function at that time equals the discounted version of the value at the following possible stopping time, if this stopping times occurs before maturity:}
        v_K(t,x,0) &= \int_0^{T-t}\lambda e^{-\lambda
        u}
        \mathbb{E}\left[v_K\left(t+u, x(1+j)e^{(\mu - \frac{1}{2}\sigma^2)u + \xi\sigma\sqrt{u}}, 1\right)\right]
        \dd{u}\,, \label{eq:bench2}
    \intertext{with $\xi \sim \mathcal{N}(0,1)$.\newline 
    \indent We now make an \emph{ansatz} to solve \eqref{eq:bench1} and \eqref{eq:bench2}:}
        v_K(t,x,0) &= x^\eta e^{\xi(T-t)}f_K(T-t)\,, \label{eq:bench-ansatz}
    \end{align}
    with $\xi$ as defined in \eqref{eq:bench-alpha}. 
    We substitute \eqref{eq:bench1} and \eqref{eq:bench-ansatz} in \eqref{eq:bench2} and see that it suffices to find a sequence of functions $(f_K)_{K\in\mathbb{N}_0}$ that satisfies
    \begin{equation}
    f_{K+1}(s) = \psi\int_0^s \left(e^{-\xi u}\vee f_K(u)\right)\,\mathrm{d}u\,,\label{eq: ansatzfint}
    \end{equation}
    where $\psi = \lambda(1+j)^\eta$ and $f_0 \equiv 0$.
    This immediately gives that, for all $s\geq 0$, we have $f_1(s)=\frac\psi \xi (1-e^{-\xi s})$.
    Let $s^*$ be the solution to $f_1(s)=e^{-\xi s}$ so $s^*=\ln(1+\xi/\psi )/\xi$ and define $\zeta:=f_1(s^*)=e^{-\xi  s^*}=\psi /(\psi +\xi)$.  
    We use a further ansatz
    \begin{equation*}
       f_{K+1}(s)=f_1(s) 1_{\{s\leq s^*\}}+\zeta\, z_{K+1}(s-s^*) 1_{\{s>s^*\}}
    \end{equation*}
    and conclude that this is a solution to \eqref{eq: ansatzfint} for $K>0$ if $$\partial_s {z}_{K+1}(s)= \partial_s {f}_{K+1}(s+s^*)=\psi  {f}_{K}(s+s^*)=\psi  z_K(s), \qquad z_{K+1}(0)=1,$$ while for $K=0$ we have 
    $$z_1(s)=\frac{1}\zeta f_1(s+s^*)=\frac\psi{\zeta\xi} (1-e^{-\xi s}e^{-\xi s^*})=\frac\psi {\zeta\xi} (1-\zeta e^{-\xi s})=1+ \frac\psi  \xi (1-e^{-\xi s}).$$ 
    A power series expansion gives
    $$z_{K+1}(s)=\sum_{m=0}^K \frac{(\psi s)^m}{m!} + (\frac{-\psi}\xi)^{K+1}\sum_{m=K+1}^\infty \frac{(-\xi s)^m}{m!}\ = \
    \frac{\Gamma_{K+1}(\psi s) }{K!}e^{\psi s}  +   \frac{ \gamma_{K+1}(-\xi s) }{K!}  (-\frac{\psi}\xi)^{K+1} e^{-\xi s},$$
    with $\Gamma_{K+1}$ and $\gamma_{K+1}$ the upper and lower Gamma functions of order $K+1$, respectively.

    Combining this with our ansatz \eqref{eq:bench-ansatz}, we arrive at an expression for the value function in terms of $t^*:=T-s^*$:
    \begin{equation}
      v_K(t, x, 0) = x^\eta e^{\xi(T-t)}\left(1_{\{t \geq t^*\}}f_1(T-t) + 1_{\{t < t^*\}}\zeta z_{K}(t^* - t)\right)\,. \label{eq:bench-sol-v-proof}
    \end{equation}
    
    For finite $K$, equations \eqref{eq:bench1} and \eqref{eq:bench-ansatz} show that it is optimal to stop if $1\geq e^{\xi(T-t)}f_{K-1}(T-t)$, which we can rewrite as $ e^{\xi s^*}f_{K-1}(s^*) \geq e^{\xi(T-t)}f_{K-1}(T-t)$ or $t\geq t^*=T-s^*$. For finite $K$, the expressions in the theorem for $C_i^{(K)}$ then follow by taking $t=\tau_i$ and $x=X_{\tau_i}$ in \eqref{eq:bench-ansatz} and  rewriting $\zeta e^{\xi(T-t)} = e^{\xi(t^* - t)}$.
    
    Taking a limit $K \to \infty$ in \eqref{eq:bench-sol-v-proof}, we find that
    \[\lim_{K\to\infty} v_K(t,x,0) = x^\eta\left(1_{\{t \geq t^*\}}e^{\xi(T-t)}f_1(T-t) + 1_{\{t < t^*\}} \zeta e^{\xi(T-t)} e^{\psi(t^*-t)}\right),\]
    and $\lim_{K\to\infty} v_K(\tau_i,X_{\tau_i},0) =  C_i^{(\infty)}$ follows from letting $\varepsilon \to 0$ in Lemma~\ref{lem:FHA}. 
    It is clear that it is optimal to stop if and only if $t \geq t^*$ for this case as well.
    This concludes the proof.
\end{proof}

\newpage\section{Additional Tables}\label{app:tables}
\begin{table}[h!]
    \centering

{\small\begin{tabular}{rr||rrrrrrrr}
$\lambda$ & $X_0$ & LSMC & (s.e.) & Dual & (s.e.) & Trivial & (s.e.) & Trivial-PI & (s.e.) \\\hline\hline 
1 & 90  & 3.07696 & (0.00578) & 3.07707 & (0.00006) & 2.75785 & (0.00516) & 3.06875 & (0.01232) \\
  & 100  & 6.04603 & (0.00792) & 6.04657 & (0.00022) & 5.32202 & (0.00695) & 6.06463 & (0.01883) \\
  & 110  & 10.63174 & (0.00986) & 10.63399 & (0.00065) & 10.02295 & (0.00903) & 10.61661 & (0.02016) \\
3 & 90  & 3.97787 & (0.00619) & 3.98025 & (0.00074) & 2.77611 & (0.00414) & 3.93131 & (0.01836) \\
  & 100  & 7.46066 & (0.00811) & 7.46832 & (0.00162) & 4.96859 & (0.00514) & 7.38342 & (0.02606) \\
  & 110  & 12.57331 & (0.00951) & 12.57662 & (0.00069) & 10.43446 & (0.00659) & 12.42274 & (0.02737) \\
5 & 90  & 4.17705 & (0.00621) & 4.18727 & (0.00236) & 2.47309 & (0.00343) & 4.09642 & (0.01923) \\
  & 100  & 7.73938 & (0.00807) & 7.75021 & (0.00290) & 4.34751 & (0.00416) & 7.54921 & (0.02649) \\
  & 110  & 12.93490 & (0.00932) & 12.95111 & (0.00460) & 10.19641 & (0.00546) & 12.58117 & (0.02676) \\
\end{tabular}
}
    \caption{Value of a max-call option on random times ($d=1$). \\ 
    {\footnotesize\textit{Notes.} 
    This table displays estimated values of $Y_0^{(\infty)}$ in the optimal stopping problem for payoff \eqref{eq:maxcall} with strike price $\mathcal{K} = 100$ and independent Poissonian arrivals of stopping opportunities defined in \eqref{eq:poisson}, with $d = 1$ independent geometric Brownian motion (without jumps, $\mu_J^j = \lambda_J^j = 0$) and with varying arrival rates $\lambda$ and initial asset prices $X_0$. 
    Column \rm{LSMC} contains the lower-biased estimates based on the random times LSMC algorithm of Section~\ref{sec: fha and ls} and column \rm{Dual} contains the upper-biased estimates \eqref{dual-mc} based on the dual algorithm in Section~\ref{sec: dual alg} and the estimated LSMC policy; standard errors for the dual estimates are the standard errors of the duality gap \eqref{gap-mc}. 
    Both algorithms use the same settings as in Table~\ref{tab:benchmark}, except for the basis functions: the continuation values in the LSMC algorithm are estimated by using the basis functions $(t, x) \mapsto t^ix^j$, $i + j \leq 5$. 
    Column \rm{Trivial} corresponds to stopping as soon as $Z_{\tau_i} > 0$ or $\tau_i \geq T$ and column \rm{Trivial-PI} contains the results of the policy iteration method \eqref{polit1} with $100\mathord{,}000$ outer paths, $500$ inner paths and a window length $\kappa = K'$. 
    All simulated paths have been generated in antithetic pairs.}}
    \label{tab: LSMC poisson 1}
\end{table}

\begin{table}[h!]
    \centering

{\small\begin{tabular}{rr||rrrrrrrr}
$\lambda$ & $X_0$ & LSMC & (s.e.) & Dual & (s.e.) & Trivial & (s.e.) & Trivial-PI & (s.e.) \\\hline\hline  
1 & 90  & 5.68421 & (0.00765) & 5.68804 & (0.00116) & 4.91301 & (0.00661) & 5.66093 & (0.01842) \\
  & 100  & 10.55799 & (0.01011) & 10.55963 & (0.00053) & 8.68818 & (0.00839) & 10.55060 & (0.02845) \\
  & 110  & 17.13788 & (0.01224) & 17.14160 & (0.00096) & 15.49008 & (0.01041) & 17.12895 & (0.03247) \\
3 & 90  & 7.33534 & (0.00818) & 7.34869 & (0.00221) & 4.58246 & (0.00498) & 7.22310 & (0.02632) \\
  & 100  & 12.99966 & (0.01048) & 13.02089 & (0.00395) & 7.24893 & (0.00580) & 12.77270 & (0.03689) \\
  & 110  & 20.33221 & (0.01235) & 20.35280 & (0.00243) & 14.78042 & (0.00707) & 19.95245 & (0.04102) \\
5 & 90  & 7.69256 & (0.00826) & 7.76684 & (0.03866) & 3.95536 & (0.00403) & 7.51463 & (0.02730) \\
  & 100  & 13.46328 & (0.01047) & 13.51186 & (0.00739) & 6.04586 & (0.00460) & 13.04278 & (0.03635) \\
  & 110  & 20.89358 & (0.01222) & 21.02209 & (0.06126) & 13.91635 & (0.00568) & 20.23306 & (0.03970) \\
\end{tabular}
}
    \caption{Value of a max-call option on random times ($d=2$). \\ 
    {\footnotesize\textit{Notes.} 
    This table displays estimated values of $Y_0^{(\infty)}$ in the optimal stopping problem for payoff \eqref{eq:maxcall} with strike price $\mathcal{K} = 100$ and independent Poissonian arrivals of stopping opportunities defined in \eqref{eq:poisson}, with $d = 2$ independent geometric Brownian motions (without jumps, $\mu_J^j = \lambda_J^j = 0$) and with varying arrival rates $\lambda$ and initial asset prices $X_0$. 
    Column \rm{LSMC} contains the lower-biased estimates based on the random times LSMC algorithm of Section~\ref{sec: fha and ls} and column \rm{Dual} contains the upper-biased estimates \eqref{dual-mc} based on the dual algorithm in Section~\ref{sec: dual alg} and the estimated LSMC policy; standard errors for the dual estimates are the standard errors of the duality gap \eqref{gap-mc}. 
    Both algorithms use the same settings as in Table~\ref{tab:benchmark}, except for the basis functions: as in \cite{AB04}, we use basis functions in the ordered asset prices $X_t^{(1)} \geq X_t^{(2)}$. 
    The continuation values in the LSMC algorithm are estimated by using the basis functions $(t, x^{(1)}, x^{(2)}) \mapsto t^i\left(x^{(1)}\right)^j\left(x^{(2)}\right)^k$, $i + j + k \leq 5$. 
    Column \rm{Trivial} corresponds to stopping as soon as $Z_{\tau_i} > 0$ or $\tau_i \geq T$ and column \rm{Trivial-PI} contains the results of the policy iteration method \eqref{polit1} with $100\mathord{,}000$ outer paths, $500$ inner paths and a window length $\kappa = K'$. 
    All simulated paths have been generated in antithetic pairs.}}
    \label{tab: LSMC poisson 2}
\end{table}

\begin{table}[h!]
    \centering

{\small\begin{tabular}{rr||rrrrrrrr}
$\lambda$ & $X_0$ & LSMC & (s.e.) & Dual & (s.e.) & Trivial & (s.e.) & Trivial-PI & (s.e.) \\\hline\hline 
1 & 90  & 3.37559 & (0.00626) & 3.37832 & (0.00080) & 3.03057 & (0.00559) & 3.37001 & (0.01349) \\
  & 100  & 6.42173 & (0.00846) & 6.42563 & (0.00128) & 5.64728 & (0.00739) & 6.40339 & (0.01983) \\
  & 110  & 10.98917 & (0.01043) & 10.99165 & (0.00061) & 10.30297 & (0.00947) & 11.00701 & (0.02226) \\
3 & 90  & 4.36383 & (0.00667) & 4.36687 & (0.00081) & 3.03411 & (0.00445) & 4.32993 & (0.01991) \\
  & 100  & 7.91895 & (0.00863) & 7.92280 & (0.00141) & 5.25481 & (0.00547) & 7.78849 & (0.02775) \\
  & 110  & 13.01151 & (0.01018) & 13.01976 & (0.00245) & 10.62043 & (0.00690) & 12.82926 & (0.02978) \\
5 & 90  & 4.57511 & (0.00678) & 4.58881 & (0.00237) & 2.70413 & (0.00370) & 4.47377 & (0.02088) \\
  & 100  & 8.20940 & (0.00865) & 8.22161 & (0.00203) & 4.60046 & (0.00443) & 7.97062 & (0.02823) \\
  & 110  & 13.37356 & (0.01002) & 13.39197 & (0.00342) & 10.32018 & (0.00570) & 13.01693 & (0.02926) \\
\end{tabular}
}

    \caption{Value of a max-call option on random times with jumps ($d=1$). \\ 
    {\footnotesize\textit{Notes.} 
    This table displays estimated values of $Y_0^{(\infty)}$ in the optimal stopping problem for payoff \eqref{eq:maxcall} with strike price $\mathcal{K} = 100$ and independent Poissonian arrivals of stopping opportunities defined in \eqref{eq:poisson} with $d = 1$ independent geometric Brownian motion, with jumps, $\mu_J^j = 0.06$ and $\lambda_J^j = 1$ and with varying arrival rates $\lambda$ and initial asset prices $X_0$. 
    The respective algorithms are carried out as described in the caption of Table~\ref{tab: LSMC poisson 1}. 
    }}
    \label{tab: LSMC jumps poisson 1}
\end{table}

\begin{table}[h!]
    \centering

{\small\begin{tabular}{rr||rrrrrrrr}
$\lambda$ & $X_0$ & LSMC & (s.e.) & Dual & (s.e.) & Trivial & (s.e.) & Trivial-PI & (s.e.) \\\hline\hline 
1 & 90  & 6.23507 & (0.00828) & 6.24121 & (0.00197) & 5.36950 & (0.00713) & 6.21446 & (0.02009) \\
  & 100  & 11.21449 & (0.01083) & 11.22446 & (0.00226) & 9.20470 & (0.00893) & 11.21218 & (0.03066) \\
  & 110  & 17.81153 & (0.01306) & 17.81885 & (0.00166) & 15.97752 & (0.01098) & 17.78208 & (0.03489) \\
3 & 90  & 8.02699 & (0.00881) & 8.04256 & (0.00331) & 4.98881 & (0.00536) & 7.92356 & (0.02882) \\
  & 100  & 13.80539 & (0.01116) & 13.82755 & (0.00459) & 7.66957 & (0.00618) & 13.58178 & (0.03947) \\
  & 110  & 21.18165 & (0.01310) & 21.21136 & (0.00554) & 15.10441 & (0.00747) & 20.82862 & (0.04389) \\
5 & 90  & 8.40994 & (0.00884) & 8.44444 & (0.00501) & 4.30081 & (0.00434) & 8.23108 & (0.02953) \\
  & 100  & 14.30547 & (0.01111) & 14.35577 & (0.00655) & 6.39494 & (0.00490) & 13.95940 & (0.03905) \\
  & 110  & 21.77564 & (0.01296) & 21.83676 & (0.00848) & 14.15769 & (0.00600) & 21.09714 & (0.04272) \\
\end{tabular}
}
    \caption{Value of a max-call option on random times with jumps ($d=2$). \\ 
    {\footnotesize\textit{Notes.} 
    This table displays estimated values of $Y_0^{(\infty)}$ in the optimal stopping problem for payoff \eqref{eq:maxcall} with strike price $\mathcal{K} = 100$ and independent Poissonian arrivals of stopping opportunities defined in \eqref{eq:poisson} with $d = 2$ independent geometric Brownian motions, with jumps, $\mu_J^j = 0.06$ and $\lambda_J^j = 1$ and with varying arrival rates $\lambda$ and initial asset prices $X_0$. 
    The respective algorithms are carried out as described in the caption of Table~\ref{tab: LSMC poisson 2}. 
    }}
    \label{tab: LSMC jumps poisson 2}
\end{table}
\end{document}